\RequirePackage{amsthm}%
\documentclass[sn-mathphys,Numbered]{sn-jnl}

\usepackage{graphicx}%
\usepackage{multirow}%
\usepackage{amssymb,amsfonts}%
\usepackage{mathrsfs}%
\usepackage[title]{appendix}%
\usepackage{xcolor}%
\usepackage{textcomp}%
\usepackage{manyfoot}%
\usepackage{booktabs}%
\usepackage{algorithm}%
\usepackage{algorithmicx}%
\usepackage{algpseudocode}%
\usepackage{listings}%
\usepackage{amsmath}%
\usepackage{upref}%
\usepackage{dsfont}%
\usepackage{stmaryrd}%
\usepackage{enumitem}%

\usepackage{amsmath,amssymb,graphicx,tikz}
\usetikzlibrary{arrows.meta}

  \let\<\langle \let\>\rangle \newcommand{\doi}[1]{doi:\,\href{http://dx.doi.org/#1}{\nolinkurl{#1}}}  \let\uml\"

\theoremstyle{thmstyleone}%
\newtheorem{theorem}{Theorem}[section]%
\newtheorem{proposition}[theorem]{Proposition}%
\newtheorem{lemma}[theorem]{Lemma}
\newtheorem{corollary}[theorem]{Corollary}

\theoremstyle{thmstyletwo}%
\newtheorem{example}[theorem]{Example}%
\newtheorem{remark}[theorem]{Remark}%

\theoremstyle{thmstylethree}%
\newtheorem{definition}[theorem]{Definition}%

\numberwithin{equation}{section}
\numberwithin{table}{section}
\numberwithin{figure}{section}

\allowdisplaybreaks
\begin{document}

\hypersetup{colorlinks=true, linkcolor=blue, citecolor=blue, urlcolor=blue}

\title[The random KKM mapping principle and its applications]{The random KKM mapping principle and its applications}

\author [2]{{\fnm{{Tiexin}} \sur{{Guo}}}}\email{{tiexinguo@csu.edu.cn}}

\author [1]{{\fnm{{Xiaohuan}} \sur{{Mu}}}}\email{{xiaohuanmu@163.com}}

\author *[1]{{\fnm{{Qiang}} \sur{{Tu}}}}\email{{qiangtu126@126.com}}

\makeatletter
{\let\@makefnmark\@empty
\footnotetext{This work was supported by the National Natural Science Foundation of China (Grant Nos.12371141).}}
\makeatother

\affil[1]{School of Mathematics and Statistics, Lingnan Normal University, Zhanjiang 524048, Guangdong Province, P. R. China}

\affil[2]{School of Mathematics and Statistics, Central South University, Changsha 410083, Hunan Province, P. R. China}

\abstract{Based on the recently developed random Brouwer fixed point theorem and the theory of stable compactness, we establish a random Knaster--Kuratowski--Mazurkiewicz (KKM) mapping principle in Hausdorff topological $L^0$-modules. As consequences of this principle, we obtain several results on fixed points, minimax problems, variational inequalities, and conditional games. These results solve a large class of noncompact problems in random functional analysis and its applications.}

\keywords{random KKM mapping principle, stable compactness, topological $L^0$-module, fixed point theorem, conditional Nash equilibrium}

\pacs[MSC Classification]{41A65, 46H25, 47H10, 46A22}

\maketitle

\tableofcontents

	\section{Introduction}\label{sec.1}

	In 1929, Knaster, Kuratowski and Mazurkiewicz \cite{KKM1929} established the Knaster--Kuratowski--Mazurkiewicz (briefly, KKM) mapping principle for simplices in Euclidean spaces (usually called the KKM lemma) by means of Sperner's lemma and applied it to give a new proof of the Brouwer fixed point theorem (the two results are equivalent; see \cite{DG1978} for a direct proof of the converse implication). In 1961, Ky Fan \cite{Fan1961} extended this principle to arbitrary subsets of Hausdorff topological vector spaces. By making use of this principle, Fan obtained a generalization of the Tychonoff fixed point theorem \cite{Fan1961}, established his minimax inequality \cite{Fan1972}, and developed coincidence and intersection theorems \cite{Fan1966}, which yield new proofs of Sion's minimax theorem \cite{Sion1958} and Nash's equilibrium theorem \cite{Nash1951,Nash1950}. The wide range of applications of the KKM mapping principle has established it as a powerful tool in nonlinear functional analysis; see \cite{Browder1967,Browder1968,DG1978,Fan1969,Ha1980,Lassonde1983,Granas2015} and the references therein.
	
    \par 
    In recent years, considerable progress has been made in fixed point theory in random functional analysis \cite{GWXYC2025,GZWG2020,MTG2025,TMG2024,TMGC2026}. Nevertheless, some basic problems remain open, including those posed by Guo et al. in \cite[Problem~6.1 and Remark~5.7]{GWXYC2025} that concern the random Tychonoff fixed point theorem and the conditional Nash equilibrium theorem for $L^0$-quasiconcave payoff functions in random locally convex modules. As was already called for in \cite[Remark~5.7]{GWXYC2025}, the establishment of a random KKM mapping principle is a fundamental and pressing problem in random functional analysis since it would enable us to follow Ky Fan's method in the study of fixed points, minimax problems, variational inequalities and conditional equilibria in the random setting. The main difficulty lies in overcoming noncompactness.
   
    \par
    To explain the background and the difficulties involved, we divide this section into three subsections. Throughout this section, $(\Omega,\mathcal{F},P)$ is a fixed probability space, $\mathbb{K}$ denotes the scalar field $\mathbb{R}$ of real numbers or $\mathbb{C}$ of complex numbers, and $\mathbb{N}$ denotes the set of positive integers. We denote by $L^{0}(\mathcal{F},\mathbb{K})$ (resp., $L^{0}(\mathcal{F},\mathbb{N})$) the set of equivalence classes of $\mathbb{K}$-valued (resp., $\mathbb{N}$-valued) random variables on $(\Omega,\mathcal{F},P)$. It is well known that $L^0(\mathcal{F},\mathbb{R})$ is a Dedekind complete lattice under the partial order $\leq$: $\xi\leq\eta$ if $\xi^0(\omega)\leq\eta^0(\omega)$ for almost all $\omega\in\Omega$, where $\xi^0$ and $\eta^0$ are arbitrarily chosen representatives of $\xi$ and $\eta$, respectively (similarly, we write $\eta<\xi$ on $1$ if $\eta^0(\omega)<\xi^0(\omega)$ for almost all $\omega\in\Omega$).

	\subsection{The challenge due to noncompactness in random functional analysis}\label{sec.1.1}

    The development of random functional analysis is based on the idea of randomizing traditional space theory of functional analysis. The basic framework of random functional analysis consists of random normed modules and random locally convex modules (see Definition~\ref{defn. RN module} for details), which are random generalizations of normed spaces and locally convex spaces, respectively, and were introduced by Guo \cite{Guo1992,Guo1993} in connection with the theory of probabilistic metric spaces \cite{SS2005}. Random normed modules (a special case of random locally convex modules) were also introduced independently by Gigli \cite{Gigli2018,GigliRCD2018} in connection with nonsmooth differential geometry on metric measure spaces (see \cite{BPS2023,CLPV2025,CGP2025,GLP2025,LP2019,LPV2024} for recent related advances). 
    
    \par 
    Two topologies are commonly used on a random locally convex module $E$, namely, the $(\varepsilon,\lambda)$-topology and the locally $L^0$-convex topology. In general, neither topology makes $E$ a locally convex space. Before 2009, the theory of random locally convex modules was developed under the $(\varepsilon,\lambda)$-topology, which is a linear topology and is metrizable in the case of a random normed module (see Proposition~\ref{prop.topology1}). Since this topology need not be locally convex, the classical duality theory for locally convex spaces does not apply directly. This led Guo to introduce the notion of a random conjugate space and to establish the random Hahn--Banach theorem \cite{Guo1992,Guo1993}. Subsequently, the theory of random conjugate spaces became one of the main subjects of research in random functional analysis \cite{Guo2008,GL2005,GXC2009}. An important advantage of the $(\varepsilon,\lambda)$-topology is that it allows one to use the normed spaces associated with random normed modules to connect random functional analysis with classical functional analysis, as shown in \cite{Guo2008,GL2005}. In 2009, motivated by financial applications, Filipovi\'c et al.~\cite{FKV2009} introduced the locally $L^0$-convex topology and used it to develop parts of random convex analysis. This topology is finer than the $(\varepsilon,\lambda)$-topology and is closer in form to the classical locally convex topology (see Proposition~\ref{prop.topology2}). It is, however, not in general a linear topology, so the connections with classical functional analysis established in \cite{Guo2008,GL2005} are difficult to obtain in the same way. To relate the results obtained under these two topologies, Guo introduced in \cite{Guo2010} the notion of $\sigma$-stability for nonempty subsets of an $L^0(\mathcal F,\mathbb K)$-module. This notion makes it possible to take full advantage of both topologies in the study of the same problem. With this approach, considerable progress has been made over the past decade in random convex analysis \cite{GZWYZ2017,GZWW2021} and especially in fixed point theory in random functional analysis \cite{GWXYC2025,GZWG2020,MTG2025,TMG2024,TMGC2026}.
   
	\par
	Compactness plays a fundamental role in the proofs of many theorems in functional analysis. In a Euclidean space, every bounded closed subset is compact. In the random setting, however, Guo's result \cite[Theorem~2.2]{Guo2008} implies that, unless $\mathcal{F}$ is essentially purely $P$-atomic, an a.s. bounded closed subset of the random Euclidean space $L^{0}(\mathcal{F},\mathbb{R}^{n})$ (see Subsection~\ref{sec.2.2} for its definition) need not be compact under either the $(\varepsilon,\lambda)$-topology or the finer locally $L^{0}$-convex topology. Thus the compactness that underlies classical results such as the Brouwer fixed point theorem and the KKM mapping principle cannot in general be taken for granted in the random setting. Similar difficulties arise when compactness is considered with respect to the random weak and random weak-star topologies, which are introduced by means of the theory of random conjugate spaces. As Guo showed in \cite[Theorems~2.3 and 2.4]{Guo2008}, when compactness is understood in the classical sense, the Banach--Alaoglu theorem for the random weak-star topology and the Banach--Bourbaki--Kakutani--\v{S}mulian theorem for the random weak topology both fail in general.

	\subsection{Approaches to overcoming noncompactness: the theory of stable compactness}\label{sec.1.2}

	The randomized Bolzano--Weierstrass theorem \cite{KS2001} provides a way to overcome the limitations of classical compactness in the random Euclidean space. It states that for every a.s. bounded sequence $\{x_n,n\in\mathbb{N}\}$ in $L^0(\mathcal{F},\mathbb{R}^d)$, there exists a sequence $\{n_k,k\in\mathbb{N}\}$ in $L^0(\mathcal{F},\mathbb{N})$ satisfying $n_k<n_{k+1}$ on $1$ for every $k\in\mathbb{N}$ such that $\{x_{n_k},k\in\mathbb{N}\}$ converges a.s. to some $x\in L^0(\mathcal{F},\mathbb{R}^d)$, where $x_{n_k}=\sum_{l=1}^{\infty}I_{(n_k=l)}x_l$ and $I_{(n_k=l)}$ denotes the equivalence class of the characteristic function of $(n_k^0=l):=\{\omega\in\Omega:n_k^0(\omega)=l\}$ for a representative $n_k^0$ of $n_k$. Kabanov and Stricker used this theorem to give a considerably simpler proof of the no-arbitrage criteria in mathematical finance \cite{KS2001}. By combining it with a random Sperner lemma, we established a general random Brouwer fixed point theorem in \cite{TMGC2026} to provide a tool for the study of stochastic differential equations \cite{Ponosov1988,Ponosov2022}. Motivated by the randomized Bolzano--Weierstrass theorem, Guo et al.~\cite{GWXYC2025} systematically developed the theory of random sequential compactness and established the noncompact Schauder fixed point theorem and the conditional Nash equilibrium theorem in random normed modules. The resulting theory provides an effective approach to noncompactness problems in random normed modules and has led to further developments in topological fixed point theory \cite{TMG2024,TMGYS2025}.
	
	\par 
	Another approach is provided by $L^0$-convex compactness. Weak compactness is of fundamental importance in infinite-dimensional convex optimization and variational analysis \cite{ET1999}, but need not retain its classical role in non-locally-convex spaces. The notion of convex compactness introduced by \v{Z}itkovi\'c provides a substitute for weak compactness in this setting and has found applications in convex optimization and mathematical economics \cite{Zitkovic2010}.  Motivated by this work, Guo et al.~\cite{GZWW2021} introduced $L^0$-convex compactness for closed $L^0$-convex subsets of random locally convex modules. By means of the theory of random conjugate spaces, they established a James-type characterization of $L^0$-convex compactness in random normed modules and applied it to random convex optimization and variational inequalities. This characterization was subsequently extended to random locally convex modules by Wu et al. \cite{WZZ2022}. Thus $L^0$-convex compactness serves as a substitute for weak compactness in the study of closed $L^0$-convex subsets of random locally convex modules. This theory has led to further developments in metric fixed point theory in random functional analysis, including the random Browder--Kirk fixed point theorem \cite{GZWG2020} and its extension to commuting families of nonexpansive mappings \cite{MTG2025}; see also \cite{SGT2025,GMTS2026} for related developments.
	
    \par
    Random sequential compactness and $L^0$-convex compactness have greatly advanced the study of fixed point theory in random functional analysis, but both have limitations. On the one hand, random sequential compactness is a random generalization of classical sequential compactness, and its theory is developed for subsets of random normed modules under the metrizable $(\varepsilon,\lambda)$-topology. In a general random locally convex module, however, neither the $(\varepsilon,\lambda)$-topology nor the locally $L^0$-convex topology is necessarily metrizable, so that random sequential compactness no longer suffices for the study of noncompactness problems in this more general setting. On the other hand, the intersection property defining $L^0$-convex compactness applies only to families of closed $L^0$-convex subsets, whereas the closed values of a KKM mapping need not be convex. Thus the present paper requires a corresponding generalization of ordinary compactness.
    
	\par 
    Both random sequential compactness and $L^0$-convex compactness rely on the theory of $\sigma$-stable sets introduced in \cite{Guo2010}. To provide a suitable framework for the further study of compactness, Guo et al.~\cite{GMT2024} generalized the notion of a $\sigma$-stable set to that of a $B$-stable set associated with a complete Boolean algebra $B$ (see \cite[Definition~5.3]{GMT2024}). Based on this work, Guo et al.~\cite{GMPTY2026} recently developed the theory of stable topological spaces and stably compact spaces. The hereditary complement introduced in \cite[Definition~2.15]{GMPTY2026} plays a central role in establishing the characterizations of stable compactness \cite[Theorem~2.33]{GMPTY2026}, since it overcomes the difficulty that the ordinary complement of a stable subset need not be stable. Within this framework, random sequential compactness of a $\sigma$-stable random metric space is equivalent to its stable compactness \cite[Theorem~2.50]{GMPTY2026}, while $L^0$-convex compactness of an $(\varepsilon,\lambda)$-closed $L^0$-convex subset of an $(\varepsilon,\lambda)$-complete random locally convex module is equivalent to its stable compactness under the stable random weak topology \cite[Corollary~2.66]{GMPTY2026}. The theory also yields the random Banach--Alaoglu, random James and random Banach--Bourbaki--Kakutani--\v{S}mulian theorems in terms of stable compactness. The work in \cite{GMPTY2026} provides the foundation for the present paper.

	\subsection{The main contributions of this paper}\label{sec.1.3}

	With the theory of stable compactness and the random Brouwer fixed point theorem \cite{TMGC2026} in hand, we can now establish the random KKM mapping principle in Hausdorff topological $L^0$-modules (random generalizations of topological vector spaces, with random locally convex modules as a special case; see Definition~\ref{defn.topological L0-module} for details). Following the classical approach, we first need to establish a random KKM lemma on the standard $L^0$-simplex in a random Euclidean space $L^{0}(\mathcal{F},\mathbb{R}^{n})$ by making use of the random Brouwer fixed point theorem. The first difficulty is to find an appropriate formulation of this lemma. Let us first recall the classical KKM lemma.
	\begin{proposition}[\cite{KKM1929}]\label{prop.KKM}
	Let $\Delta_n=\operatorname{conv}(\{e_1,\ldots,e_n\})$ be the standard $(n-1)$-dimensional simplex in $\mathbb R^n$ and $F_1,\ldots,F_n$ be closed subsets of $\Delta_n$ such that
	\begin{equation}\label{eq.covering-condition1}
	\operatorname{conv}(\{e_i\mid i\in I\})\subseteq \bigcup\{F_i\mid i\in I\}
	\end{equation}
	for each nonempty subset $I$ of $\{1,\ldots,n\}$. Then $\bigcap_{i=1}^n F_i\neq\emptyset$.
	\end{proposition}
	The covering condition \eqref{eq.covering-condition1} is the essential part of this lemma. In the random setting, we consider $\sigma$-stable closed subsets of the standard $L^0$-simplex and replace the ordinary convex hull by the $L^0$-convex hull. As shown in \cite{GWXYC2025}, $\operatorname{Conv}_{L^0}(\{e_i\mid i\in I\})$ is $\sigma$-stable. However, a finite union of $\sigma$-stable sets need not be $\sigma$-stable. To formulate a covering condition compatible with $\sigma$-stability, we adopt the following random version of \eqref{eq.covering-condition1}:
	\begin{equation}\label{eq.covering-condition2}
	\operatorname{Conv}_{L^0}(\{e_i\mid i\in I\})\subseteq \bigcup\sigma\{F_i\mid i\in I\},	
	\end{equation}
	where $\sigma\{F_i\mid i\in I\}$ denotes the $\sigma$-stable hull of the family $\{F_i\mid i\in I\}$ (it is easy to check that $\bigcup\sigma\{F_i\mid i\in I\}$ is $\sigma$-stable; see Definition~\ref{defn.stable family} for details). Since $\bigcup\{F_i\mid i\in I\}\subseteq \bigcup\sigma\{F_i\mid i\in I\}$  and the family $\sigma\{F_i\mid i\in I\}$ is not finite or even countable in general, the contradiction step in the classical proof \cite{DG1978} based on Brouwer fixed point theorem cannot be carried over directly. To overcome this difficulty, we use the barycentric coordinates of the fixed point obtained from the random Brouwer fixed point theorem to construct a countable partition of unity in $B_{\mathcal F}$ and then apply \eqref{eq.covering-condition2} to obtain the required contradiction. This gives the random KKM lemma stated in Theorem~\ref{thm.KKM1}. Motivated by \eqref{eq.covering-condition2}, we introduce a definition of a random KKM mapping (see Definition \ref{defn.KKM mapping}), which requires both the corresponding random covering property and $\sigma$-stability of the set-valued mapping. Following Ky Fan's approach and making full use of the properties of stably compact sets, we then generalize the random KKM lemma on the standard $L^0$-simplex (Theorem~\ref{thm.KKM1}) to the random KKM mapping principle in Hausdorff topological $L^0$-modules (Theorem~\ref{thm.KKM2}).
 
	\par
	We next apply Theorem~\ref{thm.KKM2} to prove a random Fan--Browder fixed point theorem for Hausdorff topological $L^0$-modules (Theorem~\ref{thm.Fan fixed point}), which provides the basis for many subsequent results. For this, we first need to construct a $\sigma$-stable set-valued mapping $T$ for a given mapping $\varphi:G\rightarrow 2^G\setminus\{\emptyset\}$. The classical construction
	$$T(y)=G\setminus\varphi^{-1}(y),\quad \forall y\in\varphi(G),$$
	does not guarantee the $\sigma$-stability of $T$ since the complement of a $\sigma$-stable set need not be $\sigma$-stable. Building on the theory of hereditary complements developed in \cite{GMPTY2026}, we establish Lemmas~\ref{lemm.fundament1}, \ref{lemm.stable1} and~\ref{lemm.stable2} and use them to construct a $\sigma$-stable mapping $T$ to which Theorem~\ref{thm.KKM2} can be applied. The resulting argument yields, for each $a\in B_{\mathcal{F}}$ with $a>0$, a partial fixed point $x_b\in G$ on some $b\in B_{\mathcal{F}}$ with $0<b\leq a$, namely, $I_bx_b\in I_b\varphi(x_b)$. To pass from these partial fixed points to a fixed point of $\varphi$, we establish an exhaustion-type lemma (Lemma~\ref{lemm.fixed point}).

	\par
	We next study several applications of Theorems~\ref{thm.KKM2} and~\ref{thm.Fan fixed point} in Hausdorff topological $L^0$-modules. The random coincidence theorem (Theorem~\ref{thm.coincidence}) and the random intersection theorem (Theorem~\ref{thm.intersection}) both follow directly from Theorem~\ref{thm.Fan fixed point}. To derive the conditional Nash equilibrium theorem (Theorem~\ref{thm.Nash equilibrium}) from the random intersection theorem, we first establish a version of the tube lemma for the stable topology (Lemma~\ref{lemm.tubular neighborhood}). To prove the random minimax inequality (Theorem~\ref{thm.Minimax inequality}) and the random minimax principle (Theorem~\ref{thm.Minimax principle}), we first establish a Weierstrass-type theorem on stably compact subsets of Hausdorff topological $L^0$-modules (Lemma~\ref{lemm.tools}). We then prove the former by Theorem~\ref{thm.KKM2} and the latter by Theorem~\ref{thm.Fan fixed point}.
	
	\par
	In random locally convex modules, different choices of topology lead to three notions of continuity (see Definition~\ref{defn.continuity}). For $\sigma$-stable mappings between $\sigma$-stable subsets of such modules, Lemma~\ref{lemm.two continuty} shows that each of the other two continuity conditions implies $\mathcal T_c(\mathcal P_{cc})$-continuity. We therefore state the random Tychonoff fixed point theorem and its extensions under this continuity assumption. By Theorem~\ref{thm.KKM2}, we prove the random Tychonoff fixed point theorem (Theorem~\ref{thm.tychonoff1}) and extend it to mappings satisfying an $L^0$-segment condition (Theorem~\ref{thm.inward1}). This extension is a random generalization of Ky Fan's fixed point theorem under a segment condition and yields a fixed point theorem for inward and outward mappings (Theorem~\ref{thm.inward2}). By Theorem~\ref{thm.Fan fixed point}, we also prove the random Schauder--Tychonoff theorem (Theorem~\ref{thm.Schauder-Tychonoff}) and the random Fan--Iokhvidov theorem (Theorem~\ref{thm.Fan-Iokhvidov}). For $\mathcal T_{\varepsilon,\lambda}$-complete real random inner product modules, we prove a random Hartman--Stampacchia theorem (Theorem~\ref{thm.Hartman-Stampacchia}) by making use of Theorem~\ref{thm.KKM2} and the properties of $L^0$-convex compactness.

	\par
	The remainder of this paper is organized as follows. Section~\ref{sec.2} recalls the necessary preliminaries and states the main results. Section~\ref{sec.proofs} gives the proofs in the order of Subsections~\ref{sec.2.2}--\ref{sec.2.7}. In Section~\ref{sec.4}, we show that, except for the random minimax principle (Theorem~\ref{thm.Minimax principle}) and the random Hartman--Stampacchia theorem (Theorem~\ref{thm.Hartman-Stampacchia}), all the main existence theorems proved in this paper are equivalent to the random Brouwer fixed point theorem.

    \section{Preliminaries and main results}\label{sec.2}

    \subsection{Preliminaries}\label{sec.preliminaries}

    \par 
    Throughout this paper, $\mathbb{K}$ denotes the scalar field $\mathbb{R}$ of real numbers or $\mathbb{C}$ of complex numbers, and $\mathbb{N}$ denotes the set of positive integers. Let $(\Omega,\mathcal{F},P)$ be a given probability space, and let $L^{0}(\mathcal{F},\mathbb{K})$ denote the algebra of equivalence classes of $\mathbb{K}$-valued random variables on $(\Omega,\mathcal{F},P)$. We write $L^{0}(\mathcal{F}):=L^{0}(\mathcal{F},\mathbb{R})$ and denote by $\bar{L}^{0}(\mathcal{F})$ the set of equivalence classes of extended real-valued random variables on $(\Omega,\mathcal{F},P)$. More generally, for any  separable metric space $M$, $L^{0}(\mathcal{F},M)$ denotes the set of equivalence classes of $M$-valued $\mathcal{F}$-measurable functions on $(\Omega,\mathcal{F},P)$.
    
    \par 
    It is well known from \cite{DS1958} that $\bar{L}^{0}(\mathcal{F})$ is a complete lattice under the partial order $\xi \leq \eta$ iff $\xi^{0}(\omega) \leq \eta^{0}(\omega)$ for almost all $\omega \in \Omega$ (briefly, $\xi^{0}(\omega) \leq \eta^{0}(\omega)$ a.s.), where $\xi^{0}$ and $\eta^{0}$ are arbitrarily chosen representatives of $\xi$ and $\eta$ in $\bar{L}^{0}(\mathcal{F})$, respectively. In particular, the sublattice $(L^{0}(\mathcal{F}),\leq)$ is a Dedekind complete lattice. For any nonempty subset $H$ of a complete lattice, we denote the supremum and infimum of $H$ by $\bigvee H$ and $\bigwedge H$, respectively.
    
    \par 
    For any $A\in\mathcal{F}$, we always use the corresponding lowercase letter $a$ for its equivalence class $[A]$, where $A$ and $B$ are equivalent if $P[(A\setminus B)\cup(B\setminus A)]=0$. The measure algebra associated with $(\Omega,\mathcal{F},P)$ is the complete Boolean algebra $B_{\mathcal{F}}=\{[A]\mid A\in\mathcal{F}\}$, with join $\vee$, meet $\wedge$, and complement $^{c}$ given by
    $$a\vee b:= [A\cup B],\quad a\wedge b:= [A\cap B],\quad a^{c}:=[\Omega\setminus A],$$
    and with least element $0:=[\emptyset]$ and greatest element $1:=[\Omega]$; see \cite{Kop1989} for details. The order on $B_{\mathcal{F}}$ is given by 
    $$a\le b~\text{if~and~only~if}~a\wedge b=a.$$
    As usual, $a<b$ means that $a\le b$ and $a\ne b$. For $a=[A]\in B_{\mathcal F}$, we write $P(a):=P(A)$.
 
    \par 
    For $\xi,\eta\in\bar{L}^{0}(\mathcal{F})$, $\xi< \eta$  means $\xi\leq \eta$ and $\xi\neq \eta$, whereas, for any $a\in B_{\mathcal{F}}$, $\xi< \eta$ on $a$ means $\xi^{0}(\omega)< \eta^{0}(\omega)$  for almost all $\omega\in A$, where $A$, $\xi^{0}$ and $\eta^{0}$ are arbitrarily chosen representatives of $a$, $\xi$ and $\eta$, respectively. Moreover, we denote $L_{+}^{0}(\mathcal{F}):=\{\xi\in L^{0}(\mathcal{F})\mid \xi\geq 0\}$ and $L_{++}^{0}(\mathcal{F}):=\{\xi\in L^{0}(\mathcal{F})\mid \xi>0~\text{on}~1\}$.

    \subsubsection{Random locally convex modules and topological $L^0$-modules}\label{sec.2.1.1}

    \begin{definition}[\cite{Guo1992,Guo1993}]\label{defn. RN module}
    An ordered pair $(E, \|\cdot\|)$ is called a random normed module (briefly, an $RN$ module) over $\mathbb{K}$ with base $(\Omega, \mathcal{F}, P)$ if $E$ is a left  module over the algebra $L^{0}(\mathcal{F},\mathbb{K})$ (briefly, an $L^{0}(\mathcal{F},\mathbb{K})$-module) and $\|\cdot\|$ is a mapping from $E$ to  $L^{0}_{+}(\mathcal{F})$ such that the following conditions are satisfied:
    	\begin{itemize} [widest={($RN$-1)}, leftmargin=*]
    		\item [($RN$-1)] $\|\xi x\|= |\xi| \|x\|$ for any $\xi \in L^{0}(\mathcal{F}, \mathbb{K})$ and any $x\in E$;
    		\item [($RN$-2)] $\|x+y\|\leq \|x\|+ \|y\|$ for any $x$ and $y$ in $E$;
    		\item [($RN$-3)] $\|x\|=0$ implies $x=\theta$ (the null element of $E$). 
    	\end{itemize}
    As usual, a mapping $p: E\rightarrow L^{0}_{+}(\mathcal{F})$ satisfying ($RN$-1) and ($RN$-2) is called an $L^{0}$-seminorm on $E$; if it further satisfies ($RN$-3), then $p$ is called an $L^{0}$-norm on $E$. More generally, an ordered pair $(E,\mathcal{P})$ is called a random locally convex module (briefly, an $RLC$ module) over $\mathbb{K}$ with base $(\Omega,\mathcal{F},P)$ if $E$ is an $L^{0}(\mathcal{F},\mathbb{K})$-module and $\mathcal{P}$ a family of $L^{0}$-seminorms on $E$ such that
    \begin{itemize}[widest={($RLC$-1)}, leftmargin=*]
    	\item [($RLC$-1)] $\bigvee\{\|x\|\mid \|\cdot\|\in\mathcal{P}\}=0$ implies $x=\theta$.
    \end{itemize}
    \end{definition}

    \par 
    If $\mathcal{P}$ consists of a single $L^{0}$-seminorm $\|\cdot\|$ on $E$, the condition ($RLC$-1) is precisely ($RN$-3), so that $\|\cdot\|$ is an $L^{0}$-norm on $E$ and $(E,\mathcal{P})$ reduces to the $RN$ module $(E,\|\cdot\|)$. A canonical example of an $RN$ module over $\mathbb{K}$ with base $(\Omega,\mathcal{F},P)$ is $(L^{0}(\mathcal{F},\mathbb{K}),|\cdot|)$, where  $|\cdot|$  denotes the absolute value mapping on $L^{0}(\mathcal{F},\mathbb{K})$.
    
    \par
    For an $RLC$ module $(E,\mathcal{P})$, let $\mathcal{P}_{f}$ denote the set of nonempty finite subfamilies of $\mathcal{P}$, and for each $Q\in\mathcal{P}_{f}$ define the $L^{0}$-seminorm $\|\cdot\|_{Q}$ by 
    $$\|x\|_{Q}=\bigvee\{\|x\|\mid \|\cdot\|\in Q\}, \forall x\in E.$$
    The family $\mathcal{P}$ induces two natural topologies on $E$.

    \begin{proposition} [\cite{Guo1999,Guo2010}]\label{prop.topology1}
    Let $(E,\mathcal{P})$ be an $RLC$ module over $\mathbb{K}$ with base $(\Omega,\mathcal{F},P)$. For any given positive numbers $\varepsilon$ and $\lambda$ with $0<\lambda<1$, and any $Q\in \mathcal{P}_{f}$, let $U(Q,\varepsilon,\lambda)=\{x\in E\mid P\{\omega\in \Omega\mid \|x\|_{Q}(\omega)<\varepsilon\}>1-\lambda\}$. Then   $\mathcal{U}(\mathcal{P}):=\{U(Q,\varepsilon,\lambda)\mid Q\in\mathcal{P}_{f},~\varepsilon>0,~0<\lambda<1\}$ is a local base of some Hausdorff linear topology on $E$, called the $(\varepsilon,\lambda)$-topology induced by $\mathcal{P}$ and denoted by $\mathcal{T}_{\varepsilon,\lambda}(\mathcal{P})$. Moreover, the following statements hold:
    \begin{itemize}[widest={(1)}, leftmargin=*]
    	\item [(1)] $(L^{0}(\mathcal{F},\mathbb{K}),\mathcal{T}_{\varepsilon,\lambda}(|\cdot|))$ is a topological algebra;
    	\item [(2)] $(E,\mathcal{T}_{\varepsilon,\lambda}(\mathcal{P}))$ is a topological module over $(L^{0}(\mathcal{F},\mathbb{K}),\mathcal{T}_{\varepsilon,\lambda}(|\cdot|))$; 
    	\item [(3)] $(E,\mathcal{T}_{\varepsilon,\lambda}(\|\cdot\|))$ is metrizable for any $RN$ module $(E,\|\cdot\|)$.
    \end{itemize}
    \end{proposition}

    \par 
    Let $G$ be a nonempty subset of an $L^{0}(\mathcal{F}, \mathbb{K})$-module $E$. Recall that $G$ is said to be $L^{0}$-convex if $\xi x+ (1-\xi) y\in G$ for any $x,y\in G$ and any $\xi\in L^{0}_{+}(\mathcal{F})$ with $0\leq \xi \leq 1$. For a nonempty subset $H$ of $E$, its $L^{0}$-convex hull is defined by
    $$\operatorname{Conv}_{L^{0}}(H)=\left\{\sum_{i=1}^{n}\xi_{i} h_{i} \middle|\ n\in\mathbb{N}, h_{i}\in H, \xi_{i}\in L^{0}_{+}(\mathcal{F})~\text{and}~\sum_{i=1}^{n}\xi_{i}=1\right\}.$$  
    Clearly, $G$ is $L^{0}$-convex if and only if $\operatorname{Conv}_{L^{0}}(G)=G$. Recall further that a mapping $f:G\to L^{0}(\mathcal{F})$ is said to be $L^{0}$-convex (resp., $L^{0}$-concave) if $G$ is $L^{0}$-convex and 
    $$f\left(\sum_{i=1}^{n}\xi_{i}x_{i}\right)\leq\sum_{i=1}^{n}\xi_{i}f(x_{i})~\left(\text{resp.,}~f\left(\sum_{i=1}^{n}\xi_{i}x_{i}\right)\geq\sum_{i=1}^{n}\xi_{i}f(x_{i})\right)$$
    for any $n\in\mathbb{N}$, $x_{1},\ldots,x_{n}\in G$ and $\xi_{1},\ldots,\xi_{n}\in L^{0}_{+}(\mathcal{F})$ with $\sum_{i=1}^{n}\xi_{i}=1$.

    \begin{proposition} [\cite{FKV2009}]\label{prop.topology2}
    Let $(E,\mathcal{P})$ be an $RLC$ module over $\mathbb{K}$ with base $(\Omega,\mathcal{F},P)$. For any $Q\in\mathcal{P}_{f}$ and any $r\in L^{0}_{++}(\mathcal{F})$, let $V(Q,r)=\{ x\in E\mid \|x\|_{Q} <r \text{ on } 1\}$ and $\mathcal{V}(\mathcal{P}):=\{ V(Q,r)\mid Q\in\mathcal{P}_{f} \text{ and }r\in L^{0}_{++}(\mathcal{F})\}$. Then $\{ x+V(Q,r)\mid x\in E, V(Q,r)\in \mathcal{V}(\mathcal{P})\}$ forms a base for some Hausdorff topology on $E$, called the locally $L^{0}$-convex topology induced by $\mathcal{P}$ and denoted by $\mathcal{T}_{c} (\mathcal{P})$. Further, we have the following statements:
    \begin{itemize}[widest={(1)}, leftmargin=*]
    	\item [(1)] $(L^{0}(\mathcal{F},\mathbb{K}),\mathcal{T}_{c}(|\cdot|))$ is a topological ring;
    	\item [(2)] $(E,\mathcal{T}_{c}(\mathcal{P}))$ is a topological module over $(L^{0}(\mathcal{F},\mathbb{K}),\mathcal{T}_{c}(|\cdot|))$; 
    	\item [(3)] $\{ V(Q,r)\mid Q\in\mathcal{P}_{f} \text{ and }r\in L^{0}_{++}(\mathcal{F})\}$ is a local base at $\theta$ of $\mathcal{T}_{c}(\mathcal{P})$.
    \end{itemize}
    \end{proposition}
    
    \par 
    Let $(E,\mathcal{P})$ be an $RLC$ module. The $(\varepsilon,\lambda)$-topology $\mathcal{T}_{\varepsilon,\lambda}(\mathcal{P})$ is an abstract generalization of the topology of convergence in probability. This topology is natural and convenient and has been used in the development of nonsmooth differential geometry on metric measure spaces \cite{BPS2023,CLPV2025,CGP2025,GLP2025,LP2019,LPV2024}. The locally $L^{0}$-convex topology $\mathcal{T}_{c}(\mathcal{P})$ is finer than $\mathcal{T}_{\varepsilon,\lambda}(\mathcal{P})$, and each $V(Q,r)$ in Proposition~\ref{prop.topology2} is $L^{0}$-convex. However, $\mathcal{T}_{c}(\mathcal{P})$ is not in general a linear topology over $\mathbb{K}$ since scalar multiplication is not necessarily continuous, see \cite{FKV2009,Guo2010} for details.

    \par 
    Let $(E,\mathcal{P})$ be an $RLC$ module. If $(\Omega,\mathcal{F},P)$ is trivial, namely, $\mathcal{F}=\{\emptyset,\Omega\}$, then $E$ reduces to a locally convex space, and the $(\varepsilon,\lambda)$-topology and the locally $L^{0}$-convex topology both reduce to the classical locally convex topology on $E$. Thus an $RLC$ module is precisely a random generalization of a classical locally convex space. Definition~\ref{defn.topological L0-module} below gives the corresponding random version of the classical notion of a topological linear space.

    \begin{definition}\label{defn.topological L0-module}
    	Let $E$ be an $L^{0}(\mathcal{F},\mathbb{K})$-module and $\mathcal{T}$ be a topology on $E$. The pair $(E,\mathcal{T})$ is called a topological $L^{0}$-module if the module operations
    	\begin{itemize}[widest={(ii)}, leftmargin=*]
    		\item[(i)] $(E,\mathcal{T})\times(E,\mathcal{T})\rightarrow(E,\mathcal{T}),\quad (x,y)\mapsto x+y$;
    		\item[(ii)] $(L^{0}(\mathcal{F},\mathbb{K}),\mathcal{T}_{c}(|\cdot|))\times(E,\mathcal{T})\rightarrow(E,\mathcal{T}),\quad (\xi,x)\mapsto \xi x$
    	\end{itemize}
    	are continuous with respect to the corresponding product topologies.
    \end{definition}

    \par 
    For a topological $L^0$-module $(E,\mathcal T)$, it is easy to check that the mapping $x\mapsto \xi x+z$ is continuous for any $\xi\in L^0(\mathcal F,\mathbb K)$ and $z\in E$. Moreover, each translation $x\mapsto x+z$ is a homeomorphism of $E$.

    \begin{remark}
    The two topologies $\mathcal{T}_{c}(|\cdot|)$ and $\mathcal{T}_{\varepsilon,\lambda}(|\cdot|)$ on $L^{0}(\mathcal{F},\mathbb{K})$ give rise to two notions of a topological $L^{0}$-module: one over the topological ring $(L^{0}(\mathcal{F},\mathbb{K}),\mathcal{T}_{c}(|\cdot|))$ and the other over the topological algebra $(L^{0}(\mathcal{F},\mathbb{K}),\mathcal{T}_{\varepsilon,\lambda}(|\cdot|))$. Since $\mathcal{T}_{c}(|\cdot|)$ is finer than $\mathcal{T}_{\varepsilon,\lambda}(|\cdot|)$, every topological $L^{0}$-module over $(L^{0}(\mathcal{F},\mathbb{K}),\mathcal{T}_{\varepsilon,\lambda}(|\cdot|))$ is also a topological $L^{0}$-module over $(L^{0}(\mathcal{F},\mathbb{K}),\mathcal{T}_{c}(|\cdot|))$. We therefore adopt the more general notion in Definition~\ref{defn.topological L0-module}. Moreover, by Propositions~\ref{prop.topology1} and~\ref{prop.topology2}, for every $RLC$ module $(E,\mathcal{P})$, both $(E,\mathcal{T}_{\varepsilon,\lambda}(\mathcal{P}))$ and $(E,\mathcal{T}_{c}(\mathcal{P}))$ are topological $L^{0}$-modules in this sense.
    \end{remark}

    \subsubsection{Stable topologies and stable compactness}\label{sec.2.1.2}

    We next recall the notions of $\sigma$-stability, stable topologies and stable compactness. We begin with partitions of unity in $B_{\mathcal F}$
    and the regularity of $L^0(\mathcal F,\mathbb K)$-modules. 
    
    \par 
    A nonempty subset $\{a_{i}\mid i\in I\}$ of $B_{\mathcal{F}}$ is called a partition of unity in $B_{\mathcal{F}}$ if $\bigvee_{i\in I} a_{i} =1$ and $a_{i} \wedge a_{j} = 0$ whenever $i\neq j$; the family of all such partitions is denoted by $p(1)$. It is well known that, for any $\{a_{j}\mid j\in J\}\in p(1)$, the set $\{j\in J \mid  a_{j}>0\}$ is at most countable, and thus any partition of unity in $B_{\mathcal{F}}$ is always expressed in a form $\{a_{n},n\in \mathbb{N}\}$. Moreover, it is clear that $\{a_{n},n\in \mathbb{N}\}\in p(1)$ iff  there exists a countable measurable partition $\{A_{n}, n\in \mathbb{N} \}$ of $\Omega$ such that $a_{n}=[A_{n}]$ for each $n\in \mathbb{N}$.  More generally, for $a\in B_{\mathcal{F}}$, let $p(a)$ denote the family of all partitions of $a$, that is, all families $\{a_{i}\mid i\in I\}$ in $B_{\mathcal{F}}$ such that
    $$\bigvee_{i\in I}a_{i}=a,\qquad a_{i}\wedge a_{j}=0\quad\text{whenever }i\neq j.$$
    
    \par 
    For any $a=[A]\in B_{\mathcal{F}}$, we use $I_{a}$ to denote the equivalence class of $I_{A}$, where $I_{A}$ denotes the characteristic function of $A$, namely, $I_{A}(\omega)=1$ for $\omega\in A$ and $I_{A}(\omega)=0$ otherwise.
    
    \par 
    Recall that an  $L^{0}(\mathcal{F},\mathbb{K})$-module $E$ is said to be regular \cite{Guo2010} if $E$ has the following property: for any given two elements $x$ and $y$ in $E$, if there exists $\{a_{n},n\in \mathbb{N}\}\in p(1)$ such that $I_{a_{n}}x=I_{a_{n}}y$ for each $n\in \mathbb{N}$, then $x=y$.
    
    \par 
    In the remainder of this paper, all the $L^{0}(\mathcal{F},\mathbb{K})$-modules under consideration are assumed to be regular. This restriction is not excessive, since all $RLC$ modules are regular.

    \begin{definition}[\cite{Guo2010}]\label{defn.stable set}
    Let $E$ be an $L^{0}(\mathcal{F},\mathbb{K})$-module. A nonempty subset $G$ of $E$ is said to be
    \begin{itemize}[widest={(1)}, leftmargin=*]
    	\item [(1)] finitely stable if $I_{a}x+I_{a^{c}}y\in G$ for any  $x,y\in G$ and any $a\in B_{\mathcal{F}}$.
    	\item [(2)] $\sigma$-stable if, for any sequence $\{x_{n}, n\in \mathbb{N} \}$ in $G$ and any $\{a_{n},n\in\mathbb{N}\} \in p(1)$, there exists $x\in G$ such that  $I_{a_{n}}x=I_{a_{n}}x_{n}$ for each $n\in \mathbb{N}$. 
    \end{itemize}
    \end{definition}

    \par 
    The property now called $\sigma$-stability was originally called the countable concatenation property in \cite{Guo2010}. By regularity, the element $x$ in part (2) of Definition~\ref{defn.stable set} is unique and is called the countable concatenation of $(x_{n})$ along $\{a_{n}, n\in\mathbb{N}\}$ and denoted by
    $\sum_{n=1}^{\infty}I_{a_{n}}x_{n}$. For a nonempty subset $H$ of a $\sigma$-stable set $G$, the set
    $$\sigma(H):=\left\{\,\sum_{n=1}^{\infty}I_{a_{n}}h_{n}\;\Big|\;\{h_{n}, n\in \mathbb{N} \}\subseteq H,\ \{a_{n}, n\in \mathbb{N} \}\in p(1)\right\}$$
    is called the $\sigma$-stable hull of $H$. Clearly, $H$ is $\sigma$-stable if and only if $\sigma(H)=H$. A simple example of a $\sigma$-stable set is  $L^{0}(\mathcal{F})$.

    \begin{definition}[\cite{GMPTY2026}]\label{defn.Cartesian product}
    Let $I$ be a nonempty index set and $\{E_{i}\mid i\in I\}$ a family of $L^{0}(\mathcal{F},\mathbb{K})$-modules. Their Cartesian product $E:=\prod_{i\in I}E_{i}$ is also an $L^{0}(\mathcal{F},\mathbb{K})$-module, where addition and scalar multiplication are defined coordinatewise. If $G_{i}$ is a nonempty $\sigma$-stable subset of $E_{i}$ for each $i\in I$, then $G:=\prod_{i\in I}G_{i}$ is a nonempty $\sigma$-stable subset of $E$ with countable concatenations given by
    $$\sum_{m=1}^{\infty}I_{a_{m}}x^{m}=\left(\sum_{m=1}^{\infty}I_{a_{m}}x_{i}^{m} \right)_{i\in I}$$
    for any sequence $\{(x_{i}^{m})_{i\in I},m\in \mathbb{N}\}$ in $G$ and any $\{a_{m}, m\in\mathbb{N}\}\in p(1)$.	
    \end{definition}

    \par 
    Let $(E,\mathcal{P})$ be an $RLC$ module over $\mathbb{K}$ with base $(\Omega,\mathcal{F},P)$. For any sequence $\{Q_{n}, n\in\mathbb{N}\}$ in $\mathcal{P}_{f}$ and any $\{a_{n}, n\in\mathbb{N}\}\in p(1)$, the mapping $\sum_{n=1}^{\infty}I_{a_{n}}\|\cdot\|_{Q_{n}}:E\to L^{0}_{+}(\mathcal{F})$ defined by 
    $$\left(\sum_{n=1}^{\infty}I_{a_{n}}\|\cdot\|_{Q_{n}}\right)(x)=\sum_{n=1}^{\infty}I_{a_{n}}\|x\|_{Q_{n}}, \forall x\in E$$ 
    is still an $L^{0}$-seminorm on $E$, and the family of all such $L^{0}$-seminorms is denoted by $\mathcal{P}_{cc}$. Clearly, $\mathcal{P}\subseteq\mathcal{P}_{cc}$ and it is known from \cite{GMPTY2026,Guo2010} that 
    $$\mathcal{T}_{\varepsilon,\lambda}(\mathcal{P}_{cc})=\mathcal{T}_{\varepsilon,\lambda}(\mathcal{P})\subseteq \mathcal{T}_{c}(\mathcal{P})\subseteq \mathcal{T}_{c}(\mathcal{P}_{cc}).$$
    We write $\mathcal{V}(\mathcal{P}_{cc})$ for the local base at $\theta$ specified in Proposition~\ref{prop.topology2}, with $\mathcal{P}$ replaced by $\mathcal{P}_{cc}$.

    \par
    It is known from \cite{Guo2010} that every $\mathcal{T}_{\varepsilon,\lambda}(\mathcal{P})$-closed $L^{0}$-convex subset of a $\mathcal{T}_{\varepsilon,\lambda}(\mathcal{P})$-complete $RLC$ module $(E,\mathcal{P})$ is $\sigma$-stable. Proposition \ref{prop.stable} below gives two relations among the three topologies $\mathcal{T}_{\varepsilon,\lambda}(\mathcal{P})$, $\mathcal{T}_{c}(\mathcal{P})$ and  $\mathcal{T}_{c}(\mathcal{P}_{cc})$.

    \begin{proposition}[\cite{Guo2010,GZZ2015}]\label{prop.stable}
    Let $(E,\mathcal{P})$ be an $RLC$ module over $\mathbb{K}$ with base $(\Omega,\mathcal{F},P)$ and $G$ be a nonempty subset of $E$. Then we have the following:
    \begin{itemize}[widest={(1)}, leftmargin=*]
    	\item [(1)] If $G$ is $\sigma$-stable, then $cl_{\mathcal{T}_{\varepsilon,\lambda}(\mathcal{P})}(G)=cl_{\mathcal{T}_{c}(\mathcal{P})}(G)=cl_{\mathcal{T}_{c}(\mathcal{P}_{cc})}(G)$, where $cl_{\mathcal{T}_{\varepsilon,\lambda}(\mathcal{P})}(G),~cl_{\mathcal{T}_{c}(\mathcal{P})}(G)$ and $cl_{\mathcal{T}_{c}(\mathcal{P}_{cc})}(G)$ denote the closure of $G$ under $\mathcal{T}_{\varepsilon,\lambda}(\mathcal{P}),~\mathcal{T}_{c}(\mathcal{P})$ and $\mathcal{T}_{c}(\mathcal{P}_{cc})$, respectively.
    	\item [(2)] If $G$ is finitely stable, then $G$ is $\mathcal{T}_{\varepsilon,\lambda}(\mathcal{P})$-complete iff $G$ is both $\sigma$-stable and $\mathcal{T}_{c}(\mathcal{P})$-complete. Moreover, if $G$ is $\sigma$-stable, then $G$ is $\mathcal{T}_{\varepsilon,\lambda}(\mathcal{P})$-complete iff $G$ is $\mathcal{T}_{c}(\mathcal{P})$-complete iff $G$ is $\mathcal{T}_{c}(\mathcal{P}_{cc})$-complete.
    	\end{itemize}
    \end{proposition}

    \begin{remark}\label{rmk.same closedness}
    By Proposition~\ref{prop.stable}, the closedness and completeness of a $\sigma$-stable subset $G$ of an $RLC$ module $(E,\mathcal{P})$ are independent of the choice among $\mathcal{T}_{\varepsilon,\lambda}(\mathcal{P})$, $\mathcal{T}_{c}(\mathcal{P})$ and $\mathcal{T}_{c}(\mathcal{P}_{cc})$. Hence, in the remainder of this paper, we shall refer to such a $\sigma$-stable subset as closed or complete without specifying the topology.
    \end{remark}

    \par
    To define stable topological spaces and stable compactness, we first recall some basic notions concerning $\sigma$-stable families of sets.

    \begin{definition}[\cite{GMPTY2026}]\label{defn.stable family}
   	 Let $G$ be a $\sigma$-stable subset of an $L^{0}(\mathcal{F},\mathbb{K})$-module.
   	 \begin{itemize}[widest={(1)}, leftmargin=*]
   	 \item [(1)] For any sequence of nonempty subsets $\{G_{n}, n \in \mathbb{N}\}$ of $G$ and any $\{a_{n}, n \in \mathbb{N}\}\in  p(1)$,  $$\sum^{\infty}_{n=1} I_{a_{n}} G_{n}:=\left\{\sum^{\infty}_{n=1} I_{a_{n}} x_{n}\mid  x_{n}\in G_{n}, \forall~ n\in\mathbb{N}\right\}$$ 
   	 is called the countable concatenation of $\{G_{n},n\in\mathbb{N}\}$ along $\{a_{n},n\in \mathbb{N}\}$. It is easy to check that $\sum^{\infty}_{n=1} I_{a_{n}} G_{n}$ is $\sigma$-stable if each $G_{n}$ is $\sigma$-stable.
   	 \item [(2)] For a family $\mathcal{E}$ of nonempty subsets of $G$, $$\sigma(\mathcal{E}):=\left\lbrace \sum_{n=1}^{\infty}I_{a_{n}}G_{n} \middle |\ \{G_{n},~n\in\mathbb{N}\}\subseteq\mathcal{E}~\text{and}~\{a_{n},~n\in\mathbb{N}\}\in  p(1)\right\rbrace $$ 
   	 is called the $\sigma$-stable hull of $\mathcal{E}$. $\mathcal{E}$ is said to be $\sigma$-stable if $\sigma(\mathcal{E})=\mathcal{E}$.
   	 \item [(3)] A family $\mathcal{E}$ of nonempty subsets of $G$ is said to be $\sigma$-bistable if $\mathcal{E}$ and each element of $\mathcal{E}$ are $\sigma$-stable.
   	 \end{itemize}
    \end{definition}

   \begin{definition}[\cite{GMPTY2026}]\label{defn.stable filter}
   	Let $X$ be a nonempty $\sigma$-stable subset of an $L^{0}(\mathcal{F},\mathbb{K})$-module $E$. A filter $\mathscr{F}$ on $X$ is said to be stable if it admits a $\sigma$-bistable base. A stable filter on $X$ is called a stable ultrafilter if it is not properly contained in any other stable filter on $X$.
   \end{definition}

   \begin{definition}[\cite{GMPTY2026}]\label{defn.stable topology}
   	Let $X$ be a nonempty subset of an $L^{0}(\mathcal{F},\mathbb{K})$-module and $\mathcal{T}$ be a topology on $X$. $(X,\mathcal{T})$ is called a stable topological space if $X$ is $\sigma$-stable and $\mathcal{T}$ has a $\sigma$-bistable base.
   \end{definition}

   \begin{example}\label{ex.stable topology space}
   	Let $(E,\mathcal{P})$ be an $RLC$ module over $\mathbb{K}$ with base $(\Omega,\mathcal{F},P)$ and $G$ be a nonempty $\sigma$-stable subset of $E$. For any $x\in G$, $q\in\mathcal{P}_{cc}$ and $r\in L_{++}^{0}(\mathcal{F})$, let $V_{G}(x,q,r):=\{y\in G\mid q(x-y)<r~\text{on}~1\}$. Then $\{V_{G}(x,q,r)\mid x\in G,~q\in\mathcal{P}_{cc},~r\in L_{++}^{0}(\mathcal{F})\}$ is a base of $\mathcal{T}_{c}(\mathcal{P}_{cc})|_{G}$ (namely, the relativization of $\mathcal{T}_{c}(\mathcal{P}_{cc})$ to $G$). It is clear that each $V_{G}(x,q,r)$ is $\sigma$-stable. Moreover, for any sequence $\{x_{n},n\in\mathbb{N}\}$ in $G$, any sequence $\{q_{n},n\in\mathbb{N}\}$ in $\mathcal{P}_{cc}$, any sequence $\{r_{n},n\in\mathbb{N}\}$ in $L_{++}^{0}(\mathcal{F})$ and any $\{a_{n},n\in\mathbb{N}\}\in p(1)$, we have $$\sum_{n=1}^{\infty}I_{a_{n}}V_{G}(x_{n},q_{n},r_{n})=V_{G}\left(\sum_{n=1}^{\infty}I_{a_{n}}x_{n},\sum_{n=1}^{\infty}I_{a_{n}}q_{n},\sum_{n=1}^{\infty}I_{a_{n}}r_{n}\right).$$
   	Thus, $(G,\mathcal{T}_{c}(\mathcal{P}_{cc})|_{G})$ (briefly, denoted by $(G,\mathcal{T}_{c}(\mathcal{P}_{cc}))$) is a stable topological space. In general, neither $(G,\mathcal{T}_{\varepsilon,\lambda}(\mathcal{P}))$ nor $(G,\mathcal{T}_{c}(\mathcal{P}))$ is stable. 
   \end{example}

   \par 
   To formulate stable compactness in terms of open covers, we first recall the notion of a stably finite family.

   \begin{definition}[\cite{GMPTY2026}]\label{defn.stably finite}
   Let $G$ be a $\sigma$-stable subset of an $L^{0}(\mathcal{F},\mathbb{K})$-module and $\mathcal{E}$ be a $\sigma$-stable family of nonempty $\sigma$-stable subsets of $G$. A nonempty subfamily $\tilde{\mathcal{E}}$ of $\mathcal{E}$ is said to be stably finite if there exist $\{a_{n},n\in \mathbb{N}\}\in  p(1)$ and a sequence $\{\mathcal{E}_{n},n\in \mathbb{N}\}$ of nonempty finite subfamilies of $\mathcal{E}$ such that $\tilde{\mathcal{E}}=\{\sum_{n} I_{a_{n}}H_{n}:H_{n}\in\sigma(\mathcal{E}_{n})\}$.
   \end{definition}

   \begin{definition}[\cite{GMPTY2026}]\label{defn.stable compactness}
   	A stable topological space $(X,\mathcal{T})$ is said to be stably compact if every $\sigma$-bistable open cover of $X$ has a stably finite subcover. A $\sigma$-stable subset $Y$ of a topological space $(X_{1},\mathcal{T}_{1})$ is said to be stably compact if $(Y,\mathcal{T}_{1}|_{Y})$ is stably compact.
   \end{definition}

   \par
   In \cite[Definition 2.10]{GWXYC2025}, Guo et al. introduced the notion of a random sequentially compact subset $G$ of an $RN$ module $(E,\|\cdot\|)$. It is proved in \cite[Theorem 2.50]{GMPTY2026} that a $\sigma$-stable subset $G$ of $E$ is random sequentially compact iff $(G,\mathcal{T}_{c}(\|\cdot\|))$ is stably compact.

   \begin{proposition}[\cite{GMPTY2026}]\label{prop.stable compactness}
   	Let $(X,\mathcal{T})$ be a stable topological space. Then the following statements are equivalent to each other:
   	\begin{itemize}[widest={(1)}, leftmargin=*]
   		\item [(1)] $(X,\mathcal{T})$ is stably compact;
   		\item [(2)] Every $\sigma$-stable family of $\sigma$-stable closed subsets of $X$ has a nonempty intersection whenever the family has the finite intersection property;
   		\item [(3)] Every family of $\sigma$-stable closed subsets of $X$ has a nonempty intersection whenever the family has the finite intersection property;
   		\item [(4)] Every stable filter on $X$ admits a cluster point.
   	\end{itemize}
   \end{proposition}

   \begin{remark}\label{rmk.stable net}
   	Let $(X,\mathcal{T})$ be a stable topological space. In parallel with the classical characterization of compactness by nets, $(X,\mathcal{T})$ is stably compact if and only if every $\sigma$-stable net in $X$ admits a convergent $\sigma$-stable subnet (see \cite[Theorem~2.33]{GMPTY2026} for details). Since the definitions of a $\sigma$-stable net and a $\sigma$-stable subnet are somewhat involved and are nowhere used in this paper, we do not reproduce them here.
   \end{remark}

   \begin{remark}\label{rmk.same stable compactness}
   	Let $(E,\mathcal{P})$ be an $RLC$ module and $G$ be a nonempty $\sigma$-stable subset of $E$. As shown in Example \ref{ex.stable topology space},  $(G,\mathcal{T}_{\varepsilon,\lambda}(\mathcal{P}))$ and $(G,\mathcal{T}_{c}(\mathcal{P}))$ need not be stable. However, (1) of Proposition \ref{prop.stable} and (3) of Proposition \ref{prop.stable compactness} ensure that the following definition is well defined: a $\sigma$-stable subset $G$ of $E$ is said to be $\mathcal{T}_{\varepsilon,\lambda}(\mathcal{P})$-stably compact (or $\mathcal{T}_{c}(\mathcal{P})$-stably compact) if every family of $\sigma$-stable closed subsets of $G$ with the finite intersection property has a nonempty intersection. Unless otherwise stated, stable compactness for subsets of an $RLC$ module $(E,\mathcal{P})$ is understood with respect to $\mathcal{T}_{c}(\mathcal{P}_{cc})$.
   \end{remark}

   \begin{remark}
   The notions recalled above are treated in a more general setting in \cite{GMPTY2026}, where $B$ is an arbitrary complete Boolean algebra and $B$-stability is defined for a nonempty set $G$ by means of a $B$-regular equivalence relation on $G\times B$, without requiring $G$ to be a subset of an $L^{0}(\mathcal{F},\mathbb{K})$-module. Since this paper deals only with $L^{0}(\mathcal{F},\mathbb{K})$-modules, we consider the special case $B=B_{\mathcal{F}}$, with the regular equivalence relation as described in \cite[Example 2.4]{GMPTY2026}. Accordingly, the notions of $\sigma$-stability, stable topological spaces and stable compactness in Definitions~\ref{defn.stable set}, \ref{defn.stable topology} and~\ref{defn.stable compactness} are special cases of Definitions 2.1, 2.25 and 2.32 in \cite{GMPTY2026}, respectively, and Proposition~\ref{prop.stable compactness} is the corresponding special case of \cite[Theorem 2.33~and~Remark 2.34]{GMPTY2026}.
   \end{remark}


   \subsection{The random KKM mapping principle}\label{sec.2.2}

   Throughout this paper, $L^{0}(\mathcal{F},\mathbb{R}^{n})$ denotes the linear space of equivalence classes of $\mathbb{R}^{n}$-valued random vectors on
   $(\Omega,\mathcal{F},P)$. It is also a free module of rank $n$ over the algebra $L^{0}(\mathcal{F},\mathbb{R})$. Equipped with the $L^{0}$-norm $\|\cdot\|$ defined by
   $$\|x\|=\left(\sum_{i=1}^{n}|\xi_{i}|^{2}\right)^{\frac{1}{2}},~\forall x=(\xi_{1},\ldots,\xi_{n}),$$
   $(L^{0}(\mathcal{F},\mathbb{R}^{n}),\|\cdot\|)$ is a complete $RN$ module over $\mathbb{R}$ with base $(\Omega,\mathcal{F},P)$ and is called a random Euclidean space. The standard simplex in $L^{0}(\mathcal{F},\mathbb{R}^{n})$ is defined by 
   $$\Delta_{n-1}=\operatorname{Conv}_{L^{0}}(\{e_{1},e_{2},\cdots,e_{n}\}),$$
   where  $e_{i}$ denotes the $i$-th standard basis vector of $L^{0}(\mathcal{F},\mathbb{R}^{n})$. 
   
   \par 
   A subset $G$ of an $RLC$ module $(E,\mathcal{P})$ is said to be a.s. bounded if $\bigvee\{\|g\|\mid g\in G\}\in L_{+}^{0}(\mathcal{F})$ for every $\|\cdot\|\in\mathcal{P}$. In particular, $\Delta_{n-1}$ is a.s. bounded, since $\|x\|\leq 1$ for every $x\in\Delta_{n-1}$. By \cite[Remark 4.4]{TMGC2026}, $\Delta_{n-1}$ is $\sigma$-stable and random sequentially compact, and hence stably compact by \cite[Theorem 2.50]{GMPTY2026}.
   
   \par  
   Now we present the random KKM lemma on the standard $L^{0}$-simplex in $L^{0}(\mathcal{F},\mathbb{R}^{n})$. This result is a random generalization of Knaster--Kuratowski--Mazurkiewicz's result \cite{KKM1929}.

   \begin{theorem}[The random KKM lemma on the standard $L^{0}$-simplex]\label{thm.KKM1}
   	Let $\Delta_{n-1}=\operatorname{Conv}_{L^{0}}(\{e_{1},e_{2},\cdots,e_{n}\})$ be the standard $(n-1)$-dimensional $L^{0}$-simplex in $L^{0}(\mathcal{F},\mathbb{R}^{n})$ and $\{F_{1},F_{2},\cdots,F_{n}\}$ be a family of nonempty $\sigma$-stable closed subsets of $\Delta_{n-1}$ such that 
   	$$\operatorname{Conv}_{L^{0}}(\{e_{i}\mid i\in I\})\subseteq \bigcup\sigma\{F_{i}\mid i\in I\}$$ 
   	for each nonempty subset $I$ of $\{1,2,\cdots,n\}$. Then $\bigcap_{i=1}^{n}F_{i}$ is nonempty and stably compact.
   \end{theorem}

   \par 
   The covering condition in Theorem \ref{thm.KKM1} motivates the following notion of a random KKM mapping (namely, Definition \ref{defn.KKM mapping}). To introduce this notion, we first define $\sigma$-stable set-valued mappings.

   \begin{definition}\label{defn.stable mapping}
  	Let $E_{1}$ and $E_{2}$ be two $L^{0}(\mathcal{F},\mathbb{K})$-modules, $X\subseteq E_{1}, Y\subseteq E_{2}$  two nonempty sets and  $T:X\rightarrow 2^{Y}\setminus \{\emptyset\}$ a set-valued mapping. $T$ is said to be $\sigma$-stable if both $X$ and $T(X):=\bigcup_{x\in X}T(x)$ are $\sigma$-stable and 
  	$$T\left(\sum_{n=1}^{\infty}I_{a_{n}}x_{n}\right)=\sum_{n=1}^{\infty}I_{a_{n}}T(x_{n})$$ 
  	for any sequence $\{x_{n},n\in \mathbb{N}\}$ in $X$ and any $\{a_{n},n\in \mathbb{N}\}\in p(1)$.
   \end{definition}

   \par 
   For a $\sigma$-stable set-valued mapping $T:X\rightarrow 2^{Y}\setminus \{\emptyset\}$, it is easy to check that $T(x)$ is $\sigma$-stable for each $x\in X$. Definition~\ref{defn.stable mapping} reduces to the usual notion of a $\sigma$-stable single-valued mapping \cite{GWXYC2025} when each $T(x)$ is a singleton.

  \begin{definition}\label{defn.KKM mapping}
  	Let $E$ be an $L^{0}(\mathcal{F},\mathbb{K})$-module over $\mathbb{K}$ with base $(\Omega,\mathcal{F},P)$, $G$ a nonempty $\sigma$-stable subset of $E$ and $T:G\rightarrow 2^{E}\setminus \{\emptyset\}$ a $\sigma$-stable set-valued mapping. $T$ is called a random KKM mapping if for any finite subset $\{g_{1},\cdots,g_{n}\}$ of $G$, it holds that $$\operatorname{Conv}_{L^{0}}\{g_{1},\cdots,g_{n}\}\subseteq\bigcup \{T(g)\mid g\in \sigma\{g_{1},\cdots,g_{n}\}\}.$$
  \end{definition}

  \par 
  In Definition~\ref{defn.KKM mapping}, it is easy to check that $\operatorname{Conv}_{L^{0}}\{g_{1},\ldots,g_{n}\}$ is $\sigma$-stable. Moreover, the $\sigma$-stability of $T$ gives 
  $$\bigcup\{T(g)\mid g\in\sigma\{g_{1},\ldots,g_{n}\}\}=\bigcup\sigma\{T(g_{1}),\ldots,T(g_{n})\}.$$
  Since  $\sigma\{T(g_{1}),\ldots,T(g_{n})\}$ is a $\sigma$-bistable family, its union is $\sigma$-stable \cite{GMPTY2026}. Hence, the covering condition in Definition~\ref{defn.KKM mapping} is compatible with $\sigma$-stability. When $\mathcal{F}=\{\emptyset,\Omega\}$, a random KKM mapping reduces to a classical KKM mapping \cite{Granas2015}.

  \par 
  Theorem \ref{thm.KKM1} leads to the general form of the random KKM mapping principle below, which is a random generalization of Fan's result \cite{Fan1961}.

  \begin{theorem}[The random KKM mapping principle on topological $L^{0}(\mathcal{F},\mathbb{K})$-modules]\label{thm.KKM2} 
  	Let $(E,\mathcal{T})$ be a Hausdorff topological $L^{0}(\mathcal{F},\mathbb{K})$-module with base $(\Omega,\mathcal{F},P)$, $G$ a nonempty $\sigma$-stable subset of $E$ and $T:G\rightarrow 2^{E}\setminus \{\emptyset\}$ a random KKM mapping with $\sigma$-stable $\mathcal{T}|_{T(G)}$-closed values. Then $\{T(g)\mid g\in G\}$ has the finite intersection property. Moreover, if there exists $g_{0}\in G$ such that $T(g_{0}) $ is stably compact, then $\bigcap \{T(g)\mid g\in G\}\neq \emptyset$.
  \end{theorem}

  \par 
  For later applications, we give a construction of random KKM mappings based on $L^0$-quasiconvexity.

  \begin{definition}\label{defn.quasiconvexity}
  	Let $E$ be an $L^{0}(\mathcal{F},\mathbb{K})$-module and $G$ be a nonempty $L^{0}$-convex subset of $E$. A mapping  $f:G\rightarrow L^{0}(\mathcal{F})$ is said to be $L^{0}$-quasiconvex (resp., $L^{0}$-quasiconcave) if $\{x\in G\mid f(x)\leq \eta\}$ (resp., $\{x\in G\mid f(x)\geq \eta\}$) is $L^{0}$-convex (if nonempty) for any $\eta\in L^{0}(\mathcal{F})$. 
  \end{definition}

  \begin{remark}\label{rmk.quasiconvexity}
  	Let $G$ and $f$ be as in Definition~\ref{defn.quasiconvexity}. Then it is easy to check that $f$ is $L^{0}$-quasiconvex (resp., $L^{0}$-quasiconcave) if and only if
  	$$f\left(\sum_{i=1}^{n}\lambda_{i} x_{i}\right)\leq\bigvee_{i=1}^{n}f(x_{i})~\left(\text{resp.}, f\left(\sum_{i=1}^{n}\lambda_{i} x_{i}\right)\geq\bigwedge_{i=1}^{n}f(x_{i})\right)$$
  	for any $n\in\mathbb{N}$, $x_{1},\ldots,x_{n}\in G$ and $\lambda_{1},\ldots,\lambda_{n}\in L^{0}_{+}(\mathcal{F})$ with $\sum_{i=1}^{n}\lambda_{i}=1$. 	In particular, every $L^{0}$-convex (resp., $L^{0}$-concave) mapping is $L^{0}$-quasiconvex (resp., $L^{0}$-quasiconcave), since $\sum_{i=1}^{n}\lambda_{i}f(x_{i})\leq\bigvee_{i=1}^{n}f(x_{i})$ (resp., $\sum_{i=1}^{n}\lambda_{i}f(x_{i})\geq\bigwedge_{i=1}^{n}f(x_{i})$).
  \end{remark}

  \par
  Lemma \ref{lemm.example KKM} below provides a class of random KKM mappings that will be used in the subsequent proofs.

  \begin{lemma}\label{lemm.example KKM}
  	Let $E$ be an $L^{0}(\mathcal{F},\mathbb{K})$-module, $G$ a $\sigma$-stable and $L^{0}$-convex subset of $E$ and $H:G\times G\rightarrow L^{0}(\mathcal{F})$ a $\sigma$-stable mapping such that $H(y,\cdot)$ is $L^{0}$-quasiconvex for any $y\in G$. Define $T:G\rightarrow 2^{G}\setminus \{\emptyset\}$ by 
  	$$T(x)=\{y\in G\mid H(y,y)\leq H(y,x)\},~\forall x\in G.$$ 
  	Then $T$ is a random KKM mapping.
  \end{lemma}

  \subsection{The random Fan--Browder fixed point and coincidence theorems}\label{sec.2.3}

  We first establish a random generalization of the Fan--Browder fixed point theorem \cite{Fan1969,Browder1968} using the random KKM mapping principle. Besides yielding the coincidence theorem below, this fixed point theorem will be used in the proofs of the intersection theorem, the minimax principle and two extensions of the random Tychonoff fixed point theorem.

  \par
  For a set-valued mapping $T:X\rightarrow 2^{Y}$, the inverse $T^{-1}:Y\rightarrow 2^{X}$ of $T$ is defined by 
  $$T^{-1}(y)=\{x\in X\mid y\in T(x)\},~\forall y\in Y.$$

  \begin{theorem}[Random Fan--Browder fixed point theorem]\label{thm.Fan fixed point}
  	Let $(E,\mathcal{T})$ be a Hausdorff topological $L^{0}(\mathcal{F},\mathbb{K})$-module with base $(\Omega,\mathcal{F},P)$, $G$ a nonempty stably compact $L^{0}$-convex subset of $E$ and $\varphi:G\rightarrow 2^{G}\setminus \{\emptyset\}$ a $\sigma$-stable set-valued mapping such that
  	\begin{itemize}[widest={(1)}, leftmargin=*]
  		\item [(1)] $\varphi^{-1}(y)$ is $\mathcal{T}|_{G}$-open and $\sigma$-stable for each $y\in \varphi(G)$;
  		\item [(2)] $\varphi (x)$ is $L^{0}$-convex for each $x\in G$.
  	\end{itemize}
  	Then there exists $g\in G$ such that $g\in \varphi(g)$.
  \end{theorem}

  \par 
  As a direct application of Theorem~\ref{thm.Fan fixed point}, we obtain the following random coincidence theorem.

  \begin{theorem}\label{thm.coincidence}
  	Let $(E_{1},\mathcal{T}_{1})$ and $(E_{2},\mathcal{T}_{2})$ be two Hausdorff topological $L^{0}(\mathcal{F},\mathbb{K})$-modules with base $(\Omega,\mathcal{F},P)$, $X\subseteq E_{1}$ and $Y\subseteq E_{2}$ two nonempty stably compact $L^{0}$-convex sets and $f,g:X\rightarrow 2^{Y}\setminus \{\emptyset\}$ two $\sigma$-stable set-valued mappings such that
  	\begin{itemize}[widest={(1)}, leftmargin=*]
  		\item [(1)] $f(x)$ is $\sigma$-stable and $\mathcal{T}_{2}|_{Y}$-open and $g(x)$ is a nonempty $L^{0}$-convex set for any $x\in X$;
  		\item [(2)] $g^{-1}(y)$ is $\sigma$-stable and $\mathcal{T}_{1}|_{X}$-open for any $y\in g(X)$ and $f^{-1}(y)$ is a nonempty $L^{0}$-convex set for any $y\in Y$;
  	\end{itemize}
  	Then there exists $x\in X$ such that $f(x)\cap g(x)\neq \emptyset$.
  \end{theorem}

  \subsection{Random minimax results}\label{sec.2.4}

  \par
  We next present two minimax results derived from Theorems \ref{thm.KKM2} and \ref{thm.Fan fixed point}. To state them, we recall some basic notions for $\bar{L}^{0}(\mathcal{F})$-valued functions. Throughout this paper, we adopt the conventions $0\cdot(\pm\infty)=0$ and $(+\infty)+(-\infty)=+\infty$.

  \begin{definition}\label{defn.semicontinuous}
  	Let $(E,\mathcal{T})$ be a topological $L^{0}(\mathcal{F},\mathbb{K})$-module with base $(\Omega,\mathcal{F},P)$ and $G$ a nonempty $\sigma$-stable subset of $E$. A mapping  $f:G\rightarrow \bar{L}^{0}(\mathcal{F})$ is said to be
  	\begin{itemize}[widest={(1)}, leftmargin=*]
  		\item [(1)] proper if $f(x)>-\infty$ on $1$ for all $x\in G$ and there exists some $x\in G$ such that $f(x)\in L^{0}(\mathcal{F})$;
  		\item [(2)] lower semicontinuous if $\{x\in G\mid f(x)\leq\xi\}$ is closed in $(G,\mathcal{T}|_{G})$ for any $\xi\in L^{0}(\mathcal{F})$;
  		\item [(3)] upper semicontinuous if $\{x\in G\mid f(x)\geq\xi\}$ is closed in $(G,\mathcal{T}|_{G})$ for any $\xi\in L^{0}(\mathcal{F})$.
  	\end{itemize}
  \end{definition}

  \begin{remark}\label{rmk.semicontinuous}
  	If the mapping $f$ in Definition \ref{defn.semicontinuous} is $L^{0}(\mathcal{F})$-valued and continuous from $(G,\mathcal{T}|_{G})$ to $(L^{0}(\mathcal{F}),\mathcal{T}_{c}(|\cdot|))$, then $f$ is both lower semicontinuous and upper semicontinuous. Indeed, for each $\xi\in L^{0}(\mathcal{F})$, the sets $\{u\in L^{0}(\mathcal{F})\mid u\leq\xi\}$ and $\{u\in L^{0}(\mathcal{F})\mid u\geq\xi\}$ are  $\mathcal{T}_{c}(|\cdot|)$-closed, and hence their inverse images under $f$ are closed in $(G,\mathcal{T}|_{G})$.
  \end{remark}

  \par 
  The first result is the following random minimax inequality.

  \begin{theorem}\label{thm.Minimax inequality}
  	Let $(E,\mathcal{T})$ be a Hausdorff topological $L^{0}(\mathcal{F},\mathbb{K})$-module with base $(\Omega,\mathcal{F},P)$, $G$ a nonempty stably compact $L^{0}$-convex subset of $E$ and $f:G\times G\rightarrow L^{0}(\mathcal{F})$ a $\sigma$-stable mapping such that
  	\begin{itemize}[widest={(1)}, leftmargin=*]
  		\item [(1)] for each $x\in G$, $f(x,\cdot)$ is lower semicontinuous on $G$;
  		\item [(2)] for each $y\in G$, $f(\cdot,y)$ is $L^{0}$-quasiconcave.
  	\end{itemize}	
  	Then the mapping $y\mapsto\bigvee_{x\in G}f(x,y)$ attains its minimum on $G$, and
  	$$\min_{y\in G}\bigvee_{x\in G}f(x,y)\leq \bigvee_{x\in G}f(x,x).$$
  \end{theorem}

  \par 
  The second result is the following random minimax principle,  obtained from Theorem~\ref{thm.Fan fixed point}.

  \begin{theorem}\label{thm.Minimax principle}
  	Let $(E_{1},\mathcal{T}_{1})$ and $(E_{2},\mathcal{T}_{2})$ be two Hausdorff topological $L^{0}(\mathcal{F},\mathbb{K})$-modules with base $(\Omega,\mathcal{F},P)$, $X\subseteq E_{1}$ and $Y\subseteq E_{2}$ two nonempty stably compact $L^{0}$-convex sets. Let $f:X\times Y\rightarrow L^{0}(\mathcal{F})$ be a $\sigma$-stable mapping such that
  	\begin{itemize}[widest={(1)}, leftmargin=*]
  		\item [(1)] for each $x\in X$, $f(x,\cdot)$ is lower semicontinuous and $L^{0}$-quasiconvex;
  		\item [(2)] for each $y\in Y$, $f(\cdot,y)$ is upper semicontinuous and $L^{0}$-quasiconcave.
  	\end{itemize}
  	Then 
  	$$\max_{x\in X}\min_{y\in Y}f(x,y)=\min_{y\in Y}\max_{x\in X}f(x,y).$$
  \end{theorem}

  \subsection{The random intersection theorem and the conditional Nash equilibrium theorem}\label{sec.2.5}

  We first derive a random intersection theorem from Theorem~\ref{thm.Fan fixed point} and then apply it to establish the existence of conditional Nash equilibria.

  \par
  For a Cartesian product $X=\prod_{i=1}^{n}X_{i}$, we write $X^{i}:=\prod_{j\neq i}X_{j}$ and, for $x=(x_{1},\ldots,x_{n})\in X$, let
  $$x^{i}:=(x_{1},\ldots,x_{i-1},x_{i+1},\ldots,x_{n}).$$
  The coordinate projections $p_{i}$ and $p^{i}$ are given by $p_{i}(x)=x_{i}$ and $p^{i}(x)=x^{i}$, respectively. For $y_{i}\in X_{i}$, we write
  $$(y_{i},x^{i}):=(x_{1},\ldots,x_{i-1},y_{i},x_{i+1},\ldots,x_{n}).$$

  \begin{definition}\label{defn.product}
  	Let $\{(E_{i},\mathcal{T}_{i}):i\in I\}$ be a family of Hausdorff topological $L^{0}(\mathcal{F},\mathbb{K})$-modules with base $(\Omega,\mathcal{F},P)$. The Cartesian product $E=\prod_{i\in I}E_{i}$ is an $L^{0}(\mathcal{F},\mathbb{K})$-module, with addition and scalar multiplication defined coordinatewise. Let $\mathcal{T}$ be the product topology on $E$. Then $(E,\mathcal{T})$ is a Hausdorff topological $L^{0}(\mathcal{F},\mathbb{K})$-module, called the product of $\{(E_{i},\mathcal{T}_{i}):i\in I\}$.
  \end{definition}

  \par
  For nonempty subsets $X_i\subseteq E_i$, $i=1,\ldots,n$, put $X=\prod_{i=1}^nX_i$. The relative topology $\mathcal T|_X$ coincides with the product of the relative topologies $\mathcal T_i|_{X_i}$. If each $X_i$ is stably compact, then $X$ is stably compact by \cite[Theorem 2.40]{GMPTY2026}, since the product topology agrees with the stable product topology for finitely many factors.

  \begin{theorem}\label{thm.intersection}
  	Let $(E,\mathcal{T})$ be the product of $n$ Hausdorff topological $L^{0}(\mathcal{F},\mathbb{K})$-modules $(E_{1},\mathcal{T}_{1}),\ldots,(E_{n},\mathcal{T}_{n})$ with base $(\Omega,\mathcal{F},P)$, and $X_{i}$ a nonempty stably compact $L^{0}$-convex subset of $E_{i}$ for each $i=1,\ldots,n$. Suppose that $G_{1},\ldots,G_{n}$ are $\sigma$-stable subsets of $X=\prod_{i=1}^{n}X_{i}$ such that
  	\begin{itemize}[widest={(1)},leftmargin=*]
  		\item[(1)] for each $x\in X$ and each $i=1,\ldots,n$,
  		$$G_{i}(x)=\{y\in X:(y_{i},x^{i})\in G_{i}\}$$
  		is $L^{0}$-convex and nonempty;
  		\item[(2)] for each $y\in X$ and each $i=1,\ldots,n$,
  		$$G^{i}(y)=\{x\in X:(y_{i},x^{i})\in G_{i}\}$$
  		is $\mathcal{T}|_X$-open and, if nonempty, $\sigma$-stable.
  	\end{itemize}
  	Then $\bigcap_{i=1}^{n}G_{i}\neq\emptyset$.
  \end{theorem}

  \par
  Following \cite{GWXYC2025}, we interpret the data in Theorem~\ref{thm.Nash equilibrium} as an $n$-person conditional game. The probability space $(\Omega,\mathcal{F},P)$ represents the random model, where $\mathcal{F}$ represents the information available when the players choose their strategies. For each player $i$, $X_{i}$ is the strategy set and $f_{i}:X=\prod_{j=1}^{n}X_{j}\to L^{0}(\mathcal{F})$ is the $L^{0}(\mathcal{F})$-valued payoff function, with $f_i(x)$ representing the payoff to player $i$ when the strategy profile $x=(x_1,\ldots,x_n)$ is chosen. A strategy profile $x\in X$ is called a conditional Nash equilibrium if
  $$f_i(x)\geq f_i(y_i,x^i),~\forall y_i\in X_i,~\forall i=1,\ldots,n,$$
  or, equivalently,
  $$f_i(x)=\bigvee\{f_i(y_i,x^i)\mid y_i\in X_i\},~\forall i=1,\ldots,n.$$

  \par
  Applying Theorem~\ref{thm.intersection}, we obtain the following existence theorem for conditional Nash equilibria, which generalizes Nash's equilibrium theorem \cite{Nash1951} to the random setting and gives an affirmative answer to the question of equilibrium existence raised in \cite[Remark~5.7]{GWXYC2025}. Indeed, each $RN$ module considered there is a Hausdorff topological $L^0$-module under its locally $L^0$-convex topology, and the $\mathcal T_c$-topology induced by the product $L^0$-norm on their finite product agrees with the product topology. Moreover, each strategy simplex $S_i$ is random sequentially compact and hence stably compact by \cite[Theorem~2.50]{GMPTY2026}. By \cite[Lemma~4.3]{GWXYC2025}, random sequential continuity of each $\sigma$-stable payoff function is equivalent to $\mathcal T_c$-continuity. Thus condition (A) gives the required $\sigma$-stability and continuity, while condition (B) is precisely the $L^0$-quasiconcavity assumption in the following theorem.

  \begin{theorem}\label{thm.Nash equilibrium}
  	Let $(E,\mathcal{T})$ be the product of $n$ Hausdorff topological $L^{0}(\mathcal{F},\mathbb{K})$-modules $(E_{1},\mathcal{T}_{1}),\ldots,(E_{n},\mathcal{T}_{n})$ with base $(\Omega,\mathcal{F},P)$, and $X_{i}$ a nonempty stably compact $L^{0}$-convex subset of $E_{i}$ for each $i=1,\ldots,n$. For each $i=1,\ldots,n$, let $f_{i}:X=\prod_{j=1}^{n}X_{j}\to L^{0}(\mathcal{F})$ be a mapping such that
  	\begin{itemize}[widest={(1)},leftmargin=*]
  		\item[(1)] $f_{i}$ is $\sigma$-stable and continuous from $(X,\mathcal{T}|_{X})$ to $(L^{0}(\mathcal{F}),\mathcal{T}_{c}(|\cdot|))$;
  		\item[(2)] for each $x\in X$, the mapping $\varphi_{i}:X_{i}\to L^{0}(\mathcal{F})$ defined by
  		$$\varphi_{i}(y_{i})=f_{i}(y_{i},x^{i}),~\forall y_{i}\in X_{i},$$
  		is $L^{0}$-quasiconcave.
  	\end{itemize}
  	Then there exists $x\in X$ such that
  	$$f_{i}(x)=\bigvee\{f_{i}(y_{i},x^{i}):y_{i}\in X_{i}\},~\forall i=1,\ldots,n.$$
  \end{theorem}

  \subsection{The random Tychonoff fixed point theorem and its extensions}\label{sec.2.6}

  We now turn to fixed point results in $RLC$ modules. We first specify the notions of continuity used below and compare them for $\sigma$-stable mappings.

  \begin{definition}\label{defn.continuity}
  	Let $(E_{1},\mathcal{P}_{1})$ and $(E_{2},\mathcal{P}_{2})$ be two $RLC$ modules over $\mathbb{K}$ with base $(\Omega,\mathcal{F},P)$ and $G_{1}\subseteq E_{1}$ and $G_{2}\subseteq E_{2}$ be two nonempty sets. The mapping $f:G_{1}\rightarrow G_{2}$ is said to be 
  	\begin{itemize}[widest={(1)}, leftmargin=*]
  		\item [(1)] $\mathcal{T}_{\varepsilon,\lambda}(\mathcal{P})$-continuous if $f$ is a continuous mapping from $(G_{1},\mathcal{T}_{\varepsilon,\lambda}(\mathcal{P}_{1})|_{G_{1}})$ to $(G_{2},\mathcal{T}_{\varepsilon,\lambda}(\mathcal{P}_{2})|_{G_{2}})$;
  		\item [(2)] $\mathcal{T}_{c}(\mathcal{P})$-continuous if $f$ is a continuous mapping from $(G_{1},\mathcal{T}_{c}(\mathcal{P}_{1})|_{G_{1}})$ to $(G_{2},\mathcal{T}_{c}(\mathcal{P}_{2})|_{G_{2}})$;
  		\item [(3)] $\mathcal{T}_{\varepsilon,\lambda}(\mathcal{P}_{cc})$-continuous if $f$ is a continuous mapping from $(G_{1},\mathcal{T}_{\varepsilon,\lambda}((\mathcal{P}_{1})_{cc})|_{G_{1}})$ to $(G_{2},\mathcal{T}_{\varepsilon,\lambda}((\mathcal{P}_{2})_{cc})|_{G_{2}})$;
  		\item [(4)] $\mathcal{T}_{c}(\mathcal{P}_{cc})$-continuous if $f$ is a continuous mapping from $(G_{1},\mathcal{T}_{c}((\mathcal{P}_{1})_{cc})|_{G_{1}})$ to $(G_{2},\mathcal{T}_{c}((\mathcal{P}_{2})_{cc})|_{G_{2}})$.
  	\end{itemize}
  \end{definition}

  \begin{remark}\label{rmk.continuous}
  	Since $\mathcal{T}_{\varepsilon,\lambda}(\mathcal{P}_{1})=\mathcal{T}_{\varepsilon,\lambda}((\mathcal{P}_{1})_{cc})$ and $\mathcal{T}_{\varepsilon,\lambda}(\mathcal{P}_{2})=\mathcal{T}_{\varepsilon,\lambda}((\mathcal{P}_{2})_{cc})$, $f$ is $\mathcal{T}_{\varepsilon,\lambda}(\mathcal{P})$-continuous if and only if $f$ is $\mathcal{T}_{\varepsilon,\lambda}(\mathcal{P}_{cc})$-continuous.
  \end{remark}

  \par 
  For $\sigma$-stable mappings, these notions of continuity are related by Lemma \ref{lemm.two continuty} below.

  \begin{lemma}\label{lemm.two continuty}
  	Let $(E_{1},\mathcal{P}_{1})$ and $(E_{2},\mathcal{P}_{2})$ be two $RLC$ modules over $\mathbb{K}$ with base $(\Omega,\mathcal{F},P)$, $G_{1}\subseteq E_{1}$ and $G_{2}\subseteq E_{2}$ two $\sigma$-stable sets, and $f:G_{1}\rightarrow G_{2}$ a $\sigma$-stable mapping. Then $f$ is $\mathcal{T}_{\varepsilon,\lambda}(\mathcal{P})$-continuous if and only if it is $\mathcal{T}_{c}(\mathcal{P}_{cc})$-continuous.  Moreover, if $f$ is $\mathcal{T}_{c}(\mathcal{P})$-continuous, then it is $\mathcal{T}_{c}(\mathcal{P}_{cc})$-continuous.
  \end{lemma}

  \begin{remark}
  Let $(E_1,\|\cdot\|_1)$ and $(E_2,\|\cdot\|_2)$ be two $RN$ modules over $\mathbb{K}$ with base $(\Omega,\mathcal{F},P)$, $G_1\subseteq E_1$ and $G_2\subseteq E_2$ two $\sigma$-stable sets, and $T:G_1\rightarrow G_2$ a $\sigma$-stable mapping. In \cite[Lemma~4.3]{GWXYC2025}, Guo et al.\ showed that $T$ is random sequentially continuous if and only if it is $\mathcal{T}_c(\mathcal{P})$-continuous and, in particular, that if $T$ is $\mathcal{T}_{\varepsilon,\lambda}(\mathcal{P})$-continuous, then it is $\mathcal{T}_c(\mathcal{P})$-continuous. Since $\mathcal{P}_i:=\{\|\cdot\|_i\}$ satisfies $(\mathcal{P}_i)_{cc}=\mathcal{P}_i$ for $i=1,2$, Lemma~\ref{lemm.two continuty} shows that $\mathcal{T}_c(\mathcal{P})$-continuity also implies $\mathcal{T}_{\varepsilon,\lambda}(\mathcal{P})$-continuity for $T$. Consequently, random sequential continuity, $\mathcal{T}_{\varepsilon,\lambda}(\mathcal{P})$-continuity and $\mathcal{T}_c(\mathcal{P})$-continuity are equivalent for $T$.
  \end{remark}

  \begin{remark}
  In the sequel, we state the fixed point results for $\mathcal{T}_{c}(\mathcal{P}_{cc})$-continuous mappings. By Lemma~\ref{lemm.two continuty}, these results can equivalently be stated for $\mathcal{T}_{\varepsilon,\lambda}(\mathcal{P})$-continuous mappings. The corresponding versions for $\mathcal{T}_{c}(\mathcal{P})$-continuous mappings also follow immediately.
  \end{remark}
  
  \par
  As an application of the random KKM mapping principle, we obtain the following random Tychonoff fixed point theorem, which solves an open problem posed in \cite{GWXYC2025}.

  \begin{theorem}[Random Tychonoff fixed point theorem]\label{thm.tychonoff1}
  	Let $(E,\mathcal{P})$ be an $RLC$ module over $\mathbb{K}$ with base $(\Omega,\mathcal{F},P)$, $G$ a nonempty stably compact $L^{0}$-convex subset of $E$ and $f:G\rightarrow G$ a $\sigma$-stable $\mathcal{T}_{c}(\mathcal{P}_{cc})$-continuous mapping. Then $f$ has a fixed point.
  \end{theorem}

  \par
  Theorems~\ref{thm.inward1} and~\ref{thm.inward2} below relax the self-mapping assumption in Theorem~\ref{thm.tychonoff1} by allowing mappings $f:G\to E$ subject to suitable geometric conditions. Theorem~\ref{thm.inward1} below is a random generalization of Fan's result \cite{Fan1969}.
  
  \begin{theorem}\label{thm.inward1}
  	Let $(E,\mathcal{P})$ be an $RLC$ module over $\mathbb{K}$ with base $(\Omega,\mathcal{F},P)$, $G$ a nonempty stably compact $L^{0}$-convex subset of $E$  and  $f:G\rightarrow E$ a  $\sigma$-stable $\mathcal{T}_{c}(\mathcal{P}_{cc})$-continuous mapping. Assume that, for each $x\in G$ with $x\neq f(x)$, the $L^{0}$-line segment
  	$$[x,f(x)]:=\left\{\lambda x+(1-\lambda)f(x)\middle | \ \lambda\in L^{0}(\mathcal{F},[0,1])\right\}$$
  	contains at least two distinct points of $G$. Then $f$ has a fixed point.
  \end{theorem}

  \par 
  Definition~\ref{defn.inward} below gives random versions of the notions of inward and outward mappings in the sense of Halpern and Bergman \cite{Hal1968}.

  \begin{definition}\label{defn.inward}
  For an $L^{0}$-convex subset $G$ of an $L^{0}(\mathcal{F},\mathbb{K})$-module $E$ and each $x\in G$, define 
  $$ \operatorname{In}_{G}(x)=\{y\in E\mid \exists~y_{0}\in G~\text{and}~\lambda\in L_{+}^{0}(\mathcal{F})~\text{such~that}~y=x+\lambda(y_{0}-x)\}$$ 
  and 
  $$\operatorname{Out}_{G}(x)=\{y\in E\mid \exists~y_{0}\in G~\text{and}~\lambda\in L_{+}^{0}(\mathcal{F})~\text{such~that}~y=x-\lambda(y_{0}-x)\}.$$ 
  A mapping $f:G\rightarrow E$ is said to be random inward (resp., outward) if $f(x)\in \operatorname{In}_{G}(x)$ (resp., $f(x)\in \operatorname{Out}_{G}(x)$) for any $x\in G$.	
  \end{definition}

  \par 
  As a consequence of Theorem~\ref{thm.inward1}, we obtain a fixed point theorem for random inward and outward mappings. This theorem generalizes Browder's results \cite[Theorems~1 and~2]{Browder1967} from locally convex spaces to random locally convex modules.

  \begin{theorem}\label{thm.inward2}
  	Let $(E,\mathcal{P})$ be an $RLC$ module over $\mathbb{K}$ with base $(\Omega,\mathcal{F},P)$, $G$ a nonempty stably compact $L^{0}$-convex subset of $E$  and  $f:G\rightarrow E$ a  $\sigma$-stable $\mathcal{T}_{c}(\mathcal{P}_{cc})$-continuous random inward (resp., outward) mapping.  Then $f$ has a fixed point.
  \end{theorem}

  \par 
  Theorems~\ref{thm.Schauder-Tychonoff} and~\ref{thm.Fan-Iokhvidov} below follow from Theorem~\ref{thm.Fan fixed point} and give two further
  generalizations of the random Tychonoff fixed point theorem (namely, Theorem~\ref{thm.tychonoff1}).

  \begin{theorem}[Random Schauder--Tychonoff fixed point theorem] \label{thm.Schauder-Tychonoff}
  	Let $(E,\mathcal{P})$ be a $\sigma$-stable $RLC$ module over $\mathbb{K}$ with base $(\Omega,\mathcal{F},P)$, $G$ a nonempty $\sigma$-stable $L^{0}$-convex subset of $E$ and $f:G\rightarrow G$ a $\sigma$-stable $\mathcal{T}_{c}(\mathcal{P}_{cc})$-continuous stably compact mapping (namely, $cl_{\mathcal{T}_{c}(\mathcal{P}_{cc})}(f(G))$ is a stably compact subset of $G$). Then $f$ has a fixed point.
  \end{theorem}

  \begin{theorem}[Random Fan--Iokhvidov theorem]\label{thm.Fan-Iokhvidov}
  	Let $(E,\mathcal{P})$ be a $\sigma$-stable $RLC$ module over $\mathbb{K}$ with base $(\Omega,\mathcal{F},P)$, $G$  a nonempty stably compact $L^{0}$-convex subset of $E$, $K$ a nonempty $\sigma$-stable $L^{0}$-convex closed subset of $E$ and $f:G\rightarrow E$ a $\sigma$-stable $\mathcal{T}_{c}(\mathcal{P}_{cc})$-continuous mapping such that $f(G)\subseteq G+K$. Then there exists $x_{0}\in G$ such that $f(x_{0})\in x_{0}+K$.
  \end{theorem}

\subsection{The random Hartman--Stampacchia theorem}\label{sec.2.7}

  \par 
  We next apply the random KKM mapping principle to variational inequalities in random inner product modules.
  
  \begin{definition}[\cite{Guo1999}]\label{defn.RIPM}
	An ordered pair $(E,\langle \cdot,\cdot \rangle)$ is called a random inner product module (briefly, an $RIP$ module) over $\mathbb{K}$ with base $(\Omega,\mathcal{F},P)$ if $E$ is an $L^{0}(\mathcal{F},\mathbb{K})$-module and $\langle \cdot,\cdot \rangle$ is a mapping from $E\times E$ to $L^{0}(\mathcal{F},\mathbb{K})$ such that the following conditions are satisfied:
  	\begin{enumerate}[label={(\textup{RIP}-\arabic*)},leftmargin=*,itemsep=1pt]
  		\item $\langle x,x \rangle\geq 0,~\forall x\in E$, and $\langle x,x \rangle=0$ iff $x=\theta$;
  		\item $\langle x,y \rangle=\overline{\langle y,x \rangle},~\forall x,y\in E$, where $\overline{\langle y,x \rangle}$ denotes the complex conjugate of $\langle y,x \rangle$;
  		\item $\langle \xi x,y \rangle=\xi \langle x,y \rangle,~\forall\xi\in L^{0}(\mathcal{F},\mathbb{K}),~x,y\in E$;
  		\item $\langle x+y,z \rangle =\langle x,z \rangle +\langle y,z \rangle,~\forall x,y,z\in E$;
  	\end{enumerate}
  	where $\langle x,y \rangle$ is called the random inner product of $x$ and $y$.
  \end{definition}
  
  \par 
  Every $RIP$ module $(E,\langle\cdot,\cdot\rangle)$ becomes an $RN$ module under the $L^{0}$-norm $\|\cdot\|$ defined by
  $$\|x\|=\sqrt{\langle x,x\rangle },~\forall  x\in E.$$
  The $(\varepsilon,\lambda)$-topology and the locally $L^{0}$-convex topology on $E$ are understood to be those of the associated $RN$ module $(E,\|\cdot\|)$.
  
  \par 
   Before stating the existence theorem for random variational inequalities (namely, Theorem \ref{thm.Hartman-Stampacchia}), we recall monotonicity and hemicontinuity for mappings on $RIP$ modules over $\mathbb{R}$. For any two distinct elements $x,y$ of an $L^{0}(\mathcal{F})$-module $E$, the set
   $$l(x,y):=\{\lambda x+(1-\lambda)y\mid \lambda\in L^{0}(\mathcal{F})\}$$
   is called an $L^{0}$-line \cite{WGL2022}.

  \begin{definition}\label{defn.monotone hemi-continuous}
  	Let $(E,\langle \cdot,\cdot\rangle)$ be an $RIP$ module over $\mathbb{R}$ with base $(\Omega,\mathcal{F},P)$ and $G$ be a nonempty subset of $E$. A mapping $f:G\rightarrow E$ is said to be
  	\begin{itemize}[widest={(1)}, leftmargin=*]
  		\item [(1)] monotone on $G$ if 
  		$$\langle f(x)-f(y),x-y\rangle\geq 0,~\forall x,y\in G.$$
  		\item [(2)] $\mathcal{T}_{c}$-hemicontinuous (resp., $\mathcal{T}_{\varepsilon,\lambda}$-hemicontinuous) if its restriction to $l(x,y)\cap G$ is $\mathcal{T}_{c}$-continuous (resp., $\mathcal{T}_{\varepsilon,\lambda}$-continuous) for every pair of distinct elements $x,y\in E$.
  	\end{itemize}
  \end{definition}

  \par 
  By Lemma \ref{lemm.two continuty}, for a $\sigma$-stable mapping, $\mathcal{T}_{\varepsilon,\lambda}$-hemicontinuity implies $\mathcal{T}_{c}$-hemicontinuity.
  We therefore formulate the following existence theorem under $\mathcal{T}_{c}$-hemicontinuity.

  \begin{theorem}[Random Hartman--Stampacchia theorem]\label{thm.Hartman-Stampacchia}
  	Let $(E,\langle\cdot,\cdot\rangle)$ be a $\mathcal{T}_{\varepsilon,\lambda}$-complete $RIP$ module over $\mathbb{R}$ with base $(\Omega,\mathcal{F},P)$, $G$ a nonempty $\mathcal{T}_{\varepsilon,\lambda}$-closed a.s. bounded $L^{0}$-convex subset of $E$ and $f:G\rightarrow E$ a $\sigma$-stable monotone $\mathcal{T}_{c}$-hemicontinuous mapping. Then there exists $y_{0}\in G$ such that $$\langle f(y_{0}),y_{0}-x\rangle\leq 0,~\forall x\in G.$$
  \end{theorem}

\section{Proofs}\label{sec.proofs}

\subsection{Proofs of the results in Subsection~\ref{sec.2.2}}\label{sec.3.1}

In this section, we prove the results of Subsection~\ref{sec.2.2}. We also obtain a set-valued enlargement result and a best approximation result in $RN$ modules.

\par 
Throughout, for any $\xi,\eta\in \bar{L}^{0}(\mathcal{F})$, we use $(\xi=\eta)$ for the equivalence class of $\{\omega\in \Omega\mid \xi^{0}(\omega)=\eta^{0}(\omega)\}$, where $\xi^{0}$ and $\eta^{0}$ are arbitrarily chosen representatives of $\xi$ and $\eta$, respectively. The notations $(\xi<\eta)$ and $(\xi\leq\eta)$ are understood similarly.

\par 
For each $x\in \Delta_{n-1}$, there exists a unique $n$-tuple $(\lambda_{1}^{x},\ldots,\lambda_{n}^{x})\in L^{0}(\mathcal{F},[0,1])^{n}$ such that $x=\sum_{i=1}^{n}\lambda_{i}^{x}e_{i}$ and $\sum_{i=1}^{n}\lambda_{i}^{x}=1$. This tuple is called the $L^{0}$-barycentric coordinate
vector of $x$. Throughout this paper, $\lambda_{i}^{x}$ denotes the $i$-th $L^{0}$-barycentric coordinate of $x$.

\begin{proof}[Proof of Theorem \ref{thm.KKM1}]
	For each $i\in \{1,2,\cdots,n\}$, let
	$$d(x,F_{i}):=\bigwedge\{\|x-y\|\mid y\in F_{i}\},~\forall x\in L^{0}(\mathcal{F},\mathbb{R}^{n}).$$
	It is easy to check that $d(\cdot,F_{i})$ is $\sigma$-stable and
	$$|d(x,F_{i})-d(z,F_{i})|\leq \|x-z\|, \forall x,z\in L^{0}(\mathcal{F},\mathbb{R}^{n}),$$
    so $d(\cdot,F_{i})$ is $\mathcal{T}_{c}(\mathcal{P}_{cc})$-continuous. Define a mapping $f:\Delta_{n-1}\rightarrow \Delta_{n-1}$ by $$f(x)=\frac{1}{\sum_{j=1}^{n}d(x,F_{j})+1}\sum_{i=1}^{n}[d(x,F_{i})+\lambda_{i}^{x}]e_{i}, \forall x\in \Delta_{n-1}.$$ 
    Since each $d(\cdot,F_{i})$ is $\mathcal{T}_{c}(\mathcal{P}_{cc})$-continuous, the mapping $f$ is clearly $\mathcal{T}_{c}(\|\cdot\|)$-continuous. Moreover, for any $\{d_{m},m\in \mathbb{N}\}\in p(1)$ and any sequence $\{x_{m},m\in \mathbb{N}\} $ in $\Delta_{n-1}$, we have
	\begin{align}
		&f\left(\sum_{m=1}^{\infty}I_{d_{m}}x_{m}\right)\nonumber\\
		&=\frac{1}{\sum_{j=1}^{n}d\left(\sum_{m=1}^{\infty}I_{d_{m}}x_{m},F_{j}\right)+1}\sum_{i=1}^{n}\left[ d\left(\sum_{m=1}^{\infty}I_{d_{m}}x_{m},F_{i}\right)+\sum_{m=1}^{\infty}I_{d_{m}}\lambda_{i}^{x_{m}}\right] e_{i}\nonumber\\
		&=\frac{1}{\sum_{j=1}^{n}\sum_{m=1}^{\infty}I_{d_{m}}d(x_{m},F_{j})+1}\sum_{i=1}^{n}\sum_{m=1}^{\infty}I_{d_{m}}[d(x_{m},F_{i})+\lambda_{i}^{x_{m}}] e_{i}\nonumber\\
		&=\sum_{m=1}^{\infty}I_{d_{m}}\left\{\frac{1}{\sum_{j=1}^{n}d(x_{m},F_{j})+1}\sum_{i=1}^{n}[d(x_{m},F_{i})+\lambda_{i}^{x_{m}}] e_{i}\right\}\nonumber\\
		&=\sum_{m=1}^{\infty}I_{d_{m}}f(x_{m}).\nonumber
	\end{align}
	Hence, $f$ is $\sigma$-stable. Since $\Delta_{n-1}$ is $\sigma$-stable, a.s. bounded, closed and $L^{0}$-convex, the random Brouwer fixed point theorem
	\cite[Theorem 2.1]{TMGC2026} yields a fixed point $x_{0}\in\Delta_{n-1}$. Furthermore, since $$x_{0}=\sum_{i=1}^{n}\lambda_{i}^{x_{0}}e_{i}=\frac{1}{\sum_{j=1}^{n}d(x_{0},F_{j})+1}\sum_{i=1}^{n}[d(x_{0},F_{i})+\lambda_{i}^{x_{0}}]e_{i},$$ 
	we have
	\begin{equation}\label{eq.KKM1}
		\lambda_{i}^{x_{0}}=\frac{d(x_{0},F_{i})+\lambda_{i}^{x_{0}}}{\sum_{j=1}^{n}d(x_{0},F_{j})+1},~\forall i\in \{1,2,\cdots,n\}.
	\end{equation}
	Let $a_{1}=\bigvee_{i=1}^{n}(\lambda_{i}^{x_{0}}=1)$ and 
	$$a_{k}=\bigvee_{L\in \mathcal{L}_{k}}\left\{\left[\bigwedge_{l\in L}(0<\lambda_{l}^{x_{0}}< 1)\right]\wedge[\bigwedge_{i\in \{1,2,\cdots,n\}\setminus L}(\lambda_{i}^{x_{0}}=0)]\right\},~\forall  k\in \{2,3,\cdots,n\},$$ 
	where $\mathcal{L}_{k}$ denotes the collection of all $k$-element subsets of $\{1,2,\cdots,n\}$. It is easy to check that $\{a_{k}, k=1,2,\cdots,n\}\in p(1)$.

	\par
	First, we prove that 
	$$I_{a_{1}}d(x_{0},F_{i})=0,~\forall i\in \{1,2,\cdots,n\}.$$ 
	Fix $i\in \{1,2,\cdots,n\}$. Since $I_{(\lambda_{j}^{x_{0}}=1)}\lambda_{i}^{x_{0}}=0$ for any $j\in \{1,2,\cdots,n\}\setminus \{i\}$, it follows from  \eqref{eq.KKM1} that
	\begin{equation}\label{eq.KKM2}
		I_{(\lambda_{j}^{x_{0}}=1)}d(x_{0},F_{i})=0,~\forall j\in \{1,2,\cdots,n\}\setminus \{i\}.
	\end{equation}
	Besides, since $$I_{(\lambda_{i}^{x_{0}}=1)}x_{0}=I_{(\lambda_{i}^{x_{0}}=1)}\sum_{l=1}^{n}\lambda_{l}^{x_{0}}e_{l}=I_{(\lambda_{i}^{x_{0}}=1)}\lambda_{i}^{x_{0}}e_{i}=I_{(\lambda_{i}^{x_{0}}=1)}e_{i}\in I_{(\lambda_{i}^{x_{0}}=1)}F_{i},$$ 
	we have
	\begin{equation}\label{eq.KKM3}
		I_{(\lambda_{i}^{x_{0}}=1)}d(x_{0},F_{i})=0.
	\end{equation}
	Combining \eqref{eq.KKM2} with \eqref{eq.KKM3}, we obtain $I_{a_{1}}d(x_{0},F_{i})=0$.

	\par
	Second, we prove that 
	$$I_{a_{k}}d(x_{0},F_{i})=0,~\forall i\in \{1,2,\cdots,n\},~k\in \{2,3,\cdots,n\}.$$ 
	Fix $k\in \{2,3,\cdots,n\}$. If there exists some $i_{1}\in \{1,2,\cdots,n\}$ such that $I_{a_{k}}d(x_{0},F_{i_{1}})\neq 0$, then there exist $L\in \mathcal{L}_{k}$ and $b\in B_{\mathcal{F}}$ with 
	$$0<b\leq [\bigwedge_{l\in L}(0<\lambda_{l}^{x_{0}}< 1)]\wedge[\bigwedge_{i\in \{1,2,\cdots,n\}\setminus L}(\lambda_{i}^{x_{0}}=0)]$$ 
	such that $I_{b}d(x_{0},F_{i_{1}})>0$ on $b$. Together with \eqref{eq.KKM1}, these relations imply
	that $i_{1}\in L$ and 
	\begin{equation}\label{eq.KKM4} 
		I_{b}d(x_{0},F_{i})>0~\text{on}~b,~\forall i\in L, 
	\end{equation}
	whereas 
	$$I_{b}d(x_{0},F_{i})=0,~\forall i\in \{1,2,\cdots,n\}\setminus L.$$ 
	However, since
	\begin{align}
		I_{b}x_{0}&=I_{b}\sum_{i=1}^{n}\lambda_{i}^{x_{0}}e_{i}\nonumber\\
		&=I_{b}\sum_{i\in L} \lambda_{i}^{x_{0}}e_{i}\nonumber\\
		&=I_{b}\left[(I_{b}\lambda_{i_{1}}^{x_{0}}+I_{b^{c}})e_{i_{1}}+\sum_{i\in L\setminus\{i_{1}\}} I_{b}\lambda_{i}^{x_{0}}e_{i}\right]\nonumber\\
		&\in I_{b} \operatorname{Conv}_{L^{0}}(\{e_{i}\mid i\in L\})\nonumber\\
		&\subseteq I_{b} \left(\bigcup\sigma\{F_{i}\mid i\in L\}\right),\nonumber	
	\end{align}
	there exist $j\in L$ and $c\in B_{\mathcal{F}}$ with $0<c\leq b$ such that $I_{c}x_{0} \in I_{c}F_{j}$. This implies that $I_{c} d(x_{0},F_{j})=0$, which contradicts \eqref{eq.KKM4}.

	\par
	To conclude, since $\{a_{k}, k=1,2,\cdots,n\}\in p(1)$, we have 
	$$d(x_{0},F_{i})=0,~\forall i\in \{1,2,\cdots,n\}.$$ 
	For each $i\in\{1,\dots,n\}$, since $F_{i}$ is $\sigma$-stable and closed, \cite[Theorem 2.12]{GWYZ2020} yields $x_{0}\in F_{i}$. Consequently,  $x_{0}\in\bigcap_{i=1}^{n}F_{i}$. Finally, $\bigcap_{i=1}^{n}F_{i}$ is a $\sigma$-stable closed subset of the stably compact set $\Delta_{n-1}$, and hence stably compact.
\end{proof}

\begin{proof}[Proof of Theorem \ref{thm.KKM2}]
	Let $\{g_{1}, g_{2}, \dots, g_{n}\}$ be an arbitrary finite subset of $G$. We show that 
	$$\bigcap_{i=1}^{n}T(g_{i})\neq \emptyset.$$ 
	Set $H=\operatorname{Conv}_{L^{0}}(\{g_{1},g_{2},\cdots,g_{n}\})$ and define a mapping $\varphi:\Delta_{n-1}\rightarrow H$ by 
	$$\varphi(x)=\sum_{i=1}^{n}\lambda_{i}^{x}g_{i},~\forall x\in \Delta_{n-1}.$$ 
	The mapping $\varphi$ is $\sigma$-stable and continuous from $(\Delta_{n-1},\mathcal{T}_{c}(\|\cdot\|))$ to $(H,\mathcal{T}|_{H})$. Since $T$ is a random KKM mapping, $H\subseteq T(G)$ and $g_i\in T(g_i)$ for each $i=1,\ldots,n$. Thus $T(g_i)\cap H$ is nonempty, $\sigma$-stable and $\mathcal{T}|_{H}$-closed  because $T(g_i)$ is $\mathcal{T}|_{T(G)}$-closed. Consequently,
	$$\{\varphi^{-1}(T(g_{i})\cap H)\mid i=1,2,\dots,n\}$$
	is a family of nonempty $\sigma$-stable closed subsets of $\Delta_{n-1}$.
	
	\par 
	For any nonempty $I\subseteq \{1,2,\cdots,n\}$ and any $x\in \Delta_{n-1}$ with $\{\lambda_{i}^{x}\mid i\in I\}\subseteq L^{0}(\mathcal{F},[0,1])$ and $\sum_{i\in I}\lambda_{i}^{x}=1$, since $T$ is a random KKM mapping, we have
	\begin{align}
		\varphi\left(\sum_{i\in I}\lambda_{i}^{x}e_{i}\right)&=\sum_{i\in I}\lambda_{i}^{x}g_{i}\nonumber\\
		&\in \operatorname{Conv}_{L^{0}}(\{g_{i}\mid i\in I\})\nonumber\\
		&\subseteq\left(\bigcup\{T(g)\mid g\in \sigma\{g_{i}\mid i\in I\}\}\right)\cap H\nonumber\\
		&=\bigcup\{T(g)\cap H\mid g\in \sigma\{g_{i}\mid i\in I\}\}.\nonumber
	\end{align}		
	Furthermore, since $\varphi$ is $\sigma$-stable and surjective,  \cite[Lemma 2.9]{TMG2024} yields
	\begin{align}
		\sum_{i\in I}\lambda_{i}^{x}e_{i}&\in \varphi^{-1}\left(\bigcup\{T(g)\cap H\mid g\in \sigma\{g_{i}\mid i\in I\}\}\right)\nonumber\\
		&=\bigcup\{\varphi^{-1}(T(g)\cap H)\mid g\in \sigma\{g_{i}\mid i\in I\}\}\nonumber\\
		&=\bigcup\left\{\sum_{i\in I}I_{a_{i}}\varphi^{-1}(T(g_{i})\cap H)\mid \{a_{i},i\in I\}\in p(1)\right\}\nonumber\\
		&=\bigcup\sigma\{\varphi^{-1}(T(g_{i})\cap H)\mid i\in I\}\nonumber\\
		&=\bigcup\sigma\{\varphi^{-1}(T(g_{i}))\mid i\in I\}.\nonumber
	\end{align}	
	Here the last equality follows from $\varphi(\Delta_{n-1})=H$. Therefore, 
	$$\operatorname{Conv}_{L^{0}}(\{e_{i}\mid i\in I\})\subseteq \bigcup\sigma\{\varphi^{-1}(T(g_{i}))\mid i\in I\}.$$ 
	By Theorem \ref{thm.KKM1}, there exists $x_{0}\in\bigcap_{i=1}^{n}\varphi^{-1}(T(g_{i}))$. Hence $\varphi(x_{0})\in\bigcap_{i=1}^{n}T(g_{i})$. Since $\{g_{1},\ldots,g_{n}\}$ was arbitrary, the family $\{T(g)\mid g\in G\}$ has the finite intersection property.
 
	\par 
	Now suppose that $T(g_{0})$ is stably compact for some $g_{0}\in G$. By the finite intersection property established above, $\{T(g)\cap T(g_{0})\mid g\in G\}$ is a family of nonempty $\sigma$-stable closed subsets of $T(g_{0})$ with the finite intersection property. By (3) of Proposition \ref{prop.stable compactness},
	$$\bigcap \{T(g)\cap T(g_{0})\mid g\in G\}=\bigcap \{T(g)\mid g\in G\}\neq \emptyset.$$
\end{proof}

\begin{proof}[Proof of Lemma \ref{lemm.example KKM}]
	Since $x\in T(x)$ for any $x\in G$, $T$ has nonempty values and $T(G)=G$. Hence, $T(G)$ is $\sigma$-stable.
	
	\par 
	We first prove that $T$ is $\sigma$-stable. Let $\{x_{n},n\in \mathbb{N}\}$ be a sequence in $G$ and $\{a_{n},n\in \mathbb{N}\}\in p(1)$. We show that 
	$$\sum_{n=1}^{\infty}I_{a_{n}}T(x_{n})= T\left(\sum_{n=1}^{\infty}I_{a_{n}}x_{n}\right).$$ 
	On the one hand, for any $y:=\sum_{n=1}^{\infty}I_{a_{n}}y_{n}\in \sum_{n=1}^{\infty}I_{a_{n}}T(x_{n})$ with $y_{n}\in T(x_{n})$ for each $n\in \mathbb{N}$, by the $\sigma$-stability of $H$, 
	\begin{align}
		H(y,y)&=\sum_{n=1}^{\infty}I_{a_{n}}H(y_{n},y_{n})\nonumber\\
		&\leq \sum_{n=1}^{\infty}I_{a_{n}}H(y_{n},x_{n})\nonumber\\
		&=H\left(y,\sum_{n=1}^{\infty}I_{a_{n}}x_{n}\right),\nonumber	
	\end{align}
	implying $y\in T\left(\sum_{n=1}^{\infty}I_{a_{n}}x_{n}\right)$. Then we have 
	$$\sum_{n=1}^{\infty}I_{a_{n}}T(x_{n})\subseteq T\left(\sum_{n=1}^{\infty}I_{a_{n}}x_{n}\right).$$ 
	On the other hand, for any chosen $z\in T\left(\sum_{n=1}^{\infty}I_{a_{n}}x_{n}\right)$, let $z_{n}=I_{a_{n}} z+I_{a_{n}^{c}}x_{n}$ for each $n\in \mathbb{N}$. For any $n\in \mathbb{N}$, since
	\begin{align}
		H(z,z)&\leq H\left(z,\sum_{n=1}^{\infty}I_{a_{n}}x_{n}\right)\nonumber\\
		&=\sum_{n=1}^{\infty}I_{a_{n}}H(z,x_{n}),\nonumber		
	\end{align}
	it follows that
	\begin{align}
		H(z_{n},z_{n})&=I_{a_{n}}H(z,z)+I_{a_{n}^{c}}H(x_{n},x_{n})\nonumber\\
		&\leq I_{a_{n}}\left(\sum_{k=1}^{\infty}I_{a_{k}}H(z,x_{k})\right)+I_{a_{n}^{c}}H(x_{n},x_{n})\nonumber\\
		&=I_{a_{n}}H(z,x_{n})+I_{a_{n}^{c}}H(x_{n},x_{n})\nonumber\\
		&=H(z_{n},x_{n}),\nonumber
	\end{align} 
	implying $z_{n}\in T(x_{n})$. Consequently, 
	$$z=\sum_{n=1}^{\infty}I_{a_{n}}z=\sum_{n=1}^{\infty}I_{a_{n}}z_{n}\in \sum_{n=1}^{\infty}I_{a_{n}}T(x_{n}),$$ 
	that is, 
	$$T\left(\sum_{n=1}^{\infty}I_{a_{n}}x_{n}\right)\subseteq\sum_{n=1}^{\infty}I_{a_{n}}T(x_{n}).$$
	Therefore,
	$$T\left(\sum_{n=1}^{\infty}I_{a_{n}}x_{n}\right)=\sum_{n=1}^{\infty}I_{a_{n}}T(x_{n}).$$
	
	\par 
	Second, we prove that
	$$\operatorname{Conv}_{L^{0}}\{g_{1},\cdots,g_{n}\}\subseteq\bigcup\{T(g)\mid g\in \sigma\{g_{1},\cdots,g_{n}\}\}$$ 
	for any finite subset $\{g_{1},\cdots,g_{n}\}$ of $G$.  For this, let $y=\sum_{i=1}^{n}\lambda_{i} g_{i}\in \operatorname{Conv}_{L^{0}}\{g_{1},\ldots,g_{n}\}$, where $\lambda_{i}\in L^{0}_{+}(\mathcal{F})$ for each $i=1,\ldots,n$ and $\sum_{i=1}^{n}\lambda_{i}=1$. Let $a_{i}=( H(y,y)\leq H(y,g_{i})) $ for each $i\in \{1,2,\cdots,n\}$. 
	By the $L^{0}$-quasiconvexity of $H(y,\cdot)$ and Remark \ref{rmk.quasiconvexity}, we have 
	$$H(y,y)= H\left(y,\sum_{i=1}^{n}\lambda_{i}g_{i}\right)\leq\bigvee_{i=1}^{n}H(y,g_{i}),$$
	which implies $\bigvee_{i=1}^{n}a_{i}=1$. Now, let $b_{1}=a_{1}$ and
	$$b_{i}=a_{i}\wedge \left(\bigvee_{j=1}^{i-1}a_{j}\right)^{c},\qquad \forall i=2,\ldots,n.$$
	Then $\{b_{i},i=1,\ldots,n\}\in p(1)$ and $b_{i}\leq a_{i}$ for each $i\in \{1,2,\cdots,n\}$. Let $g=\sum_{i=1}^{n}I_{b_i}g_i$. Then $g\in\sigma\{g_1,\ldots,g_n\}$. By the $\sigma$-stability of $H$, we have 
	$$H(y,g)=H\left(y,\sum_{i=1}^{n}I_{b_{i}}g_{i}\right)=\sum_{i=1}^{n}I_{b_{i}}H(y,g_{i})\geq\sum_{i=1}^{n}I_{b_{i}}H(y,y)=H(y,y).$$
	Thus, $y\in T\left(\sum_{i=1}^{n}I_{b_{i}}g_{i}\right)\subseteq \bigcup\{T(g)\mid g\in \sigma\{g_{1},\cdots,g_{n}\}\}$. Since $y$ was arbitrary, we conclude that 
	$$\operatorname{Conv}_{L^{0}}\{g_{1},\ldots,g_{n}\}\subseteq\bigcup\{T(g)\mid g\in \sigma\{g_{1},\cdots,g_{n}\}\}.$$
\end{proof}

\par 
The following corollary allows the random KKM mapping principle to be applied through a set-valued enlargement.

\begin{corollary}\label{coro.generalization}
	Let $(E,\mathcal{T})$ be a Hausdorff topological $L^{0}(\mathcal{F},\mathbb{K})$-module with base $(\Omega,\mathcal{F},P)$, $G$ a nonempty $\sigma$-stable subset of $E$, and $f:G\to 2^{E}\setminus\{\emptyset\}$ a random KKM mapping. Suppose that there exists a $\sigma$-stable mapping $g:G\to 2^{E}\setminus\{\emptyset\}$ with $\sigma$-stable closed values such that
	$$f(x)\subseteq g(x),~\forall x\in G,$$
	and
	$$\bigcap_{x\in G}f(x)=\bigcap_{x\in G}g(x).$$
	If $g(x_{0})$ is stably compact for some $x_{0}\in G$, then
	$$\bigcap_{x\in G}f(x)\neq\emptyset.$$
\end{corollary}
\begin{proof}
	Since $f(x)\subseteq g(x)$ for every $x\in G$ and $g$ is $\sigma$-stable, $g$ is also a random KKM
	mapping. By Theorem \ref{thm.KKM2},
	$$\emptyset\neq\bigcap_{x\in G}g(x)=\bigcap_{x\in G}f(x).$$
\end{proof}

\par 
The random KKM mapping principle also yields the following best approximation result in $RN$ modules.

\begin{corollary}
	Let $(E,\|\cdot\|)$ be an $RN$ module over $\mathbb{K}$ with base $(\Omega,\mathcal{F},P)$, $G$ a nonempty stably compact $L^{0}$-convex subset of $E$, and $f:G\rightarrow E$ a $\sigma$-stable $\mathcal{T}_{c}(\|\cdot\|)$-continuous mapping. Then there exists $y_{0}\in G$ such that $$\|y_{0}-f(y_{0})\|=\min_{x\in G}\|x-f(y_{0})\|. $$
\end{corollary}
\begin{proof}
	Define a mapping $T:G\to 2^{G}\setminus\{\emptyset\}$ by
	$$T(x)=\{y\in G \mid \|y-f(y)\|\leq\|x-f(y)\|\},~\forall x\in G.$$
	By Lemma \ref{lemm.example KKM}, $T$ is a random KKM mapping. Since $f$ is $\sigma$-stable and
	$\mathcal{T}_{c}$-continuous, each $T(x)$ is a $\sigma$-stable closed subset of $G$, and hence stably compact. By Theorem \ref{thm.KKM2}, $\bigcap_{x\in G}T(x)\neq\emptyset$. For any $y_{0}\in \bigcap_{x\in G}T(x)$, we have 
	$$\|y_{0}-f(y_{0})\|\leq\|x-f(y_{0})\|,~\forall x\in G$$
	and the desired equality follows.
\end{proof}

\subsection{Proofs of the results in Subsection~\ref{sec.2.3}}\label{sec.3.2}

We first collect the localization properties needed for the random Fan--Browder fixed point theorem. The coincidence theorem is then obtained by applying this theorem to a set-valued mapping on $X\times Y$.

\par
We first recall Lemma~\ref{lemm.hereditarily disjoint degree}, which will be frequently used in the sequel.

\begin{lemma}[\cite{GMPTY2026}]\label{lemm.hereditarily disjoint degree}
	Let $E$ be an $L^{0}$-module and $M$ and $G$ be two $\sigma$-stable subsets of $E$ such that $M\cap G=\emptyset$. Then there exists a unique $b\in B_{\mathcal{F}}$ with $b>0$ such that the following are satisfied:
	\begin{enumerate}[label=(\arabic*)]
		\item $I_{a}M \cap I_{a}G=\emptyset$ for any $a\in B_{\mathcal{F}}$ with $0<a\leq b$;
		\item $I_{a}M \cap I_{a}G \neq \emptyset$ for any $a\in B_{\mathcal{F}}$ with $a\leq b^{c}$.
	\end{enumerate}
\end{lemma}

\par 
In \cite[Definition 3.7]{GMPTY2026}, $b$ as guaranteed in Lemma \ref{lemm.hereditarily disjoint degree} is called the hereditarily disjoint degree between $M$ and $G$, denoted by $h(M,G)$.  When $M=\{x\}$, we simply denote $h(\{ x\}, G)$ by $h(x,G)$. Finally, we make the convention that $h(M,G)=0$ if $M$ and $G$ are two $\sigma$-stable subsets of $E$ such that $M\cap G\neq\emptyset$.

\par 
The following lemma will be used in the proofs of Theorem~\ref{thm.Fan fixed point} and Lemma~\ref{lemm.two continuty}.

\begin{lemma}\label{lemm.fundament1}
	Let $G$ and $H$ be two $\sigma$-stable subsets of an $L^{0}(\mathcal{F},\mathbb{K})$-module and $a\in B_{\mathcal{F}}$ with $a>0$. If $I_{a'}G\nsubseteq I_{a'}H$ for any $0<a'\leq a$, then there exists $g\in G$ such that $I_{a'}g\notin I_{a'}H$ for any $0<a'\leq a$. 
\end{lemma}
\begin{proof}
	If not, then for any $g\in G$, there exists $0<a'_{g}\leq a$ such that $I_{a'_{g}}g\in I_{a'_{g}}H$, this implies that $(h(g,H))^{c}\wedge a>0$. For each $g\in G$, let $a_{g}=(h(g,H))^{c}\wedge a$. Then $I_{a_{g}}g\in I_{a_{g}}H$ and
	\begin{equation}\label{eq.lemm.fundament1}
		I_{a''_{g}}g\notin I_{a''_{g}}H,~\forall 0<a''_{g}\leq a\wedge a_{g}^{c}.
	\end{equation}
	We claim that $\bigwedge_{g\in G}a_{g}=0$. Otherwise, suppose $\bigwedge_{g\in G}a_{g}>0$. Then, for each $g\in G$, $I_{\bigwedge_{g\in G}a_{g}}g\in I_{\bigwedge_{g\in G}a_{g}}H$, and hence $I_{\bigwedge_{g\in G}a_{g}}G\subseteq I_{\bigwedge_{g\in G}a_{g}}H$, which contradicts the assumption. Therefore, we must have $\bigwedge_{g\in G}a_{g}=0$, equivalently $\bigvee_{g\in G}(a\wedge a_{g}^{c})=a$. Then, there exists a sequence $\{g_{n}, n\in\mathbb{N}\}$ in $G$ such that $\bigvee_{n=1}^{\infty}(a\wedge a_{g_{n}}^{c})=a$. Let $b_{1}=a\wedge a_{g_{1}}^{c}$ and 
	$$b_{n}=(a\wedge a_{g_{n}}^{c})\wedge[\bigvee_{i=1}^{n-1}(a\wedge a_{g_{i}}^{c})]^{c},~\forall n\geq 2.$$ 
	Then $\{a^{c}\}\cup\{b_{n},n\in \mathbb{N}\}\in p(1)$. Arbitrarily choose a $g_{0}\in G$ and let $g_{*}=\sum_{n=1}^{\infty}I_{b_{n}}g_{n}+I_{a^{c}}g_{0}$, then $g_{*}\in G$. By \eqref{eq.lemm.fundament1}, we have $$I_{a'}g_{*}=\sum_{n=1}^{\infty}I_{b_{n}\wedge a'}g_{n}\notin I_{a'}H,~\forall 0<a'\leq a,$$ 
	which contradicts our assumption.
\end{proof}

\begin{lemma}\label{lemm.stable1}
	Let $E_{1}$ and $E_{2}$ be two $L^{0}(\mathcal{F},\mathbb{K})$-modules, $X\subseteq E_{1}, Y\subseteq E_{2}$  two nonempty sets and  $T:X\rightarrow 2^{Y}\setminus \{\emptyset\}$ a $\sigma$-stable set-valued mapping. Then the following statements hold:
	\begin{itemize}[widest={(1)}, leftmargin=*]
		\item [(1)] $T^{-1}(y)$ is $\sigma$-stable for any $y\in T(X)$.
		\item [(2)] For any $\{a_{n},n\in \mathbb{N}\}\in p(1)$ and any sequence $\{y_{n},n\in \mathbb{N}\}$ in $T(X)$, 
		$$\sum_{n=1}^{\infty}I_{a_{n}}T^{-1}(y_{n})=T^{-1}\left(\sum_{n=1}^{\infty}I_{a_{n}} y_{n}\right).$$
		\item [(3)] For any $a\in B_{\mathcal{F}}$ and any $x\in X$, if there exists a finitely stable subset $G$ of $X$ such that $I_{a}x\in I_{a}G$, then $I_{a}T(x)\subseteq I_{a}T(G)$.
		\item [(4)] For any $a\in B_{\mathcal{F}}$, $x\in X$ and any $y\in T(X)$,  $I_{a}y\in I_{a}T(x)$ if and only if $I_{a} x\in I_{a} T^{-1}(y)$.
	\end{itemize}
\end{lemma}
\begin{proof}
	(1) Let $y\in T(X)$. Then $T^{-1}(y)\neq\emptyset$. For any $\{b_{m},m\in \mathbb{N}\}\in p(1)$ and any sequence $\{x_{m},m\in \mathbb{N}\}$ in $T^{-1}(y)$, we have
	$$y=\sum_{m=1}^{\infty}I_{b_{m}}y
	\in \sum_{m=1}^{\infty}I_{b_{m}}T(x_{m})
	=T\left(\sum_{m=1}^{\infty}I_{b_{m}}x_{m}\right),$$
	which implies that $\sum_{m=1}^{\infty}I_{b_{m}}x_{m}\in T^{-1}(y)$. Hence,  $T^{-1}(y)$ is $\sigma$-stable.
	
	\par 
	(2) Let $x=\sum_{n=1}^{\infty} I_{a_{n}}x_{n}$ with $x_{n}\in T^{-1}(y_{n})$ for each $n\in\mathbb{N}$. Then
	$$\sum_{n=1}^{\infty}I_{a_{n}}y_{n} \in\sum_{n=1}^{\infty}I_{a_{n}}T(x_{n})=T\left(\sum_{n=1}^{\infty}I_{a_{n}}x_{n}\right)=T(x),$$
	which implies that $x\in T^{-1}\left(\sum_{n=1}^{\infty} I_{a_{n}}y_{n}\right)$. Hence, $$\sum_{n=1}^{\infty}I_{a_{n}}T^{-1}(y_{n})\subseteq T^{-1}\left(\sum_{n=1}^{\infty}I_{a_{n}} y_{n}\right).$$
	
	\par 
	Conversely, let $x\in T^{-1}\left(\sum_{n=1}^{\infty} I_{a_{n}}y_{n}\right)$. For each $n\in\mathbb{N}$, choose an arbitrary $x_{n}\in T^{-1}(y_{n})$ and let $z_{n}=I_{a_{n}} x +I_{a_{n}^{c}}x_{n}$. Then
	$$y_{n}=I_{a_{n}} y_{n} +I_{a_{n}^{c}}y_{n}\in I_{a_{n}} T(x) +I_{a_{n}^{c}}T(x_{n})=T(I_{a_{n}} x +I_{a_{n}^{c}}x_{n})=T(z_{n}),$$
	so that $z_{n}\in T^{-1}(y_{n})$ for each $n\in \mathbb{N}$. Therefore, $$x=\sum_{n=1}^{\infty}I_{a_{n}} z_{n}\in \sum_{n=1}^{\infty}I_{a_{n}} T^{-1}(y_{n}),$$ which yields $$T^{-1}\left(\sum_{n=1}^{\infty}I_{a_{n}} y_{n}\right)\subseteq\sum_{n=1}^{\infty}I_{a_{n}}T^{-1}(y_{n}).$$
	
	\par 
	(3) Arbitrarily choose $z\in G$ and let $x_{1}=I_{a}x+I_{a^{c}}z$. Then $x_{1}\in G$, and we have $$I_{a}T(x)+I_{a^{c}}T(z)=T(I_{a}x+I_{a^{c}}z)=T(x_{1})\subseteq T(G),$$ implying $I_{a}T(x)\subseteq I_{a}T(G)$.
	
	\par 
	(4) If $I_{a}y\in I_{a}T(x)$, take an arbitrary $y_{0}\in T(x)$ and let $z=I_{a}y+I_{a^{c}}y_{0}$. Then $x\in T^{-1}(z)=I_{a}T^{-1}(y)+I_{a^{c}}T^{-1}(y_{0})$, which implies that $I_{a}x\in I_{a} T^{-1}(y)$.
	
	\par 
	If $I_{a}x\in I_{a} T^{-1}(y)$,  take an arbitrary $x_{0}\in T^{-1}(y)$ and let $x_{1}=I_{a}x+I_{a^{c}}x_{0}$. Then $y\in T(x_{1})=I_{a}T(x)+I_{a^{c}}T(x_{0})$, which implies that $I_{a}y\in I_{a} T(x)$.
	
\end{proof}

\par 
For $L^{0}(\mathcal{F},\mathbb{K})$-modules, we introduce the following localized version of the hereditary expansion defined in \cite[Definition 2.15]{GMPTY2026}.

\begin{definition}
	Let $E$ be an $L^{0}(\mathcal{F},\mathbb{K})$-module, $G$ a subset of $E$ and $a\in B_{\mathcal{F}}$ with $a>0$. Then $$ \llbracket G \rrbracket_{a}=\{y\in E:~\text{there~exists}~0<a'\leq a~\text{such~that}~I_{a'}y\in I_{a'}G\}$$ is called the hereditary $a$-expansion of $G$.
\end{definition}

\par 
It is easy to check that $I_{a}G\subseteq \llbracket G \rrbracket_{a}$, $\llbracket G \rrbracket_{a}\subseteq \llbracket G \rrbracket_{b} $ whenever $0<a\leq b$, and $\llbracket G \rrbracket_{a}=\emptyset$ iff $G=\emptyset$. When $a=1$, we write $\llbracket G\rrbracket$ for $\llbracket G\rrbracket_{1}$. If $E$ is $\sigma$-stable, $\llbracket G\rrbracket$ coincides with the hereditary expansion defined in \cite[Definition 2.15]{GMPTY2026} with $B=B_{\mathcal{F}}$.

\par
The properties in Lemma \ref{lemm.stable2} below are localized counterparts of the properties of hereditary expansion and hereditary complements established in \cite[Theorems 2.16 and 2.29]{GMPTY2026}. For the reader's convenience, we give a direct proof.

\begin{lemma}\label{lemm.stable2}
	Let $E$ be an $L^{0}(\mathcal{F},\mathbb{K})$-module, $G$ and $H$ two nonempty subsets of $E$ such that $H$ is $\sigma$-stable and $a\in B_{\mathcal{F}}$ with $a>0$. Then we have the following statements:
	\begin{itemize}[widest={(1)}, leftmargin=*]
		\item [(1)] $H\setminus \llbracket G \rrbracket_{a}$ is either empty or $\sigma$-stable. 
		\item [(2)] $\sum_{n=1}^{\infty}I_{a_{n}}(H\setminus \llbracket G_{n} \rrbracket_{a})=H\setminus \llbracket \sum_{n=1}^{\infty}I_{a_{n}}G_{n} \rrbracket_{a}$ for any $\{a_{n},n\in \mathbb{N}\}\in p(1)$ and any sequence $\{G_{n},n\in \mathbb{N}\}$ of nonempty subsets of $H$ with $H\setminus \llbracket G_{n} \rrbracket_{a}\neq \emptyset$ for each $n\in \mathbb{N}$.
		\item [(3)] If $(E,\mathcal{T})$ is a Hausdorff topological space, $(H,\mathcal{T}|_{H})$ is a stable topological space and $G\subseteq H$ is open in $H$, then $H\cap\llbracket G\rrbracket_{a}$ is open in $H$.
		\item [(4)] If $(E,\mathcal{T})$ is a Hausdorff topological space such that $(H,\mathcal{T}|_{H})$ is a stable topological space, and $F$ is a nonempty $\sigma$-stable closed subset of $H$, then $H\cap\llbracket F\rrbracket_{a}$ is closed in $H$.
	\end{itemize}
\end{lemma}
\begin{proof}
	(1) Assume that $H\setminus \llbracket G \rrbracket_{a}\neq \emptyset$. For any $\{a_{n},n\in \mathbb{N}\}\in p(1)$ and any sequence $\{x_{n},n\in \mathbb{N}\}$ in $H\setminus \llbracket G \rrbracket_{a}$, if $\sum_{n=1}^{\infty}I_{a_{n}}x_{n}\in \llbracket G \rrbracket_{a}$, there exists $0<a'\leq a$ such that $$I_{a'} \left(\sum_{n=1}^{\infty}I_{a_{n}}x_{n}\right) \in I_{a'} G.$$ Then, there exists $a_{n}$ with $a'\wedge a_{n}>0$ such that $I_{a'\wedge a_{n}}x_{n}\in I_{a'\wedge a_{n}} G $, which contradicts $x_{n}\in H\setminus \llbracket G \rrbracket_{a}$. Therefore, 
	$$\sum_{n=1}^{\infty}I_{a_{n}}x_{n}\in H\setminus \llbracket G \rrbracket_{a},$$ 
	and so $H\setminus \llbracket G \rrbracket_{a}$ is $\sigma$-stable.
	
	\par 
	(2) First, for any $x=\sum_{n=1}^{\infty}I_{a_{n}}x_{n}\in \sum_{n=1}^{\infty}I_{a_{n}}(H\setminus \llbracket G_{n} \rrbracket_{a})$ with $x_{n}\in H\setminus \llbracket G_{n} \rrbracket_{a}$ for each $n\in \mathbb{N}$, we have $$I_{a'}\left(\sum_{n=1}^{\infty}I_{a_{n}}x_{n}\right)\notin I_{a'}\left(\sum_{n=1}^{\infty}I_{a_{n}}G_{n}\right),~\forall 0<a'\leq a,$$ which implies that $x\in H\setminus \llbracket \sum_{n=1}^{\infty}I_{a_{n}}G_{n} \rrbracket_{a}$. Then 
	$$\sum_{n=1}^{\infty}I_{a_{n}}(H\setminus \llbracket G_{n} \rrbracket_{a})\subseteq H\setminus \llbracket \sum_{n=1}^{\infty}I_{a_{n}}G_{n} \rrbracket_{a}.$$
	
	\par 
	On the other hand, for any $x\in H\setminus \llbracket \sum_{n=1}^{\infty}I_{a_{n}}G_{n} \rrbracket_{a}$, we have $$I_{a'}x\notin I_{a'}\left(\sum_{n=1}^{\infty}I_{a_{n}}G_{n}\right),~\forall 0<a'\leq a,$$ 
	which implies that 
	$$I_{a'}x\notin I_{a'}G_{n},~\forall 0<a'\leq a\wedge a_{n}$$ 
	for each  $n\in \mathbb{N}$. Take an arbitrary $y_{n}\in H\setminus \llbracket G_{n} \rrbracket_{a}$ and let $z_{n}=I_{a_{n}}x+I_{a_{n}^{c}}y_{n}$ for each $n\in \mathbb{N}$. Then 
	$$I_{a'}z_{n}\notin I_{a'}G_{n},~\forall 0<a'\leq a,$$ 
	which implies that $z_{n}\in H\setminus \llbracket G_{n} \rrbracket_{a} $ and $I_{a_{n}}x=I_{a_{n}}z_{n}\in I_{a_{n}}(H\setminus \llbracket G_{n} \rrbracket_{a} )$ for each $n\in \mathbb{N}$. Then $x\in \sum_{n=1}^{\infty}I_{a_{n}}(H\setminus \llbracket G_{n} \rrbracket_{a})$ and so 
	$$H\setminus \llbracket \sum_{n=1}^{\infty}I_{a_{n}}G_{n} \rrbracket_{a}\subseteq \sum_{n=1}^{\infty}I_{a_{n}}(H\setminus \llbracket G_{n} \rrbracket_{a}). $$
	
	\par 
	(3) Since $H$ is $\sigma$-stable,
	$$H\cap\llbracket G\rrbracket_{a}=\bigcup_{0<c\leq a}\ \bigcup_{u\in G}\{x\in H\mid I_{c}x+I_{c^{c}}u\in G\}.$$
	For any fixed $0<c\leq a$ and $u\in G$, define $\psi:H\to H$ by 
	$$\psi(x)=I_cx+I_{c^c}u,~\forall x\in H.$$
	Let $\mathcal{B}$ be a $\sigma$-bistable base of $\mathcal{T}|_H$. Fix $x\in H$. For any $O\in \mathcal{B}$ with $\psi(x)\in O$, choose $V\in \mathcal{B}$ with $x\in V$. Then $I_{c}O+I_{c^{c}}V\in \mathcal{B}$ and $x\in I_{c}O+I_{c^{c}}V$. Since $I_{c}\psi(y)=I_{c}y$ and $I_{c^c}\psi(y)=I_{c^c}u$ for any $y\in H$, we have
	$$\psi(I_{c}O+I_{c^{c}}V)=I_{c}O+I_{c^{c}}u=I_{c}O+I_{c^{c}}\psi(x)\subseteq I_{c}O+I_{c^{c}}O=O.$$
	Hence $\psi$ is continuous, and $\psi^{-1}(G)$ is open in $H$. Consequently, $H\cap\llbracket G\rrbracket_a$ is open in $H$.
	
	\par 
	(4) Let $\mathcal{B}$ be a $\sigma$-bistable base of $\mathcal{T}|_{H}$ and fix $x\in cl_{\mathcal{T}|_{H}}(H\cap\llbracket F\rrbracket_{a})$, where $cl_{\mathcal{T}|_{H}}(H\cap\llbracket F\rrbracket_{a})$ is the closure of $H\cap\llbracket F\rrbracket_{a}$ under $\mathcal{T}|_{H}$. For each $O\in\mathcal{B}$ containing $x$, define
	$$\mathcal{A}_{O}=\{0<c\leq a\mid I_{c}O\cap I_{c}F\neq\emptyset\}$$
	and let $a_{O}:=\bigvee\mathcal{A}_{O}$. Then $\mathcal{A}_{O}\neq \emptyset$. Since $O$ and $F$ are $\sigma$-stable, $a_{O}$ is attained. We claim that $$d:=\bigwedge \{a_{O} \mid x\in O,~O\in\mathcal{B}\}>0.$$
	Otherwise, $d=0$, namely, $\bigvee \{a_{O}^{c} \mid x\in O,~O\in\mathcal{B}\}=1$. Then, there exists a sequence $\{O_{n},n\in \mathbb{N}\}$ in $\mathcal{B}$ such that $x\in O_{n}$ for each $n\in \mathbb{N}$ and $\bigvee_{n=1}^{\infty}a_{O_{n}}^{c}=1$. Let $b_{1}=a_{O_{1}}^{c}$ and 
	$$b_{i}=a_{O_{i}}^{c}\wedge (\bigvee_{j=1}^{i-1}a_{O_{j}}^{c})^{c},~\forall i\geq 2.$$ 
	Then $\{b_{n}:n\in\mathbb{N}\}\in p(1)$ and  $b_{n}\leq a_{O_{n}}^{c}$ for each $n\in \mathbb{N}$. Let $V=\sum_{n=1}^{\infty}I_{b_{n}}O_{n}$. By the $\sigma$-bistability of $\mathcal{B}$, we have $V\in\mathcal{B}$ and $x\in V$. Suppose that $I_{c}V\cap I_{c}F\neq\emptyset$ for some $0<c\leq a$. Then, there exists $n\in \mathbb{N}$ such that $c\wedge b_{n}>0$ and 
	$$I_{c\wedge b_{n}}V\cap  I_{c\wedge b_{n}}F=I_{c\wedge b_{n}}O_{n}\cap  I_{c\wedge b_{n}}F\neq\emptyset,$$
	which implies that $c\wedge b_{n}\in \mathcal{A}_{O_{n}}$. This contradicts the fact that $0<c\wedge b_{n}\leq a_{O_{n}}^{c}$. Thus $(\sum_{n=1}^{\infty}I_{b_{n}}O_{n})\cap\llbracket F\rrbracket_{a}=\emptyset$, contrary to the choice of $x$.
	
	\par 
	Choose $y\in F$ and $O_{x}\in\mathcal{B}$ with $x\in O_{x}$, and set $z=I_{d}x+I_{d^{c}}y$. Given $O_{z}\in\mathcal{B}$ with $z\in O_{z}$, let $U=I_{d}O_{z}+I_{d^{c}}O_{x}$. Then $U\in\mathcal{B}$ and $x\in U$. Since $d\leq a_{U}$ and $a_{U}$ is attained, there exists $x_{0}\in F$ such that $I_{d}x_{0}\in I_{d}U=I_{d}O_{z}$. Moreover, $I_{d^{c}}y=I_{d^{c}}z\in I_{d^{c}}O_{z}$. By the $\sigma$-stability of $O_{z}$ and $F$, 
	$$I_{d}x_{0}+I_{d^{c}}y\in O_{z}\cap F.$$
	Hence $z\in cl_{\mathcal{T}|_{H}}(F)=F$, and therefore $I_{d}x=I_{d}z\in I_{d}F$. Since $0<d\leq a$, we obtain $x\in H\cap\llbracket F\rrbracket_{a}$. This proves the assertion.
\end{proof}

\begin{lemma}\label{lemm.fixed point}
	Let $E$ be an $L^{0}(\mathcal{F},\mathbb{K})$-module, $G$ a  nonempty $\sigma$-stable subset of $E$ and  $T:G\rightarrow 2^{G}\setminus \{\emptyset\}$ a $\sigma$-stable set-valued mapping. If for any $a\in B_{\mathcal{F}}$ with $a>0$, there exist $0<b\leq a$ and $x_{b}\in G$ such that $I_{b}x_{b}\in I_{b}T(x_{b})$, then $T$ has a fixed point.
\end{lemma}
\begin{proof}
	Let $$\mathcal{A}=\left\lbrace a>0:\ \text{there exists }x\in G\text{ such that }I_{a}x\in I_{a}T(x) \right\rbrace$$ and let $d=\bigvee\mathcal{A}$. It is easy to check that $d$ is attained, that is, there exists $x_{d}\in G$ such that $I_{d}x_{d}\in I_{d}T(x_{d})$. We claim that $d=1$, namely, $x_{d}$ is a fixed point of $T$.
	
	\par 
	Assume to the contrary that $d^{c}>0$. By the hypothesis, there exist $0<b\leq d^{c}$ and $x_{b}\in G$ such that $I_{b}x_{b}\in I_{b}T(x_{b})$. This implies that $b\in\mathcal{A}$, which contradicts the maximality of $d$.
\end{proof}

\begin{proof}[Proof of Theorem \ref{thm.Fan fixed point}]
	Let $a\in B_{\mathcal{F}}$ with $a>0$ be given. By Lemma \ref{lemm.fixed point}, it suffices to find $0<b\leq a$ and $x\in G$ such that $I_{b}x\in I_{b}\varphi(x)$.
	
	\par
	Suppose first that $G\setminus \llbracket\varphi^{-1}(y_{0})\rrbracket_{a}=\emptyset$ for some $y_{0}\in \varphi(G)$. Then  $I_{a}G\subseteq \llbracket \varphi^{-1}(y_{0}) \rrbracket_{a}$. By Lemma \ref{lemm.fundament1}, there exists $0<a'\leq a$ such that $$I_{a'}G\subseteq I_{a'}\varphi^{-1}(y_{0}).$$ 
	Since $\varphi^{-1}(y_{0})\subseteq G$, we obtain $I_{a'}G=I_{a'}\varphi^{-1}(y_{0})$. As $y_{0}\in G$, part (4) of Lemma \ref{lemm.stable1} gives
	$I_{a'}y_{0}\in I_{a'}\varphi(y_{0})$, as required.
	
	\par 
	We may therefore assume that $G\setminus\llbracket\varphi^{-1}(y)\rrbracket_{a}\neq\emptyset$ for any $y\in\varphi(G)$. Define a mapping $T:\varphi(G)\rightarrow 2^{G}\setminus \{\emptyset\}$ by 
	$$T(y)=G\setminus \llbracket \varphi^{-1}(y) \rrbracket_{a},~\forall y\in \varphi(G).$$ 
	It follows from part (2) of Lemma~\ref{lemm.stable1} and part (2) of Lemma~\ref{lemm.stable2} that $T$ is $\sigma$-stable. We claim that 
	$$\bigcap\{T(y):y\in \varphi(G)\}=\emptyset.$$ 
	In fact, if $x\in \bigcap\{T(y):y\in \varphi(G)\}$, then we have 
	$$x\notin \varphi^{-1}(y),~\forall y\in \varphi(G),$$ 
	which is a contradiction. By parts (1) and (3) of Lemma~\ref{lemm.stable2}, each $T(y)$ is a $\sigma$-stable closed subset of $G$, and hence stably compact. Since $T(\varphi(G))\subseteq G$, each $T(y)$ is also closed in $T(\varphi(G))$. If $T$ is a random KKM mapping, then Theorem~\ref{thm.KKM2} yields $\bigcap_{y\in\varphi(G)}T(y)\neq\emptyset$, a contradiction. Thus $T$ is not a random KKM mapping. Therefore, there exist $\{g_{1},\cdots,g_{n}\}\subseteq \varphi(G)$ and  $\{\lambda_{i}\mid i=1,\ldots,n\}\subseteq L^{0}_{+}(\mathcal{F})$ with $\sum_{i=1}^{n}\lambda_{i}=1$ such that $$g_{0}=\sum_{i=1}^{n}\lambda_{i}g_{i}\notin\bigcup\{T(g):g\in \sigma\{g_{1},\cdots,g_{n}\}\}.$$ 
	Then we have 
	$$	g_{0}\in\bigcap \{\llbracket\varphi^{-1}(g)\rrbracket_{a}:g\in \sigma\{g_{1},\cdots,g_{n}\}\},$$ 
	which, combined with part (4) of Lemma \ref{lemm.stable1}, implies that
	\begin{equation}\label{eq.KFan1}
		\sigma\{g_{1},\cdots,g_{n}\}\subseteq \llbracket\varphi(g_{0})\rrbracket_{a}.
	\end{equation}
	
	\par
	Let
	$$I=\{i\in\{1,\ldots,n\}\mid I_{a}g_{i}\notin I_{a}\varphi(g_{0})\}.$$
	If $I=\emptyset$, then $I_{a}g_{i}\in I_{a}\varphi(g_{0})$ for every $i=1,\ldots,n$. By the $L^{0}$-convexity of $\varphi(g_{0})$,
	$$I_{a}g_{0}=\sum_{i=1}^{n}\lambda_{i}I_{a}g_{i}\in I_{a}\varphi(g_{0}).$$
	Thus the required conclusion holds with $b=a$. We may therefore assume that $I\neq\emptyset$ and write
	$I=\{i_{1},\ldots,i_{k}\}$, where $1\leq k\leq n$.
	Then
	$$h(g_{i_{j}},\varphi(g_{0}))\wedge a>0,~\forall j=1,\ldots,k.$$
	We claim that
	$$\bigvee_{j=1}^{k}(h(g_{i_{j}},\varphi(g_{0}))\wedge a)<a.$$
	Otherwise,
	$$\bigvee_{j=1}^{k}(h(g_{i_{j}},\varphi(g_{0}))\wedge a)=a.$$
	Define $a_{1}=h(g_{i_{1}},\varphi(g_{0}))\wedge a$
	and
	$$a_{j}=(h(g_{i_{j}},\varphi(g_{0}))\wedge a)\wedge\left[\bigvee_{l=1}^{j-1}(h(g_{i_{l}},\varphi(g_{0}))\wedge a)\right]^{c},~\forall j=2,\ldots,k.$$
	Then $\{a_{1},\ldots,a_{k}\}\in p(a)$ and
	$$a_{j}\leq h(g_{i_{j}},\varphi(g_{0}))\wedge a,~\forall j=1,\ldots,k.$$
	Let
	$$g=\sum_{j=1}^{k} I_{a_{j}}g_{i_{j}}+I_{a^{c}}g_{i_{1}}.$$
	Then $g\in\sigma\{g_{1},\ldots,g_{n}\}$. For any $0<c\leq a$, there exists $j\in\{1,\ldots,k\}$ such that $c\wedge a_{j}>0$. Since $c\wedge a_{j}\leq h(g_{i_{j}},\varphi(g_{0}))$, we have
	$$I_{c\wedge a_{j}}g=I_{c\wedge a_{j}}g_{i_{j}}\notin I_{c\wedge a_{j}}\varphi(g_{0}).$$
	Consequently, $I_{c}g\notin I_{c}\varphi(g_{0})$ for any $0<c\leq a$, and hence $g\notin\llbracket\varphi(g_{0})\rrbracket_{a}$.
	This contradicts \eqref{eq.KFan1}.
	
	\par
	Let
	$$b=a\wedge\left[\bigvee_{j=1}^{k}(h(g_{i_{j}},\varphi(g_{0}))\wedge a)\right]^{c}.$$
	Then $0<b\leq a$ and $I_{b}g_{i}\in I_{b}\varphi(g_{0})$ for any $i=1,\ldots,n$.
	By the $L^{0}$-convexity of $\varphi(g_{0})$,
	$$I_{b}g_{0}=\sum_{i=1}^{n}\lambda_{i}I_{b}g_{i}\in\sum_{i=1}^{n}\lambda_{i}I_{b}\varphi(g_{0})=I_{b}\varphi(g_{0}).$$
\end{proof}

\begin{proof}[Proof of Theorem \ref{thm.coincidence}]
	Equip $E=E_{1}\times E_{2}$ with the product topology $\mathcal{T}=\mathcal{T}_{1}\times\mathcal{T}_{2}$. Then $(E,\mathcal{T})$ is a Hausdorff topological $L^{0}(\mathcal{F},\mathbb{K})$-module and $X\times Y$ is a nonempty stably compact $L^{0}$-convex subset of $E$ (see \cite[Theorem 2.40]{GMPTY2026}). Define a set-valued mapping $\varphi:X\times Y\rightarrow 2^{X\times Y}\setminus \{\emptyset\}$ by
	$$\varphi(x,y)=f^{-1}(y)\times g(x),~\forall (x,y)\in X\times Y.$$	
	By part (2) of Lemma \ref{lemm.stable1} and the $\sigma$-stability of $g$, $\varphi$ is $\sigma$-stable. 
	
	\par 
	Moreover, $\varphi(X\times Y)=X\times g(X)$ and, for any $(x,y)\in X\times g(X)$,  
	$$\varphi^{-1}(x,y)=g^{-1}(y)\times f(x)$$
	is a $\sigma$-stable $\mathcal{T}|_{X\times Y}$-open set. Further, $\varphi(x,y)=f^{-1}(y)\times g(x)$ is a nonempty $L^{0}$-convex set for any $(x,y)\in X\times Y$.
	
	\par 
	By Theorem~\ref{thm.Fan fixed point}, there exists $(x^{*},y^{*})\in X\times Y$ such that $(x^{*},y^{*})\in \varphi(x^{*},y^{*})$, which implies that $y^{*}\in f(x^{*})\cap g(x^{*})$.
\end{proof}

\subsection{Proofs of the results in Subsection~\ref{sec.2.4}}\label{sec.3.3}

  \par 
  To prove Theorems \ref{thm.Minimax inequality} and \ref{thm.Minimax principle}, we first give Lemma \ref{lemm.tools}.

  \begin{lemma}\label{lemm.tools}
  	Let $(E,\mathcal{T})$ be a Hausdorff topological $L^{0}(\mathcal{F},\mathbb{K})$-module and $G$ be a nonempty $\sigma$-stable subset of $E$. Assume that $(G,\mathcal{T}|_{G})$ is a stable topological space and $f:G\to\bar L^0(\mathcal F)$ is a $\sigma$-stable mapping
  	such that $f(z_0)\in L^0(\mathcal F)$ for some $z_0\in G$.
  	\begin{itemize}[widest={(1)},leftmargin=*]
  		\item [(1)] If $f$ is lower semicontinuous (resp., upper semicontinuous), then for any $a\in B_{\mathcal{F}}$ and any $\eta\in L^{0}(\mathcal{F})$,
  		$$\{x\in G\mid I_{a}f(x)\leq I_{a}\eta\}~( \text{resp.,}~\{x\in G\mid I_{a}f(x)\geq I_{a}\eta\}) $$
  		is $\sigma$-stable (if nonempty) and $\mathcal{T}|_{G}$-closed.
  		\item [(2)] If $f$ is lower semicontinuous (resp., upper semicontinuous), then for any $a\in B_{\mathcal{F}}$ and any $\eta\in L^{0}(\mathcal{F})$,
  		$$\{x\in G\mid f(x)>\eta~\text{on}~a\}~( \text{resp.,}~\{x\in G\mid f(x)<\eta~\text{on}~a\}) $$
  		is $\sigma$-stable (if nonempty) and $\mathcal{T}|_{G}$-open.
  		\item [(3)] If $G$ is stably compact and $f$ is lower semicontinuous (resp., upper semicontinuous), then there exists $x_{0}\in G$ such that 
  		$$f(x_{0})=\min_{x\in G}f(x)~\left( \text{resp.,}~ f(x_{0})=\max_{x\in G}f(x).\right) $$
  	\end{itemize}
  \end{lemma}
  \begin{proof}
  	It suffices to prove all three assertions when $f$ is lower semicontinuous, since the upper semicontinuous case follows by applying the corresponding results to $-f$.
  	
    \par 
    (1) Define a mapping $p_{a}:G\to G$ by
    $$p_a(x)=I_{a}x+I_{a^c}z_0,~\forall x\in G.$$
    Then $p_{a}$ is continuous, and the $\sigma$-stability of $f$ gives
    $$\{x\in G\mid I_{a}f(x)\leq I_{a}\eta\} =p_{a}^{-1}\bigl(\{x\in G:f(x)\leq I_a\eta+I_{a^c}f(z_0)\}\bigr).$$
    Since $I_{a}\eta+I_{a^c}f(z_0)\in L^0(\mathcal F)$, the lower semicontinuity of $f$ implies that $\{x\in G\mid I_{a}f(x)\leq I_{a}\eta\}$ is closed in $G$. Its $\sigma$-stability, when nonempty, follows from the $\sigma$-stability of $G$ and $f$.
    
  	\par 
  	(2) Fix $a\in B_{\mathcal{F}}$ and $\eta\in L^{0}(\mathcal{F})$, and let 
  	$$H=\{x\in G:f(x)>\eta~\text{on}~a\}.$$
  	If $H=\emptyset$, its openness is immediate. Assume that $H\neq\emptyset$. The $\sigma$-stability of $f$ implies that $H$ is $\sigma$-stable. Fix $x_{0}\in H$ and let $\mathcal{B}$ be a $\sigma$-bistable base for $\mathcal{T}|_{G}$. 
  	
  	\par 
  	For each $O\in \mathcal{B}$ with $x_{0}\in O$, let $b_{O}=a\wedge(\bigvee_{x\in O}(f(x)\leq \eta))$. Since both $f$ and $O$ are $\sigma$-stable, $b_{O}$ is attained, namely, there exists $y_{O}\in O$ such that $b_{O}=a\wedge (f(y_{O})\leq \eta)$ . Now define 
  	$$b:=\bigwedge\{b_{O}\mid O\in \mathcal{B}~\text{with}~x_{0}\in O\}=\bigwedge\{a\wedge (f(y_{O})\leq \eta)\mid O\in \mathcal{B}~\text{with}~x_{0}\in O\}.$$
  	We claim that $b=0$. Indeed, let 
  	$$C_{b}=\{x\in G:I_{b}f(x)\le I_{b}\eta\}.$$
  	By (1), $C_{b}$ is closed in $G$. For each $O\in\mathcal{B}$ with $x_{0}\in O$, the relation $b\leq b_{O}$ together with $b_{O}=a\wedge(f(y_{O})\leq \eta)$ gives $y_{O}\in O\cap C_{b}$, so $O\cap C_{b}\neq\emptyset$. Hence $x_{0}\in C_{b}$. Since $b\leq a$ and $x_{0}\in\{x\in G\mid f(x)>\eta~\text{on}~a\}$, the relations $I_{b}f(x_{0})\leq I_{b}\eta$ and $f(x_{0})>\eta$ on $b$ are incompatible unless $b=0$. Thus $b=0$.
  	
  	\par 
  	Consequently $b^{c}=1$, equivalently
  	$$\bigvee\{a^{c}\vee(f(y_{O})>\eta)\mid O\in\mathcal{B},\ x_{0}\in O\}=1.$$
  	Then there exists a sequence $\{O_{n},n\in \mathbb{N}\}$ in $\mathcal{B}$ with $x_{0}\in O_{n}$ for each $n\in \mathbb{N}$ such that 
  	$$\bigvee_{n=1}^{\infty} (a^{c}\vee (f(y_{O_{n}})> \eta))=1.$$
  	Let $a_{1}=a^{c}\vee (f(y_{O_{1}})> \eta)$ and 
  	$$a_{n}=(a^{c}\vee (f(y_{O_{n}})> \eta))\wedge [\bigvee_{k=1}^{n-1}(a^{c}\vee (f(y_{O_{k}})> \eta))]^{c},~\forall n\geq 2.$$
  	Then $\{a_{n},n\in \mathbb{N}\}\in p(1)$ and $a_{n}\leq a^{c}\vee (f(y_{O_{n}})> \eta)$ for each $n\in \mathbb{N}$. 
  	
    \par 
  	Let 
  	$$U=\sum_{n=1}^{\infty}I_{a_{n}}O_{n}.$$
  	Then $U\in \mathcal{B}$ and $x_{0}\in U$. For any $x=\sum_{n=1}^{\infty} I_{a_{n}}x_{n}\in U$ with $x_{n}\in O_{n}$ for each $n\in \mathbb{N}$, the definition of $b_{O_{n}}$ gives $a\wedge[f(x_{n})\leq \eta]\leq b_{O_{n}}$. Since $a_{n}\le b_{O_{n}}^{c}$, we obtain
  	$$(a\wedge a_{n})\wedge[f(x_{n})\leq \eta]=0,$$
  	so $f(x_{n})>\eta$ on $a\wedge a_{n}$. By the $\sigma$-stability of $f$,
  	$$f(x)=\sum_{n} I_{a_{n}}f(x_{n})>\eta~\text{on}~a,$$
  	and therefore $x_{0}\in U\subseteq H$. Consequently, $H$ is $\mathcal{T}|_{G}$-open.
  	
  	\par
  	(3) Let
  	$$G_0=\{x\in G:f(x)\leq f(z_0)\}.$$
  	Then $G_0$ is a nonempty $\sigma$-stable closed subset of $G$, and hence stably compact. For each $x\in G_{0}$, let 
  	$$F_x=\{y\in G:f(y)\leq f(x)\}.$$
  	Since $f(x)\leq f(z_0)\in L^0(\mathcal F)$,
  	$f(x)\vee(-m)\in L^0(\mathcal F)$ for each $m\in\mathbb N$.
  	Moreover,
  	\begin{align*}
  		F_x&=\{y\in G:f(y)\leq f(x)\}\\
  		&=\{y\in G:f(y)\leq \bigwedge_{m=1}^{\infty}(f(x)\vee(-m))\}\\ 		
  		&=\bigcap_{m=1}^{\infty}\{y\in G:f(y)\leq f(x)\vee(-m)\}.
  	\end{align*}
  	By the lower semicontinuity of $f$, $F_x$ is closed in $G$. It is nonempty and $\sigma$-stable by the $\sigma$-stability of $G$ and $f$. 
  	
  	\par For any $x_1,\ldots,x_n\in G_{0}$, let $b_i=\left(f(x_i)=\bigwedge_{j=1}^n f(x_j)\right)$ for each $i=1,\ldots,n$. Then $\bigvee_{i=1}^n b_i=1$. Define $a_1=b_1$ and 
  	$$ a_i=b_i\wedge\left(\bigvee_{j=1}^{i-1}b_j\right)^c,~\forall i=2,\ldots,n.$$
  	Then $\{a_{i},i=1,\ldots,n\}\in p(1)$ and $a_i\leq b_i$ for each $i$. Since $G_{0}$ is $\sigma$-stable, $\sum_{i=1}^{n}I_{a_{i}}x_{i}\in G_{0}$. Moreover, by the $\sigma$-stability of $f$,  
  	$$f\left( \sum_{i=1}^{n}I_{a_{i}}x_{i}\right) =\sum_{i=1}^{n}I_{a_{i}}f(x_i)=\sum_{i=1}^{n}I_{a_{i}}\left( \bigwedge_{j=1}^n f(x_j)\right)=\bigwedge_{j=1}^n f(x_j),$$
  	which implies that
  	$$\sum_{i=1}^{n}I_{a_{i}}x_{i}\in \left\lbrace y\in G:f(y)\leq \bigwedge_{j=1}^n f(x_j)\right\rbrace =\bigcap_{j=1}^{n}F_{x_{j}}.$$
  	Thus, $\{F_x\mid x\in G_{0}\}$ has the finite intersection property. By part (3) of Proposition \ref{prop.stable compactness}, there exists $x_0\in\bigcap_{x\in G_{0}}F_x$. Since $z_0\in G_0$, we have $f(x_0)\leq f(z_0)$, so $x_0\in G_0$. Consequently,
  	$$f(x_0)=\min_{x\in G_{0}}f(x).$$
  	
  	\par 
  	For any $x\in G$, let $a=(f(x)\leq f(z_0))$ and $z=I_ax+I_{a^c}z_0$. Since $G$ is $\sigma$-stable, $z\in G$. Moreover, 
  	by the $\sigma$-stability of $f$, 
  	$$f(z)=I_af(x)+I_{a^c}f(z_0)=f(x)\wedge f(z_0)\leq f(z_0).$$
  	Thus $z\in G_{0}$, and consequently
  	$$f(x_0)\leq f(z)=f(x)\wedge f(z_0)\leq f(x).$$
  	Therefore,
  	$$f(x_0)=\min_{x\in G}f(x).$$
  \end{proof}

 \begin{proof}[Proof of Theorem \ref{thm.Minimax inequality}]
 	Arbitrarily choose a point $g_{0}\in G$ and let $G_{1}=G-g_{0}$. Define a mapping $f_{1}:G_{1}\times G_{1}\rightarrow L^{0}(\mathcal{F})$ by 
 	$$f_{1}(x_{1},y_{1})=f(x_{1}+g_{0},y_{1}+g_{0})-f(g_{0},g_{0}),~\forall (x_{1},y_{1})\in G_{1}\times G_{1}.$$ 
 	Then, it is easy to check that
 	\begin{itemize}[widest={(1)}, leftmargin=*]
 		\item [(1)] $G_{1}$ is a nonempty stably compact $L^{0}$-convex subset of $E$ containing $\theta$; 
 		\item [(2)] $f_{1}$ is a $\sigma$-stable mapping such that $f_{1}(\theta,\theta)=0$;
 		\item [(3)] for each $x_{1}\in G_{1}$, $f_{1}(x_{1},\cdot)$ is  lower semicontinuous on $G_{1}$;
 		\item [(4)] for each $y_{1}\in G_{1}$, $f_{1}(\cdot,y_{1})$ is  $L^{0}$-quasiconcave;
 		\item [(5)] for each $y_{1}\in G_{1}$, 
 		$$\bigvee_{x_{1}\in G_{1}}f_{1}(x_{1},y_{1})=\bigvee_{x\in G}f(x,y_{1}+g_{0})-f(g_{0},g_{0})$$ 
 		and 
 		$$\bigvee_{x_{1}\in G_{1}}f_{1}(x_{1},x_{1})=\bigvee_{x\in G}f(x,x)-f(g_{0},g_{0}).$$ 
 	\end{itemize}	 
 	Thus, without loss of generality, we may assume that $\theta\in G$ and $f(\theta,\theta)=0$. Since $f$ is $\sigma$-stable, for any $a\in B_{\mathcal{F}}$ with $a>0$, we have
 	\begin{equation}\label{eq.Minimax1}
 		I_{a}\bigvee_{x\in G}f(x,y)=\bigvee_{x\in G}f(I_{a}x,I_{a}y)=\bigvee_{I_{a}x\in I_{a}G}f(I_{a}x,I_{a}y),~\forall y\in G	
 	\end{equation}
 	and 
 	\begin{equation}\label{eq.Minimax2}
 		I_{a}\bigvee_{x\in G}f(x,x)=\bigvee_{x\in G}f(I_{a}x,I_{a}x)=\bigvee_{I_{a}x\in I_{a}G}f(I_{a}x,I_{a}x).	
 	\end{equation}
 	
 	\par
 	First, we prove that the mapping $\bigvee_{x\in G}f(x,\cdot)$ attains its minimum on $G$. Let 
 	$$\mathcal{A}=\left\{c\in B_{\mathcal{F}}\mid ~\text{there~exists}~y\in G~\text{such~that}~I_{c}\bigvee_{x\in G}f(x,y)\in L^{0}(\mathcal{F})\right\}$$
 	and $b:=\bigvee \mathcal{A}$. If $b=0$, then $\bigvee_{x\in G}f(x,y)=+\infty$ for any $y\in G$, and the minimum of $\bigvee_{x\in G}f(x,\cdot)$ is attained. We then suppose that $b>0$. Since $G$, $L^{0}(\mathcal{F})$ and $f$ are $\sigma$-stable, it is easy to check that $b$ is attained, namely, there exists $y_{b}\in G$ such that 
 	\begin{equation}\label{eq.Minimax3}
 		I_{b}\bigvee_{x\in G}f(x,y_{b})=\bigvee_{I_{b}x\in I_{b}G}f(I_{b}x,I_{b}y_{b})\in L^{0}(\mathcal{F}).
 	\end{equation}
 	Moreover, by the definition of $b$,
 	$$I_{b^{c}}\bigvee_{x\in G}f(x,y)=\bigvee_{x\in G}f(I_{b^{c}}x,I_{b^{c}}y)=+\infty~\text{on }b^{c},~\forall y\in G.$$
 	Let $G_b=I_bG$. Since $\theta\in G$,
 	$$G_b=\{y\in G:I_{b^c}y=\theta\}.$$
 	Since $E$ is Hausdorff and the mapping $y\mapsto I_{b^c}y$ is continuous, $G_b$ is closed in $G$. Thus $G_b$ is a nonempty $\sigma$-stable closed subset of the stably compact set $G$, and hence stably compact. Moreover, $G_b$ is $L^0$-convex. Define a mapping $\psi:G_{b}\rightarrow \bar{L}^{0}(\mathcal{F})$ by
 	$$\psi(u)=\bigvee_{v\in G_{b}}f(v,u),~\forall u\in G_{b}.$$ 
 	Then $\psi$ is $\sigma$-stable. For any $\eta\in L^{0}(\mathcal{F})$,
 	\begin{align*}
 		\{u\in G_{b}:\psi(u)\leq\eta\}&=\bigcap_{v\in G_{b}}\{u\in G_{b}:f(v,u)\leq\eta\}\\
 		&=\bigcap_{v\in G_{b}}(G_{b}\cap\{y\in G:f(v,y)\leq\eta\})\\
 		&=\left( \bigcap_{v\in G_{b}}\{y\in G:f(v,y)\leq\eta\}\right) \cap G_{b}.
 	\end{align*}
 	For each $v\in G_{b}\subseteq G$, the lower semicontinuity of $f(v,\cdot)$ implies that $\{y\in G:f(v,y)\leq\eta\}$ is closed in $G$, and hence $\left( \bigcap_{v\in G_{b}}\{y\in G:f(v,y)\leq\eta\}\right) \cap G_{b}$ is closed in $G_{b}$. Therefore, $\psi$ is lower semicontinuous. Moreover,   $\psi(I_{b}y_{b})\in L^{0}(\mathcal{F})$ by \eqref{eq.Minimax3}. By part (3) of Lemma \ref{lemm.tools}, there exists $y_{*}\in G_{b}$ such that
 	$$\psi(y_{*})=\min_{u\in G_{b}}\psi(u).$$
 	Choose $y_{0}\in G$. Then $I_{b}y_{*}+I_{b^{c}}y_{0}\in G$. By \eqref{eq.Minimax1}, 
 	\begin{align*}
 		\bigvee_{x\in G}f(x,I_{b}y_{*}+I_{b^{c}}y_{0})&=I_{b}\bigvee_{x\in G}f(x,y_{*})+I_{b^{c}}\bigvee_{x\in G}f(x,y_{0})\\
 		&=I_{b}\bigvee_{I_{b}x\in G_{b}}f(I_{b}x,y_{*})+I_{b^{c}}\bigvee_{x\in G}f(x,y)\\
 		&\leq I_{b}\bigvee_{I_{b}x\in G_{b}}f(I_{b}x,I_{b}y)+I_{b^{c}}\bigvee_{x\in G}f(x,y)\\
 		&=I_{b}\bigvee_{x\in G}f(x,y)+I_{b^{c}}\bigvee_{x\in G}f(x,y)\\
 		&=\bigvee_{x\in G}f(x,y) 
 	\end{align*}
 	for any $y\in G$. Therefore, in both cases, $\bigvee_{x\in G}f(x,\cdot)$ attains its minimum on $G$.
 	
 	\par 
 	Second, we prove the desired inequality.
 	Let $\xi=\bigvee_{x\in G}f(x,x)$. If $\xi=+\infty$, then  the conclusion is valid. Suppose now that $\xi\neq+\infty$. Then $I_{(\xi<+\infty)}\xi \in L^{0}(\mathcal{F})$ and $I_{(\xi<+\infty)}G$ is a nonempty stably compact $L^{0}$-convex subset of $G$. By \eqref{eq.Minimax1} and \eqref{eq.Minimax2}, 
 	$$I_{(\xi<+\infty)}\bigvee_{x\in G}f(x,x)=\bigvee_{x\in I_{(\xi<+\infty)}G}f(x,x)$$
 	and 
 	$$I_{(\xi<+\infty)}\min_{y\in G}\bigvee_{x\in G}f(x,y)=\min_{y\in I_{(\xi<+\infty)}G}\bigvee_{x\in I_{(\xi<+\infty)}G}f(x,y).$$ 
    It therefore remains to show that
 	\begin{equation}\label{eq.Minimax4}
 		\min_{y\in I_{(\xi<+\infty)}G}\bigvee_{x\in I_{(\xi<+\infty)}G}f(x,y)\leq \bigvee_{x\in I_{(\xi<+\infty)}G}f(x,x).
 	\end{equation}
 	Define  $H:I_{(\xi<+\infty)}G\times I_{(\xi<+\infty)}G\rightarrow L^{0}(\mathcal{F})$ by
 	$$H(y,x)=(I_{(\xi<+\infty)}\xi-f(x,y))\wedge0,~\forall (x,y)\in I_{(\xi<+\infty)}G\times I_{(\xi<+\infty)}G$$
 	and $T:I_{(\xi<+\infty)}G\rightarrow 2^{I_{(\xi<+\infty)}G}\setminus \{\emptyset\}$ by 
 	$$T(x):=\{y\in I_{(\xi<+\infty)}G:f(x,y)\leq I_{(\xi<+\infty)} \xi\},~\forall x\in I_{(\xi<+\infty)}G.$$ 
 	Since $f$ is $\sigma$-stable, so is $H$. Moreover, the $L^{0}$-quasiconcavity of $f(\cdot,y)$ implies that $H(y,\cdot)$ is $L^{0}$-quasiconvex for each $y\in I_{(\xi<+\infty)}G$. Since $f(y,y)\leq I_{(\xi<+\infty)}\xi$, we have $H(y,y)=0$. Thus
 	$$H(y,y)\leq H(y,x) \Leftrightarrow f(x,y)\leq I_{(\xi<+\infty)}\xi.$$
 	By Lemma~\ref{lemm.example KKM}, $T$ is a random KKM mapping. Furthermore, for any $x\in I_{(\xi<+\infty)}G$, the lower semicontinuity of $f(x,\cdot)$ on $G$ together with the $\sigma$-stability of $f$ implies that $T(x)$ is a $\sigma$-stable closed subset of $G$, and is therefore stably compact. By Theorem \ref{thm.KKM2}, we have $$\bigcap_{x\in I_{(\xi<+\infty)}G}T(x)\neq \emptyset.$$ For any $y_{0}\in \bigcap_{x\in I_{(\xi<+\infty)}G}T(x)$, it follows that $$f(x,y_{0})\leq I_{(\xi<+\infty)} \xi,~\forall x\in I_{(\xi<+\infty)}G, $$ which implies that \eqref{eq.Minimax4} holds.
 \end{proof}

 \begin{proof}[Proof of Theorem~\ref{thm.Minimax principle}]
 	For any $(x,y)\in X\times Y$, by part (3) of Lemma~\ref{lemm.tools}, the mappings $f(x,\cdot)$ and $f(\cdot,y)$ attain their minimum and maximum, respectively. Define $m: X\rightarrow L^{0}(\mathcal{F})$ by
 	$$m(x)=\min_{y\in Y}f(x,y),~\forall x\in X,$$ 
 	and $M: Y\rightarrow L^{0}(\mathcal{F})$ by 
 	$$M(y)=\max_{x\in X}f(x,y),~\forall y\in Y.$$
 	Clearly, both $m$ and $M$ are $\sigma$-stable. For any $\xi\in L^{0}(\mathcal{F})$,  
 	$$\{x\in X\mid m(x)\ge\xi\}=\bigcap_{y\in Y}\{x\in X\mid f(x,y)\ge\xi\}$$
 	and 
 	$$\{y\in Y\mid M(y)\le\xi\}=\bigcap_{x\in X}\{y\in Y\mid f(x,y)\le\xi\}.$$
 	Since $f(x,\cdot)$ is lower semicontinuous for each $x\in X$ and $f(\cdot,y)$ is upper semicontinuous for each $y\in Y$, $m$ is upper semicontinuous on $X$ and $M$ is lower semicontinuous on $Y$. Again, by part (3) of Lemma~\ref{lemm.tools}, $m$ attains its maximum and $M$ attains its minimum. Let $\alpha:=\max_{x\in X}m(x)$ and $\beta:=\min_{y\in Y}M(y)$. Then $\alpha,\beta\in L^{0}(\mathcal{F})$ and $\alpha\leq\beta$ since $m(x)\leq f(x,y)\leq M(y)$ for any $(x,y)\in X\times Y$. Put $a:=(\alpha<\beta)$. It remains to show that $a=0$.
 	
 	\par 
 	Suppose that $a>0$, and let $r=\dfrac{\alpha+\beta}{2}\in L^{0}(\mathcal{F})$. Then 
 	\begin{equation}\label{eq.Minimax principle1}
 		m(x)\leq\alpha<r<\beta\leq M(y)~\text{on}~a, \forall (x,y)\in X\times Y.
 	\end{equation}
 	Define $S:X\rightarrow 2^{Y}$ by
 	$$S(x)=\{y\in Y\mid f(x,y)<r~\text{on}~a\},~\forall x\in X$$
 	and $T:Y\rightarrow 2^{X}$ by 
 	$$T(y)=\{x\in X\mid f(x,y)>r~\text{on}~a\},~\forall y\in Y.$$
 	By \eqref{eq.Minimax principle1}, a minimizing point for $f(x,\cdot)$ belongs to $S(x)$, and a maximizing point for $f(\cdot,y)$ belongs to $T(y)$. Thus $S(x)$ and $T(y)$ are nonempty. Clearly, $S$ and $T$ are both $\sigma$-stable. By the $L^{0}$-quasiconvexity of $f(x,\cdot)$ and the $L^{0}$-quasiconcavity of $f(\cdot,y)$, respectively, $S(x)$ and $T(y)$ are $L^{0}$-convex. Define a mapping $\Phi:X\times Y\rightarrow 2^{X\times Y}\setminus\{\emptyset\}$ by 
 	$$\Phi(x,y)=T(y)\times S(x),~\forall (x,y)\in X\times Y.$$
 	It is easy to check that $\Phi$ is $\sigma$-stable and $\Phi(x,y)$ is nonempty and $L^{0}$-convex for any $(x,y)\in X\times Y$. Moreover, for any $(u,v)\in\Phi(X\times Y)$,
 	$$\Phi^{-1}(u,v)=\{x\in X\mid f(x,v)<r~\text{on}~a\}\times\{y\in Y\mid f(u,y)>r~\text{on}~a\}.$$
 	By part (2) of Lemma~\ref{lemm.tools}, $\Phi^{-1}(u,v)$ is $\sigma$-stable and open in $(X\times Y,(\mathcal{T}_{1}\times\mathcal{T}_{2})|_{X\times Y})$. Since $X\times Y$ is a nonempty stably compact $L^{0}$-convex subset of the Hausdorff topological $L^{0}$-module $(E_{1}\times E_{2},\mathcal{T}_{1}\times\mathcal{T}_{2})$, Theorem~\ref{thm.Fan fixed point} yields $(x_{0},y_{0})\in\Phi(x_{0},y_{0})$. Consequently,
 	$$r<f(x_{0},y_{0})<r~\text{on}~a,$$
 	which contradicts $a>0$.
 	
 	\par 
 	Thus $a=0$, and hence $\alpha=\beta$.
 \end{proof}

\subsection{Proofs of the results in Subsection \ref{sec.2.5}}\label{sec.3.4}

\begin{proof}[Proof of Theorem \ref{thm.intersection}]
	Define $\varphi:X\to2^{X}\setminus\{\emptyset\}$ by
	$$\varphi(x)=\bigcap_{i=1}^{n}G_{i}(x),~\forall x\in X.$$
	By the $\sigma$-stability of $G_{1},\ldots,G_{n}$, it is easy to check that $\varphi$ is $\sigma$-stable. For any $x\in X$, since 
	$$\varphi(x)=\prod_{i=1}^{n}\{u_{i}\in X_{i}:(u_{i},x^{i})\in G_{i}\},$$
	$\varphi(x)$ is nonempty and $L^{0}$-convex. Moreover, for each $y\in\varphi(X)$
	$$\varphi^{-1}(y)=\bigcap_{i=1}^{n}G^{i}(y)$$
	is $\sigma$-stable and $\mathcal{T}|_{X}$-open. Since $X$ is a nonempty stably compact $L^{0}$-convex subset of $E$, by Theorem~\ref{thm.Fan fixed point}, there exists $x_{0}\in X$ such that $x_{0}\in\varphi(x_{0})$. By the definitions of $G_{i}(x_{0})$, $x_{0}\in G_{i}$ for every $i=1,\ldots,n$. Hence $\bigcap_{i=1}^{n}G_{i}\neq\emptyset$.
\end{proof}

\par 
To prove Theorem \ref{thm.Nash equilibrium}, we first establish a stable version of the tube lemma and then use it to prove the continuity of the maximum-value mapping.

\begin{lemma}\label{lemm.tubular neighborhood}
	Let $(K,\mathcal{T}_{1})$ and $(Y,\mathcal{T}_{2})$ be two stable topological spaces and $K$ be stably compact. Suppose that $y_{0}\in Y$ and $U\subseteq K\times Y$ is a $\sigma$-stable open set such that $K\times\{y_{0}\}\subseteq U$. Then there exists a $\sigma$-stable open set $V\subseteq Y$ containing $y_{0}$ such that $K\times V\subseteq U$.
\end{lemma}
\begin{proof}
	Let $\mathcal{O}$ be the family of nonempty $\sigma$-stable open subsets of $K$ and $\mathcal{W}$ be the family of nonempty $\sigma$-stable open subsets of $Y$. Since $(K,\mathcal{T}_{1})$ and $(Y,\mathcal{T}_{2})$ are two stable topological spaces, $\mathcal{O}$ and $\mathcal{W}$ are $\sigma$-stable. Define 
	$$\mathcal{A}:=\{O\in \mathcal{O}\mid \text{there~exists}~W\in \mathcal{W}~\text{with}~y_{0}\in W~\text{such~that}~O\times W\subseteq U\}$$
	It is easy to check that $\mathcal{A}$ is a $\sigma$-stable subfamily of $\mathcal{O}$ such that $K=\bigcup \mathcal{A}$. Since $K$ is stably compact, there exists $\{a_{n},n\in \mathbb{N}\}\in p(1)$ and a sequence $\{\mathcal{A}_{n}, n\in \mathbb{N}\}$ of nonempty finite subfamilies of $\mathcal{A}$ such that 
	\begin{equation}\label{eq.tubular neighborhood1}
		K=\bigcup \sum_{n=1}^{\infty}I_{a_{n}}\sigma(\mathcal{A}_{n}).
	\end{equation}
	For the sake of convenience, let $\mathcal{A}_{n}=\{O_{n,l}\mid l=1,\cdots k_{n}\}$ for each $n\in \mathbb{N}$. For each $O_{n,l}$, there exists $W_{n,l}\in \mathcal{W}$ with $y_{0}\in W_{n,l}$ such that  
	$$O_{n,l}\times W_{n,l} \subseteq U.$$
	Define $W_{n}=\bigcap_{l=1}^{k_{n}}W_{n,l}$ and 
	$$V:=\sum_{n=1}^{\infty}I_{a_{n}}W_{n}.$$ 
	Then $y_{0}\in W_{n}\in \mathcal{W}$ and $O_{n,l}\times W_{n} \subseteq U$ for each $n\in \mathbb{N}$ and each $l=1,\cdots, k_{n}$, and hence $y_{0}\in V\in \mathcal{W}$. 
	
	\par 
	For any $(x,y)\in K\times V$, by \eqref{eq.tubular neighborhood1}, there exists a sequence $\{\{b_{n,l},l=1,\cdots,k_{n}\},n\in \mathbb{N}\}$ in $p(1)$ such that 
	$$x\in \sum_{n=1}^{\infty}I_{a_{n}}\sum_{l=1}^{k_{n}}I_{b_{n,l}}O_{n,l}\subseteq K.$$
	Then 
	\begin{align*}
		(x,y)&\in \left(\sum_{n=1}^{\infty}I_{a_{n}}\sum_{l=1}^{k_{n}}I_{b_{n,l}}O_{n,l}\right) \times \left( \sum_{n=1}^{\infty}I_{a_{n}}\sum_{l=1}^{k_{n}}I_{b_{n,l}}W_{n}\right)\\
		&=\sum_{n=1}^{\infty}I_{a_{n}}\sum_{l=1}^{k_{n}}I_{b_{n,l}}(O_{n,l}\times W_{n})\\
		&\subseteq \sum_{n=1}^{\infty}I_{a_{n}}\sum_{l=1}^{k_{n}}I_{b_{n,l}} U\\
		&=U,
	\end{align*}
	which implies that $K\times V\subseteq U$.
\end{proof}

\begin{lemma}\label{lemm.maximum-value}
	Let $(E_{1},\mathcal{T}_{1})$ and $(E_{2},\mathcal{T}_{2})$ be two Hausdorff topological $L^{0}(\mathcal{F},\mathbb{K})$-modules with base $(\Omega,\mathcal{F},P)$, $K$ a nonempty stably compact subset of $E_{1}$ and $Y$ a nonempty subset of $E_{2}$ such that $(Y,\mathcal{T}_{2}|_{Y})$ is a stable topological space. Let $f:K\times Y\rightarrow L^{0}(\mathcal{F})$ be a $\sigma$-stable and continuous mapping from $(K\times Y,\mathcal{T}_{1}|_{K}\times\mathcal{T}_{2}|_{Y})$ to $(L^{0}(\mathcal{F}),\mathcal{T}_{c}(|\cdot|))$. Then, for each $y\in Y$, the function $f(\cdot,y)$ attains its maximum on $K$. Moreover, the maximum-value mapping $g:Y\rightarrow L^{0}(\mathcal{F})$ defined by
	$$g(y)=\max_{x\in K}f(x,y),$$
	is $\sigma$-stable and continuous from $(Y,\mathcal{T}_{2}|_{Y})$ to $(L^{0}(\mathcal{F}),\mathcal{T}_{c}(|\cdot|))$.
\end{lemma}
\begin{proof}
	For each fixed $y\in Y$, the mapping $x\mapsto f(x,y)$ is continuous and $\sigma$-stable. By part (3) of Lemma \ref{lemm.tools} and Remark \ref{rmk.semicontinuous}, there exists $x_{y}\in K$ such that
	$$g(y)=f(x_{y},y)\in L^{0}(\mathcal{F}).$$
	By the $\sigma$-stability of $K$ and $f$, it is easy to check that $g$ is $\sigma$-stable. 
	
	\par 
	Fix $y_{0}\in Y$ and $r\in L^{0}_{++}(\mathcal{F})$. Since the mapping $(x,y)\longmapsto f(x,y)-f(x,y_{0})$ is continuous and $\sigma$-stable, the set
	$$U:=\left\{(x,y)\in K\times Y\mid |f(x,y)-f(x,y_{0})|<\frac{r}{2}~\text{on}~1\right\}$$
	is $\sigma$-stable and open, and $K\times\{y_{0}\}\subseteq U$. By Lemma \ref{lemm.tubular neighborhood}, there exists a $\sigma$-stable open neighborhood $V$ of $y_{0}$ such that $K\times V\subseteq U$. Hence, 
	$$f(x,y_{0})-\frac{r}{2}<f(x,y)<f(x,y_{0})+\frac{r}{2}~\text{on}~1, ~\forall x\in K, y\in V,$$
	which implies that 
	$$g(y_{0})-\frac{r}{2}\leq g(y)\leq g(y_{0})+\frac{r}{2},~\forall y\in V.$$
	Therefore
	$$|g(y)-g(y_{0})|\leq\frac{r}{2}<r~\text{on}~1,~\forall y\in V.$$
	Thus, $g$ is continuous at $y_{0}$. Since $y_{0}$ was arbitrary, $g$ is continuous.
\end{proof}

\par 
We now prove Theorem~\ref{thm.Nash equilibrium}.

\begin{proof}[Proof of Theorem \ref{thm.Nash equilibrium}]
	Since each $X_{i}$ is stably compact, $X=\prod_{i=1}^{n}X_{i}$ is a stably compact subset of $(E,\mathcal{T})$ by \cite[Theorem 2.40]{GMPTY2026}.
	
	\par 
	Fix $i\in \{1,\ldots,n\}$. Define $v_{i}:X^{i}\rightarrow L^{0}(\mathcal{F})$ by 
	$$v_{i}(y):=\max_{u_{i}\in X_{i}}f_{i}(u_{i},y),~\forall y\in X^{i}.$$
	By Lemma~\ref{lemm.maximum-value} applied to $K=X_{i}$ and $Y=X^{i}$ with the mapping $(u_{i},y)\mapsto f_{i}(u_{i},y)$, $v_{i}$ is well-defined, $\sigma$-stable and continuous. Further, define $g_{i}:X\rightarrow L^{0}(\mathcal{F})$ by 
	$$g_{i}(x):=(v_{i}\circ p^{i})(x)=v_{i}(x^{i})=\max_{u_{i}\in X_{i}}f_{i}(u_{i},x^{i}),~\forall x\in X.$$ 
	Since the projection $p^{i}:X\to X^{i}$ is $\sigma$-stable and continuous, so is $g_{i}$. Moreover,
	\begin{equation}\label{eq.Nash1}
		f_{i}(x)\leq g_{i}(x)~~\text{and}~~g_{i}(u_{i},x^{i})=g_{i}(x), ~\forall x\in X,\ u_{i}\in X_{i}.	
	\end{equation}
	Fix $r\in L^{0}_{++}(\mathcal{F})$ and define  
	$$G_{i}^{r}=\{x\in X\mid f_{i}(x)-g_{i}(x)+r>0~\text{on}~1\}$$ 
	for each $i=1,\ldots,n$. Since $f_{i}-g_{i}$ is $\sigma$-stable and continuous, $G_{i}^{r}$ is $\sigma$-stable and $\mathcal{T}|_{X}$-open. The remainder of the proof is divided into two steps.
	
	\par 
	\noindent\textbf{Step 1.} We prove that $\bigcap_{i=1}^{n}G_{i}^{r}\neq\emptyset$ by verifying the hypotheses of Theorem~\ref{thm.intersection}.
	
	\par
	First, fix $i\in\{1,\ldots,n\}$ and $x\in X$. Since $g_{i}(u_{i},x^{i})=g_{i}(x)$ for any $u_{i}\in X_{i}$,	we have
	$$(G_{i}^{r})_{i}(x):=\{y\in X\mid (y_{i},x^{i})\in G_{i}^{r}\}=\{y\in X\mid f_{i}(y_{i},x^{i})>g_{i}(x)-r~\text{on}~1\}.$$
	This set is nonempty because the maximum defining $g_{i}(x)$ is attained. By the $L^{0}$-quasiconcavity of $f_{i}(\cdot,x^{i})$, the set
	$$\{u_{i}\in X_{i}\mid f_{i}(u_{i},x^{i})>g_{i}(x)-r~\text{on}~1\}$$
	is $L^{0}$-convex. Together with the $L^{0}$-convexity of the sets $X_{j}$, this implies that $(G_{i}^{r})_{i}(x)$ is $L^{0}$-convex.
	
	\par
	Second, fix $i\in\{1,\ldots,n\}$ and $y\in X$. The mapping $h_{i,y}:X\rightarrow X$ defined by 
	$$h_{i,y}(x):=(y_{i},x^{i}),~\forall x\in X$$
	is $\sigma$-stable and continuous. Therefore,
	$$(G_{i}^{r})^{i}(y):=\{x\in X\mid (y_{i},x^{i})\in G_{i}^{r}\}=h_{i,y}^{-1}(G_{i}^{r})$$
	is open in $(X,\mathcal{T}|_{X})$ and, if nonempty, $\sigma$-stable. 
	
	\par
	Thus, by Theorem~\ref{thm.intersection},
	$$\bigcap_{i=1}^{n}G_{i}^{r}\neq\emptyset.$$
	
	\par
	\noindent\textbf{Step 2.} We prove that there exists $x_{*}\in X$ such that 
	$$f_{i}(x_{*})=\bigvee\{f_{i}(y_{i},x_{*}^{i})\mid y_{i}\in X_{i}\},~\forall i\in \{1,2,\cdots,n\}.$$
	
	\par
	Define $H:X\to L^{0}(\mathcal{F})$ by
	$$H(x)=\bigvee_{i=1}^{n}\bigl(g_{i}(x)-f_{i}(x)\bigr),~\forall x\in X.$$
	By the continuity and $\sigma$-stability of $f_{i}$ and $g_{i}$, $H$ is continuous and $\sigma$-stable and, by \eqref{eq.Nash1}, satisfies $H(x)\geq 0$ for every $x\in X$. As $X$ is nonempty and $H$ is $L^{0}(\mathcal{F})$-valued, $H$ is proper; by Remark~\ref{rmk.semicontinuous} it is also lower semicontinuous. By part (3) of Lemma~\ref{lemm.tools}, there exists $x_{*}\in X$ such that $H(x_{*})=\bigwedge_{x\in X}H(x)$.
	
	\par
	For each $m\in\mathbb{N}$, Step 1 provides
	$$x_{m}\in\bigcap_{i=1}^{n}G_{i}^{\frac{1}{m}}.$$
	Consequently,
	$$0\leq H(x_{*})\leq H(x_{m})<\frac1m~\text{on}~1,~\forall m\in\mathbb{N}.$$
	It follows that $H(x_{*})=0$. Thus, for each $i=1,\ldots,n$,
	$$0\leq g_{i}(x_{*})-f_{i}(x_{*})\leq H(x_{*})=0.$$
	Therefore,
	$$f_{i}(x_{*})=g_{i}(x_{*})=\bigvee\{f_{i}(u_{i},x_{*}^{i})\mid u_{i}\in X_{i}\},~\forall i=1,\ldots,n,$$
	as required.
\end{proof}

\subsection{Proofs of the results in Subsection \ref{sec.2.6}}\label{sec.3.5}

\begin{proof}[Proof of Lemma \ref{lemm.two continuty}]
	Throughout this proof, for $i=1,2$, $x\in E_{i}$, an $L^{0}$-seminorm $p$ on $E_{i}$ and $r\in L^{0}_{++}(\mathcal{F})$, we write
	$$V_{i}(x,p,r):=\left\{y\in E_{i}\mid p(x-y)<r~\text{on}~1\right\};$$
	and, for any positive numbers $\varepsilon$ and $\lambda$ with $0<\lambda<1$, we write
	$$U_{i}(x,p,\varepsilon,\lambda):=\left\{y\in E_{i}\mid P\left\{\omega\in\Omega\mid p(x-y)(\omega)<\varepsilon\right\}>1-\lambda\right\}.$$
	
	(1) Given $x_{0}\in G_{1}$, assume that $f$ is $\mathcal{T}_{\varepsilon,\lambda}(\mathcal{P})$-continuous at $x_{0}$. By Remark \ref{rmk.continuous},  $f$ is also $\mathcal{T}_{\varepsilon,\lambda}(\mathcal{P}_{cc})$-continuous at $x_{0}$. We next show that $f$ is $\mathcal{T}_{c}(\mathcal{P}_{cc})$-continuous at $x_{0}$.
	
	\par 
	Let $V$ be a $\mathcal{T}_{c}((\mathcal{P}_{2})_{cc})$-neighborhood of $f(x_{0})$. Then there exist $q\in (\mathcal{P}_{2})_{cc}$ and $r\in L_{++}^{0}(\mathcal{F})$ such that $$V_{2}(f(x_{0}),q,r)\subseteq V.$$ It remains to show that there exists a $\mathcal{T}_{c}((\mathcal{P}_{1})_{cc})$-neighborhood $W$ of $x_{0}$ such that $$f(W\cap G_{1})\subseteq V_{2}(f(x_{0}),q,r)\cap G_{2}\subseteq V\cap G_{2}.$$ For each $n\in \mathbb{N}$, let
	$$\mathcal{A}_{n}:=\left\{ a\in B_{\mathcal{F}}
	\begin{array}{l}\left|
		\begin{array}{l}
			~\text{there~exists}~Q\in (\mathcal{P}_{1})_{cc}~\text{such~that}~ \\
			I_{a}f\left[V_{1}\left(x_{0},Q,\frac{1}{n}\right)\cap G_{1}\right]\subseteq I_{a}[V_{2}(f(x_{0}),q,r)\cap G_{2}]
		\end{array}\right.
	\end{array}
	\right\}. $$
	Then 
	$$\mathcal{A}_{n} \neq \emptyset,~\forall n\in \mathbb{N}.$$ 
	Furthermore, let 
	$$a_{n}=\bigvee \mathcal{A}_{n},~\forall n\in \mathbb{N}.$$ 
	We claim that each $a_{n}$ is attained.  In fact, for any fixed $n\in \mathbb{N}$, there exist a sequence  $\{b_{n,m},m\in \mathbb{N}\}$ in $\mathcal{A}_{n}$ and a sequence $\{Q_{n,m},m\in \mathbb{N}\}$ in $(\mathcal{P}_{1})_{cc}$ such that $a_{n}=\bigvee_{m=1}^{\infty}b_{n,m}$ and 
	$$I_{b_{n,m}}f\left[V_{1}\left(x_{0},Q_{n,m},\frac{1}{n}\right)\cap G_{1}\right]\subseteq I_{b_{n,m}}[V_{2}(f(x_{0}),q,r)\cap G_{2}], \forall m\in \mathbb{N}.$$ 
	Let $c_{n,1}=b_{n,1}$ and 
	$$c_{n,m}=b_{n,m}\wedge \left(\bigvee_{i=1}^{m-1}b_{n,i}\right)^{c},~\forall m\geq 2.$$ 
	Then $\{c_{n,m},m\in \mathbb{N}\}\in p(a_{n})$. Arbitrarily choose a $Q_{0}\in (\mathcal{P}_{1})_{cc}$ and let $$Q_{n}=\sum_{m=1}^{\infty}I_{c_{n,m}}Q_{n,m}+I_{a_{n}^{c}}Q_{0}, \forall n\in \mathbb{N},$$ 
	then $Q_{n}\in (\mathcal{P}_{1})_{cc}$ and
	\begin{align*}
		&I_{a_{n}}f\left[ V_{1}\left(x_{0},Q_{n},\frac{1}{n}\right)\cap G_{1}\right] \\
		&=I_{a_{n}} f \left[ \left( \sum_{m=1}^{\infty}I_{c_{n,m}}V_{1}\left(x_{0},Q_{n,m},\frac{1}{n}\right)+I_{a_{n}^{c}}V_{1}\left(x_{0},Q_{0},\frac{1}{n}\right)\right)  \cap G_{1}\right] \\
		&=\sum_{m=1}^{\infty}I_{c_{n,m}}f\left[ V_{1}\left(x_{0},Q_{n,m},\frac{1}{n}\right)\cap G_{1}\right]\\
		&\subseteq \sum_{m=1}^{\infty}I_{c_{n,m}}[ V_{2}(f(x_{0}),q,r)\cap G_{2}]\\
		&=I_{a_{n}}[V_{2}(f(x_{0}),q,r)\cap G_{2}],
	\end{align*}
	which implies that $a_{n}\in \mathcal{A}_{n}$. Therefore, each $a_{n}$ is attained, namely, there exists a sequence $\{Q_{n},n\in \mathbb{N}\}$ in $(\mathcal{P}_{1})_{cc}$ such that
	\begin{equation}\label{eq.two continuty}
		I_{a_{n}}f\left[ V_{1}\left(x_{0},Q_{n},\frac{1}{n}\right)\cap G_{1}\right] \subseteq I_{a_{n}}[ V_{2}(f(x_{0}),q,r)\cap G_{2}] ,~\forall n\in \mathbb{N}.	
	\end{equation}
	
	\par 
	If $\bigvee_{n=1}^{\infty}a_{n}<1$, then $b:=\bigwedge_{n=1}^{\infty}a_{n}^{c}>0$, namely, for any $Q\in (\mathcal{P}_{1})_{cc}$ and any $n\in \mathbb{N}$, $$I_{b'} f\left[V_{1}\left(x_{0},Q,\frac{1}{n}\right)\cap G_{1}\right]\nsubseteq I_{b'}[V_{2}(f(x_{0}),q,r)\cap G_{2}], \forall 0<b'\leq b.$$ By Lemma \ref{lemm.fundament1}, for any $Q\in (\mathcal{P}_{1})_{cc}$ and any $n\in \mathbb{N}$, there exists $x_{Q,n}\in V_{1}\left(x_{0},Q,\frac{1}{n}\right)\cap G_{1}$  such that
	\begin{equation}\label{eq.two continuty1}
		I_{b'} f(x_{Q,n}) \notin   I_{b'}[V_{2}(f(x_{0}),q,r)\cap G_{2}],~\forall 0<b'\leq b.
	\end{equation}
	Define a partial order $\preceq$ on $(\mathcal{P}_{1})_{cc}\times \mathbb{N}$ as follows: 
	$$(Q_{1}, n_{1})\preceq (Q_{2}, n_{2})~\text{iff}~Q_{1}(x)\leq Q_{2}(x)~\text{for~any}~x\in E_{1}, ~\text{and}~n_{1}\leq n_{2}.$$ 
	This order makes $(\mathcal{P}_{1})_{cc}\times\mathbb{N}$ directed. Then $\{x_{Q ,n},(Q, n)\in (\mathcal{P}_{1})_{cc}\times \mathbb{N}\}$ is a net in $G_{1}$. Furthermore, for any positive real numbers $\varepsilon,\lambda$ with $0<\lambda<1$ and any $Q_{0} \in (\mathcal{P}_{1})_{cc}$, there exists $(Q_{0}, n_{0})\in (\mathcal{P}_{1})_{cc}\times \mathbb{N}$ with $\frac{1}{n_{0}}< \varepsilon$ such that $$Q_{0}(x_{Q ,n}-x_{0})\leq Q(x_{Q ,n}-x_{0})<\frac{1}{n}\leq \frac{1}{n_{0}}<\varepsilon~\text{on}~1,~\forall (Q, n) \succeq (Q_{0}, n_{0}),$$ which implies that $\{x_{Q ,n},(Q, n)\in (\mathcal{P}_{1})_{cc}\times \mathbb{N}\}$ converges in $\mathcal{T}_{\varepsilon,\lambda}(\mathcal{P}_{cc})$ to $x_{0}$. However, by \eqref{eq.two continuty1} we have
	\begin{equation}\label{eq.two continuty2}
		q(f(x_{Q,n})-f(x_{0}))\geq r~\text{on}~b,~ \forall (Q, n)\in (\mathcal{P}_{1})_{cc}\times \mathbb{N}.
	\end{equation}
	Let $c_{1}=(r\geq 1)$ and $$c_{m}=\left(\frac{1}{m}\leq r<\frac{1}{m-1}\right),~\forall m\geq 2.$$ It is clear that $\{c_{m},m\in \mathbb{N}\}\in p(1)$. If there exists $c_{m}$ such that $b\wedge c_{m}=1$, then \eqref{eq.two continuty2} contradicts the assumption that $f$ is $\mathcal{T}_{\varepsilon,\lambda}(\mathcal{P}_{cc})$-continuous. Hence, there exists $c_{m}$ such that $0<b\wedge c_{m}<1$. It follows from \eqref{eq.two continuty2} that $U_{2}\left(f(x_{0}),q,\frac{1}{m},P(b\wedge c_{m})\right)$ is a $\mathcal{T}_{\varepsilon,\lambda}((\mathcal{P}_{2})_{cc})$-neighborhood  of $f(x_{0})$ such that $$f(x_{Q,n})\notin U_{2}\left(f(x_{0}),q,\frac{1}{m},P(b\wedge c_{m})\right),~\forall (Q, n)\in (\mathcal{P}_{1})_{cc}\times \mathbb{N}.$$ 
	This again contradicts the assumption that $f$ is $\mathcal{T}_{\varepsilon,\lambda}(\mathcal{P}_{cc})$-continuous.
	
	\par
	Thus, we conclude that $\bigvee_{n=1}^{\infty}a_{n}=1$. Let $b_{1}=a_{1}$ and 
	$$b_{n}=a_{n}\wedge \left(\bigvee_{i=1}^{n-1}a_{i}\right)^{c},~\forall n\geq 2.$$ 
	It is clear that $\{b_{n},n\in \mathbb{N}\}\in p(1)$. Furthermore, let  $u=\sum_{n=1}^{\infty}I_{b_{n}}n$ and $Q_{u} =\sum_{n=1}^{\infty}I_{b_{n}}Q_{n} $. Then $u\in L^{0}(\mathcal{F},\mathbb{N})$ and  $Q_{u} \in (\mathcal{P}_{1})_{cc}$. We claim that $V_{1}\left(x_{0},Q_{u}, \frac{1}{u}\right)$ is the desired $\mathcal{T}_{c}((\mathcal{P}_{1})_{cc})$-neighborhood of $x_{0}$. Indeed, for any given $y\in  V_{1}\left(x_{0},Q_{u},\frac{1}{u}\right)\cap G_{1}$, define 
	$$y_{n}=I_{b_{n}}y+I_{b_{n}^{c}}x_{0},~\forall n\in \mathbb{N}.$$ 
	Then $\{y_{n},n\in \mathbb{N}\}\subseteq V_{1}\left(x_{0},Q_{u},\frac{1}{u}\right)\cap G_{1}$. For each $n\in \mathbb{N}$, since 
	$$Q_{n} (y_{n}-x_{0}) =I_{b_{n}}Q_{n}(y-x_{0})=I_{b_{n}}Q_{u}(y-x_{0})<\frac{1}{n}~\text{on}~1,$$ 
	we have $y_{n}\in V_{1}\left(x_{0},Q_{n},\frac{1}{n}\right)\cap G_{1}$, which, combined with \eqref{eq.two continuty}, implies that $$I_{b_{n}}f(y_{n})\in I_{b_{n}}[V_{2}(f(x_{0}),q,r)\cap G_{2}].$$ 
	Consequently, $$f(y)=\sum_{n=1}^{\infty}I_{b_{n}}f(y_{n})\in \sum_{n=1}^{\infty}I_{b_{n}}[V_{2}(f(x_{0}),q,r)\cap G_{2}]=V_{2}(f(x_{0}),q,r)\cap G_{2}.$$
	
	\par 
	(2) Let $x_0\in G_1$ and assume that $f$ is $\mathcal{T}_c(\mathcal{P}_{cc})$-continuous at $x_0$. Given $q\in(\mathcal{P}_2)_{cc}$, $\varepsilon>0$ and $0<\lambda<1$, choose $p\in(\mathcal{P}_1)_{cc}$ and $r\in L^0_{++}(\mathcal{F})$ such that
	\begin{equation}\label{eq.lemm.continuity1}
		f[V_1(x_0,p,r)\cap G_1]\subseteq V_2(f(x_0),q,\varepsilon)\cap G_2.
	\end{equation}
	For $\lambda$, choose $\delta>0$ such that $P((r\leq\delta))<\frac{\lambda}{2}$. For any $y\in U_1(x_0,p,\delta,\frac{\lambda}{2})\cap G_1$, since $G_{1}$ and $f$ are $\sigma$-stable, $I_{(p(y-x_{0})<r)}y+I_{(p(y-x_{0})\geq r)}x_{0}\in G_{1}$ and 
	$$p(I_{(p(y-x_{0})<r)}y+I_{(p(y-x_{0})\geq r)}x_{0}-x_{0})=I_{(p(y-x_{0})<r)}p(y-x_{0})<r~\text{on}~1.$$
	Thus $I_{(p(y-x_{0})<r)}y+I_{(p(y-x_{0})\geq r)}x_{0}\in V_1(x_0,p,r)\cap G_1$. By \eqref{eq.lemm.continuity1}, we have
	$$q(f(I_{(p(y-x_{0})<r)}y+I_{(p(y-x_{0})\geq r)}x_{0})-f(x_{0}))=I_{(p(y-x_{0})<r)}q(f(y)-f(x_0))<\varepsilon~on~1,$$
	which implies that $q(f(y)-f(x_0))<\varepsilon$ on $(p(y-x_{0})<r)$. Hence, $(q(f(y)-f(x_0))<\varepsilon)\geq (p(y-x_{0})<r)$. It follows that 
		\begin{align*}
		P(q(f(y)-f(x_0))\geq\varepsilon)&\leq P(p(y-x_0)\geq r)\\
		& \leq P(p(y-x_0)\geq \delta)+P(\delta\geq r)\\
		& < \frac{\lambda}{2}+\frac{\lambda}{2}\\
		& =\lambda.
	\end{align*}
	Therefore,
	$$f[U_1(x_0,p,\delta,\lambda/2)\cap G_1]\subseteq U_2(f(x_0),q,\varepsilon,\lambda)\cap G_2.$$
	This proves that $f$ is $\mathcal{T}_{\varepsilon,\lambda}(\mathcal{P}_{cc})$-continuous at $x_0$, and hence $\mathcal{T}_{\varepsilon,\lambda}(\mathcal{P})$-continuous at $x_0$ by Remark~\ref{rmk.continuous}.
	
	\par 
	(3) Given $g\in G_{1}$, assume that $f$ is $\mathcal{T}_{c}(\mathcal{P})$-continuous at $g$. We next show that $f$ is $\mathcal{T}_{c}(\mathcal{P}_{cc})$-continuous at $g$.
	
	\par 
	Let $q\in (\mathcal{P}_{2})_{cc}$ and $\varepsilon\in L_{++}^{0}(\mathcal{F})$ be arbitrary. It suffices to show that there exist $Q\in (\mathcal{P}_{1})_{cc}$ and $r\in L_{++}^{0}(\mathcal{F})$ such that 
	$$f[V_{1}(g,Q,r)\cap G_{1}]\subseteq V_{2}(f(g),q,\varepsilon)\cap G_{2}.$$ 
	For $q$, there exist a sequence $\{\alpha_{n},n\in \mathbb{N}\}$ in $(\mathcal{P}_{2})_{f}$ and $\{a_{n},n\in \mathbb{N}\}\in p(1)$ such that $$q=\sum_{n=1}^{\infty}I_{a_{n}}\|\cdot\|_{\alpha_{n}}.$$ For each $\alpha_{n}\in (\mathcal{P}_{2})_{f}$, since $f$ is $\mathcal{T}_{c}(\mathcal{P})$-continuous at $g$, there exist $\beta_{n}\in (\mathcal{P}_{1})_{f}$ and $r_{n}\in L_{++}^{0}(\mathcal{F})$ such that 
	$$f[V_{1}(g,\|\cdot\|_{\beta_{n}},r_{n})\cap G_{1}]\subseteq V_{2}(f(g),\|\cdot\|_{\alpha_{n}},\varepsilon)\cap G_{2}.$$ 
	Let $Q=\sum_{n=1}^{\infty}I_{a_{n}}\|\cdot\|_{\beta_{n}}$ and $r=\sum_{n=1}^{\infty}I_{a_{n}}r_{n}$. Then $Q\in (\mathcal{P}_{1})_{cc}$ and $r\in L_{++}^{0}(\mathcal{F})$. It follows that
	\begin{align}
		f[V_{1}(g,Q,r)\cap G_{1}]&=f\left[V_{1}\left(g,\sum_{n=1}^{\infty}I_{a_{n}}\|\cdot\|_{\beta_{n}},\sum_{n=1}^{\infty}I_{a_{n}}r_{n}\right)\cap G_{1}\right]\nonumber\\
		&=f\left[\left(\sum_{n=1}^{\infty}I_{a_{n}}V_{1}(g,\|\cdot\|_{\beta_{n}},r_{n})\right)\cap G_{1}\right]\nonumber\\
		&=f\left[\sum_{n=1}^{\infty}I_{a_{n}}(V_{1}(g,\|\cdot\|_{\beta_{n}},r_{n})\cap G_{1})\right]\nonumber\\
		&=\sum_{n=1}^{\infty}I_{a_{n}}f[V_{1}(g,\|\cdot\|_{\beta_{n}},r_{n})\cap G_{1}]\nonumber\\
		&\subseteq \sum_{n=1}^{\infty}I_{a_{n}}[V_{2}(f(g),\|\cdot\|_{\alpha_{n}},\varepsilon)\cap G_{2}]\nonumber\\
		&=V_{2}(f(g),q,\varepsilon)\cap G_{2}.\nonumber
	\end{align}
\end{proof}

\par 
Theorem \ref{thm.tychonoff1} follows immediately from Theorem \ref{thm.inward1}, since, for every $x\in G$ with $x\neq f(x)$, the segment $[x,f(x)]$
contains the two distinct points $x$ and $f(x)$ of $G$. Here, we also give a direct proof based on the random KKM mapping principle (namely, Theorem \ref{thm.KKM2}).

\begin{proof}[Proof of Theorem \ref{thm.tychonoff1}]
	Fix $Q\in\mathcal{P}_{f}$. Let $F_{\|\cdot\|_{Q}}=\{y\in G:\|f(y)-y\|_{Q}=0\}$ and define $H_{\|\cdot\|_{Q}}:G\times G\rightarrow L^{0}(\mathcal{F})$ by 
	$$H_{\|\cdot\|_{Q}}(y,x)=\|f(y)-x\|_{Q},~\forall (y,x)\in G\times G.$$
	Since $f$ is $\sigma$-stable, it is easy to check that $H_{\|\cdot\|_{Q}}$ is also $\sigma$-stable. Moreover, for any $y,x_{1}, x_{2}\in G$ and any $\lambda\in L^{0}(\mathcal{F},[0,1])$, 
	\begin{align*}
		H_{\|\cdot\|_{Q}}(y,\lambda x_{1}+(1-\lambda)x_{2})&=\|f(y)-\lambda x_{1}-(1-\lambda)x_{2}\|_{Q}\\
		&\leq \lambda \|f(y)-x_{1}\|_{Q}+(1-\lambda) \|f(y)-x_{2}\|_{Q}\\
		&=\lambda H_{\|\cdot\|_{Q}}(y,x_{1}) +(1-\lambda) H_{\|\cdot\|_{Q}}(y,x_{2}).
	\end{align*}
	Hence $H_{\|\cdot\|_{Q}}(y,\cdot):G\rightarrow L^{0}(\mathcal{F})$ is $L^{0}$-convex for each $y\in G$. By Lemma \ref{lemm.example KKM}, the mapping $T_{\|\cdot\|_{Q}}:G\rightarrow 2^{G}\setminus \{\emptyset\}$ defined by 
	$$T_{\|\cdot\|_{Q}}(x)=\{y\in G:\|f(y)-y\|_{Q}\leq \|f(y)-x\|_{Q}\},~\forall x\in G $$
	is a random KKM mapping. For any $x\in G$, since $f$ is $\sigma$-stable and $\mathcal{T}_{c}(\mathcal{P}_{cc})$-continuous, $T_{\|\cdot\|_{Q}}(x)$ is a $\sigma$-stable closed subset of $G$, and hence stably compact. By Theorem \ref{thm.KKM2}, $\bigcap_{x\in G}T_{\|\cdot\|_{Q}}(x)$ is a nonempty $\sigma$-stable closed subset of $G$. We claim that $\bigcap_{x\in G}T_{\|\cdot\|_{Q}}(x)=F_{\|\cdot\|_{Q}}$. Indeed, for any $y\in\bigcap_{x\in G}T_{\|\cdot\|_{Q}}(x)$, 
	$$\|f(y)-y\|_{Q}\leq\|f(y)-x\|_{Q},\qquad \forall x\in G,$$
	and taking $x=f(y)\in G$ yields $\|f(y)-y\|_{Q}=0$. Hence $y\in F_{{\|\cdot\|_{Q}}}$, and so $\bigcap_{x\in G}T_{\|\cdot\|_{Q}}(x)\subseteq F_{\|\cdot\|_{Q}}$. Conversely,  the inclusion $F_{\|\cdot\|_{Q}}\subseteq\bigcap_{x\in G}T_{\|\cdot\|_{Q}}(x)$ follows immediately from the definition of $T_{\|\cdot\|_{Q}}$. 
	
	\par 
	For every $Q\in\mathcal{P}_{f}$, $\bigcap_{\|\cdot\|\in Q}F_{\|\cdot\|}=F_{\|\cdot\|_{Q}}\neq\emptyset$. Thus $\{F_{\|\cdot\|}\mid \|\cdot\|\in\mathcal{P}\}$ is a family of nonempty $\sigma$-stable closed subsets of $G$ which has the finite intersection property. By part (3) of Proposition \ref{prop.stable compactness}, $\bigcap_{\|\cdot\|\in\mathcal{P}}F_{\|\cdot\|}\neq\emptyset$.
	Consequently, every point in $\bigcap_{\|\cdot\|\in\mathcal{P}}F_{\|\cdot\|}$ is a fixed point of $f$. 
\end{proof}

\par 
In \cite{TMG2024}, Tu et al. established a strict separation theorem for stably compact and $\sigma$-stable closed sets in random locally convex modules (see Lemma~4.1 of \cite{TMG2024} for details). Lemma \ref{lemm.separation1} below is an equivalent form of \cite[Lemma~4.1]{TMG2024}.

\begin{lemma}\label{lemm.separation1}
	Let $(E,\mathcal{P})$ be an $RLC$ module over $\mathbb{K}$ with base $(\Omega,\mathcal{F},P)$, $K$ and $F$ be two nonempty subsets of $E$ such that $K$ is stably compact and $F$ is $\sigma$-stable and closed. Suppose that $a\in B_{\mathcal{F}}$ with $a>0$ satisfies $I_{b}K\cap
	I_{b}F=\emptyset$ for any $0<b\leq a$. Then there exists $V\in \mathcal{V}(\mathcal{P}_{cc})$ such that
	$$I_{b}(K+V)\cap I_{b}(F+V)=\emptyset$$
	for all $0<b\leq a$.
\end{lemma}
\begin{proof}
	If $a=1$, then the same argument as in \cite[Lemma~4.1]{TMG2024} applies, and the proof is complete. 
	
	\par 
	If $a<1$, let $A$ be an arbitrarily chosen representative of $a$, put $\Omega'=A$, $\mathcal{F}'=\{\Omega'\cap B\mid B\in \mathcal{F}\}$ and let $P':\mathcal{F}'\rightarrow[0,1]$ be defined by 
	$$P'(\Omega'\cap B)=\frac{P(\Omega'\cap B)}{P(\Omega')},~\forall B\in \mathcal{F}.$$
	Take $E'=I_{a}E$, $\mathcal{P}'=\{\|\cdot\|_{E'}:\|\cdot\|\in\mathcal{P}\}$, $K'=I_{a}K$, $F'=I_{a}F$. Then $(E',\mathcal{P}')$ is an $RLC$ module with base $(\Omega',\mathcal{F}',P')$ and the topology $\mathcal{T}_c(\mathcal{P}'_{cc})$ is just $\mathcal{T}_c(\mathcal{P}_{cc})|_{E'}$. Moreover, $K'$ is stably compact and $F'$ is $\sigma$-stable and closed. By \cite[Lemma~4.1]{TMG2024}, there exists $V'\in \mathcal{V}(\mathcal{P}'_{cc})$ such that
	$$I_{a'}(K'+V')\cap I_{a'}(F'+V')=\emptyset,~\forall 0<a'\leq a.$$
	For this $V'$, there exists $V\in \mathcal{V}(\mathcal{P}_{cc})$ such that $I_{a}V\subseteq V'$. Therefore, for every $0<b\leq a$, 
	$$I_{b}(K+V)\subseteq I_{b}(K'+V'),\qquad I_{b}(F+V)\subseteq I_{b}(F'+V'),$$
	which implies that 
	$$I_{b}(K+V)\cap I_{b}(F+V)=\emptyset.$$
\end{proof}

\begin{lemma}\label{lemm.separation2}
	Let $(E,\mathcal{P})$ be an $RLC$ module over $\mathbb{K}$ with base $(\Omega,\mathcal{F},P)$ and  $K$ be a nonempty stably compact subset of $E$. Suppose that $a\in B_{\mathcal{F}}$ with $a>0$ satisfies $\theta\notin I_{b}K$ for any $0<b\leq a$. Then there exist $p\in\mathcal{P}_{cc}$ and $r\in L^{0}_{++}(\mathcal{F})$ such that
	$$p(z)\geq r ~\text{on}~a,~\forall z\in K.$$
\end{lemma}
\begin{proof}
	Applying Lemma \ref{lemm.separation1} with $F=\{\theta\}$, there exists $V\in\mathcal{V}(\mathcal{P}_{cc})$ such that
	$$I_{b}(K+V)\cap I_{b}V=\emptyset,~\forall\,0<b\leq a.$$
	Choose $p\in\mathcal{P}_{cc}$ and $r\in L^{0}_{++}(\mathcal{F})$ such that $U:=\{x\in E\mid p(x)<r~\text{on}~1\}\subseteq V$. Since $\theta\in V$, it follows that
	$$ I_{b}K\cap I_{b}U=\emptyset,~\forall\,0<b\leq a.$$
	
	Fix $z\in K$ and let $c=a\wedge(p(z)<r)$. If $c>0$,  $p(I_{c}z)=I_{c}p(z)<r$ on $1$, so $I_{c}z\in U$ and hence $I_{c}z\in I_{c}U$. Since $z\in K$, we have $I_{c}z\in I_{c}K\cap I_{c}U$, contradicting $I_{c}K\cap I_{c}U=\emptyset$. Thus $c=0$, which implies that $p(z)\geq r$ on $a$. Since $z\in K$ was arbitrary, the proof is complete. 
\end{proof}

\begin{proof}[Proof of Theorem \ref{thm.inward1}]
	Let 
	$$\mathcal{A}:=\{a\in B_{\mathcal{F}} \mid I_{a}f(x)=I_{a}x~\text{for~some}~x\in G\}.$$
	It is clear that $\mathcal{A}\neq \emptyset$ since $0\in \mathcal{A}$. Let $a_{*}=\bigvee \mathcal{A}$. We claim that $a_{*}$ is attained. Indeed, there exists a sequence $\{a_{n},n\in \mathbb{N}\}$ in $\mathcal{A}$ and a sequence $\{x_{n},n\in \mathbb{N}\}$ in $G$ such that $a_{*}=\bigvee_{n=1}^{\infty}a_{n}$ and $I_{a_{n}}f(x_{n})=I_{a_{n}}x_{n}$ for each $n\in \mathbb{N}$. Set 
	$$ b_{1}=a_{1},\qquad b_{n}=a_{n}\wedge(\bigvee_{j=1}^{n-1}a_{j})^{c}, \quad \forall n\geq2,\qquad b_{0}=a_{*}^{c}.$$
	Then $\{b_{n},n\in \mathbb{N}\}\cup \{b_{0}\}\in p(1)$, and $b_{n}\leq a_{n} $ for each $n\in \mathbb{N}$. Fix $x_{0}\in G$. By the $\sigma$-stability of $G$ and $f$, we have $x_{*}:=\sum_{n=0}^{\infty}I_{b_{n}}x_{n}\in G$ and 
	\begin{equation}\label{eq.inward.local.fixed}
		I_{a_{*}}f(x_{*})=I_{a_{*}}\sum_{n=0}^{\infty}I_{b_{n}}f(x_{n})=I_{a_{*}}\sum_{n=0}^{\infty}I_{b_{n}}x_{n}=I_{a_{*}}x_{*}.
	\end{equation}
	Hence, $a_{*}\in \mathcal{A}$. 
	
	\par 
	If $a_{*}=1$, then $x_{*}$ is a fixed point of $f$. Suppose, therefore, that $a_{*}^{c}>0$. By the maximality of $a_{*}$, we have 
	\begin{equation}\label{eq.inward.no.local.fixed}
		I_{c}f(x)\neq I_{c}x, ~\forall x\in G,~\forall 0<c\leq a_{*}^{c}.
	\end{equation}
	
	\par 	 
	Next, consider the set $K=\{x-f(x):x\in G\}$. Since the mapping $x\mapsto x-f(x)$ is $\sigma$-stable and $\mathcal{T}_{c}(\mathcal{P}_{cc})$-continuous, $K$ is stably compact by \cite[Corollary 2.35]{GMPTY2026}. Moreover, \eqref{eq.inward.no.local.fixed} implies that $\theta\notin I_{c}K$ for any $0<c\leq a_{*}^{c}$. By Lemma \ref{lemm.separation2}, there exist $p\in\mathcal{P}_{cc}$ and $r\in L^{0}_{++}(\mathcal{F})$ such that
	\begin{equation}\label{eq.inward.positive.residual}
		p(x-f(x))\geq r~\text{on}~a_{*}^{c},~\forall x\in G.
	\end{equation}
	
	\par 	
	Define a mapping $T_{p}:G\to 2^{G}\setminus\{\emptyset\}$ by
	$$T_{p}(x)=\{y\in G\mid p(y-f(y))\leq p(x-f(y))\},~\forall x\in G.$$
	By Lemma \ref{lemm.example KKM}, $T_{p}$ is a random KKM mapping. For each $x\in G$, the set $T_{p}(x)$ is $\sigma$-stable and closed in $G$, and hence stably compact. Theorem \ref{thm.KKM2} yields $\bigcap_{x\in G}T_{p}(x)\neq\emptyset$. Choose $y\in \bigcap_{x\in G}T_{p}(x)$. Then
	\begin{equation}\label{eq.inward.best.approximation}
		p(y-f(y))\leq p(x-f(y)),\qquad \forall x\in G.
	\end{equation}
	Let $y_{0}=I_{a_{*}^{c}}y+I_{a_{*}}x_{*}$. By the $\sigma$-stability of $G$ and $f$, together with \eqref{eq.inward.local.fixed} and \eqref{eq.inward.best.approximation}, we have $y_{0}\in G$ and 
	\begin{align}\label{eq.inward.patched.residual}
		p(y_{0}-f(y_{0}))&=p(I_{a_{*}^{c}}(y-f(y))+I_{a_{*}}(x_{*}-f(x_{*})))\nonumber\\
		&= I_{a_{*}^{c}}p(y-f(y))\nonumber\\
		&\leq I_{a_{*}^{c}}p(x-f(y))+I_{a_{*}}p(x-f(x_{*}))\nonumber\\
		&=p(x-f(y_{0}))
	\end{align}
	for any $x\in G$. Moreover, by \eqref{eq.inward.local.fixed} and \eqref{eq.inward.positive.residual}, we have
	\begin{equation}\label{eq.inward1.10}
		p(y_{0}-f(y_{0}))>0\quad\text{on }a_{*}^{c},\qquad I_{a_{*}}y_{0}=I_{a_{*}}f(y_{0}).
	\end{equation}
	Hence, $y_{0}\neq f(y_{0})$. By the hypothesis, there exists $\lambda\in L^{0}(\mathcal{F},[0,1])$ such that $z=\lambda y_{0}+(1-\lambda)f(y_{0})\in G$ and $z\neq y_{0}$.  Taking $x=z$ in \eqref{eq.inward.patched.residual}, we obtain
	$$p(y_{0}-f(y_{0}))\leq p(z-f(y_{0}))=\lambda p(y_{0}-f(y_{0})).$$
	Thus $(1-\lambda)p(y_{0}-f(y_{0}))=0$, which, combined with \eqref{eq.inward1.10}, yields that $I_{a_{*}^{c}}(1-\lambda)=0$. It follows that $I_{a_{*}^{c}}z=I_{a_{*}^{c}}y_{0}$. Since $I_{a_{*}}z=I_{a_{*}}(\lambda y_{0}+(1-\lambda)f(y_{0}))=I_{a_{*}}y_{0}$ and
	$I_{a_{*}^{c}}z=I_{a_{*}^{c}}y_{0}$, we have $z=y_{0}$, a contradiction.
	
	\par 
	It follows that $a_{*}^{c}=0$. By \eqref{eq.inward.local.fixed}, $f(x_{*})=x_{*}$, completing the proof.
\end{proof}

\begin{proof}[Proof of Theorem \ref{thm.inward2}] 
	Suppose first that $f$ is inward. Fix $x\in G$ with $f(x)\neq x$. There exist $y\in G$ and
	$\alpha\in L^{0}_{+}(\mathcal{F})$ such that $f(x)=x+\alpha(y-x)$. Let 
	$$z=\frac{1}{1+\alpha}x+\frac{\alpha}{1+\alpha}y.$$
	By the $L^{0}$-convexity of $G$, $z\in G$.
	Moreover,
	$$z=\frac{\alpha}{1+\alpha}x+\frac{1}{1+\alpha}f(x)\in[x,f(x)],$$
	and
	$$z-x=\frac{1}{1+\alpha}(f(x)-x)\neq0.$$
	Thus $[x,f(x)]$ contains two distinct points of $G$. By Theorem \ref{thm.inward1}, $f$ has a fixed point.
	
	\par 
	If $f$ is a $\sigma$-stable $\mathcal{T}_{c}(\mathcal{P}_{cc})$-continuous outward mapping, define $g:G\rightarrow E$ by $$g(x)=2x-f(x),~\forall x\in G.$$ It is easy to check that $g$ is a $\sigma$-stable  $\mathcal{T}_{c}(\mathcal{P}_{cc})$-continuous inward mapping, and therefore, $g$ has a fixed point $x_{0}\in G$, namely, $$x_{0}=g(x_{0})=2x_{0}-f(x_{0}).$$ Consequently, we have $f(x_{0})=x_{0}$.
\end{proof}

\par 
To prove Theorems~\ref{thm.Schauder-Tychonoff} and~\ref{thm.Fan-Iokhvidov}, we first establish Lemmas~\ref{lemm.T1} and \ref{lemm.compact convex}.

\begin{lemma}\label{lemm.T1}
	Let $(E,\mathcal{P})$ be an $RLC$ module over $\mathbb{K}$ with base $(\Omega,\mathcal{F},P)$, $G$ a nonempty stably compact $L^{0}$-convex subset of $E$ and $V\in\mathcal{V}(\mathcal{P}_{cc})$. Let $f:G\to E$ be a $\sigma$-stable $\mathcal{T}_{c}(\mathcal{P}_{cc})$-continuous mapping such that $f(G)\subseteq G+V$. Then there exists $x_{0}\in G$ such that $f(x_{0})\in x_{0}+V$.
\end{lemma}
\begin{proof} 
	Define $T:G\rightarrow 2^{G}\setminus\{\emptyset\}$ by 
	$$T(x):=\{y\in G:f(x)-y\in V\}.$$ 
	By the $L^{0}$-convexity of $G$ and $V$, each $T(x)$ is $L^{0}$-convex. Since $V$ is $\mathcal{T}_{c}(\mathcal{P}_{cc})$-open and $f$ is continuous, each $T^{-1}(y)$ is $\mathcal{T}_{c}(\mathcal{P}_{cc})$-open in $G$. Moreover, the $\sigma$-stability of $G$ and $f$ implies that $T$ is $\sigma$-stable and $T^{-1}(y)$ is $\sigma$-stable for each $y\in T(G)$. By Theorem~\ref{thm.Fan fixed point}, there exists $x_{0}\in G$ such that $x_{0}\in T(x_{0})$, that is, $f(x_{0})\in x_{0}+V$.
\end{proof}

\begin{lemma}\label{lemm.compact convex}
	Let $(E,\mathcal{T})$ be a Hausdorff topological $L^{0}$-module and $G$ be a nonempty $L^{0}$-convex subset of $E$ such that $(G,\mathcal{T}|_{G})$ is a stable topological space.  The following statements hold:
	\begin{itemize}[widest={(1)}, leftmargin=*]
		\item [(1)] If $G_{1},\ldots,G_{n}$ are nonempty stably compact $L^{0}$-convex subsets of $G$, then $\operatorname{Conv}_{L^{0}}(G_{1}\cup\cdots\cup G_{n})$ is stably compact.
		\item [(2)] Let $\{a_{n}, n\in\mathbb{N}\}\in p(1)$, and let $\{K_{n},n\in\mathbb{N}\}$ be a sequence of nonempty stably compact subsets of $G$. Then $S:=\sum_{n=1}^{\infty}I_{a_{n}}K_{n}$ is stably compact. If each $K_{n}$ is also $L^{0}$-convex, then so is $S$.
	\end{itemize} 
\end{lemma}
\begin{proof}
	(1) Let $\Delta_{n-1}:=\operatorname{Conv}_{L^{0}}(\{e_{1},\ldots,e_{n}\})$ be the standard $L^{0}$-simplex in $L^{0}(\mathcal{F},\mathbb{R}^{n})$, and set $H=G_{1}\times\cdots\times G_{n}$. Define a mapping $f:\Delta_{n-1}\times H\to G$ by
	$$f(\xi,g)=\sum_{i=1}^{n}\xi_{i} g_{i}$$
	where $\xi=\sum_{i=1}^{n}\xi_{i} e_{i}\in\Delta_{n-1}$ with $\xi_{i}\in L^{0}_{+}(\mathcal{F})$ and $\sum_{i=1}^{n}\xi_{i}=1$, and $g=(g_{1},\ldots,g_{n})\in H$. Equip $\Delta_{n-1}$ with $\mathcal{T}_{c}(\|\cdot\|)$ and $H$ with the product topology. It is easy to check that $f$ is a $\sigma$-stable continuous mapping from $\Delta_{n-1}\times H$ to $G$. By the finite-product case of \cite[Theorem 2.40]{GMPTY2026}, $\Delta_{n-1}\times H$ is stably compact. Since $(G,\mathcal{T}|_{G})$ is a stable topological space, \cite[Corollary 2.35]{GMPTY2026} implies that $f(\Delta_{n-1}\times H)$ is stably compact. Thus it suffices to show that
	$$f(\Delta_{n-1}\times H)=\operatorname{Conv}_{L^{0}}(G_{1}\cup\cdots\cup G_{n}).$$
	
	\par
	First, it is clear that
	$$f(\Delta_{n-1}\times H)\subseteq \operatorname{Conv}_{L^{0}}(G_{1}\cup\cdots\cup G_{n}).$$
	
	\par
	On the other hand, take $(\xi,g),(\eta,h)\in\Delta_{n-1}\times H$ and $\alpha,\beta\in L^{0}_{+}(\mathcal{F})$ with $\alpha+\beta=1$.
	For each $i=1,\ldots,n$, let $s_{i}=\alpha\xi_{i}+\beta\eta_{i}$ and 
	$$v_{i}=\frac{\alpha\xi_{i}+\frac12 I_{(s_{i}=0)}}{s_{i}+I_{(s_{i}=0)}}g_{i}+\frac{\beta\eta_{i}+\frac12 I_{(s_{i}=0)}}{s_{i}+I_{(s_{i}=0)}}h_{i}.$$
	Since $\frac{\alpha\xi_{i}+\frac12 I_{(s_{i}=0)}}{s_{i}+I_{(s_{i}=0)}}, \frac{\beta\eta_{i}+\frac12 I_{(s_{i}=0)}}{s_{i}+I_{(s_{i}=0)}}\in L^{0}_{+}(\mathcal{F})$ and $\frac{\alpha\xi_{i}+\frac12 I_{(s_{i}=0)}}{s_{i}+I_{(s_{i}=0)}}+\frac{\beta\eta_{i}+\frac12 I_{(s_{i}=0)}}{s_{i}+I_{(s_{i}=0)}}=1$, $v_{i}\in G_{i}$ and $v=(v_{1},\ldots,v_{n})\in H$. Moreover,
	$$s_{i}v_{i}=\alpha\xi_{i} g_{i}+\beta\eta_{i} h_{i}.$$
	Let $\zeta=\alpha\xi+\beta\eta\in\Delta_{n-1}$. Then $\zeta\in\Delta_{n-1}$ and 
	$$\begin{aligned}
		\alpha f(\xi,g)+\beta f(\eta,h)
		&=\sum_{i=1}^{n}
		(\alpha\xi_{i} g_{i}+\beta\eta_{i} h_{i})\\
		&=\sum_{i=1}^{n} s_{i}v_{i}\\
		&=f(\zeta,v)
		\in f(\Delta_{n-1}\times H).
	\end{aligned}$$
	Thus $f(\Delta_{n-1}\times H)$ is $L^{0}$-convex. Since $G_{i}\subseteq f(\Delta_{n-1}\times H)$ for any $i=1,\ldots,n$, we obtain
	$$\operatorname{Conv}_{L^{0}}(G_{1}\cup\cdots\cup G_{n})\subseteq f(\Delta_{n-1}\times H).$$
	This proves the desired equality.
	
	\par 
	(2)  Let $\mathscr{F}$ be a stable filter on $S$ with a $\sigma$-bistable base $\mathscr{B}$.
	For each $n\in\mathbb{N}$, choose $z_{n}\in K_{n}$ and define a mapping $f_{n}:S\rightarrow K_{n}$ by
	$$f_{n}(x)=I_{a_{n}}x+I_{a_{n}^{c}}z_{n},~\forall x\in S.$$
	It is clear that each $f_{n}$ is $\sigma$-stable and continuous. Moreover, $\{f_{n}(A)\mid A\in\mathscr{B}\}$ generates a stable filter on $K_{n}$. By the stable compactness of $K_{n}$, this filter has a cluster point $y_{n}\in K_{n}$. Let $y=\sum_{n=1}^{\infty}I_{a_{n}}y_{n}\in S$. We next show that $y$ is a cluster point of $\mathscr{F}$. 
	
	\par 
	Let $A\in\mathscr{B}$ and let $U$ be a $\sigma$-stable $\mathcal{T}|_{S}$-open set such that $y\in U$. For each $n\in\mathbb{N}$, the mapping $g_{n}:K_{n}\to S$ defined by 
	$$g_{n}(z)=I_{a_{n}}z+I_{a_{n}^{c}}y,~\forall z\in K_{n}$$
	is continuous and satisfies $g_{n}(y_{n})=y$. Thus $g_{n}^{-1}(U)$ is a neighbourhood of $y_{n}$. By the choice of $y_{n}$, $g_{n}^{-1}(U)\cap f_{n}(A)\neq\emptyset$. Then, there exists $x_{n}\in A$ such that $f_{n}(x_{n})\in g_{n}^{-1}(U)$, which implies that $g_{n}(f_{n}(x_{n}))=I_{a_{n}}x_{n}+I_{a_{n}^{c}}y\in U$. Since $I_{a_{n}}g_{n}(f_{n}(x_{n}))=I_{a_{n}}x_{n}$ for each $n\in\mathbb{N}$, the $\sigma$-stability of $A$ and $U$ gives
	$$\sum_{n=1}^{\infty}I_{a_{n}}x_{n}=\sum_{n=1}^{\infty}I_{a_{n}}g_{n}(f_{n}(x_{n}))\in A\cap U.$$
	
	\par
	Hence $y$ is a cluster point of $\mathscr{F}$. By part (4) of Proposition~\ref{prop.stable compactness}, $S$ is stably compact. The last assertion follows directly from the $L^{0}$-convexity of the $K_{n}$.	
\end{proof}

\begin{proof}[Proof of Theorem \ref{thm.Schauder-Tychonoff}]
	First, we prove that, for each $V\in \mathcal{V}(\mathcal{P}_{cc})$, there exists $x_{V}\in G$ such that
	$$f(x_{V})-x_{V}\in V.$$
	For this, fix $V\in \mathcal{V}(\mathcal{P}_{cc})$. Since $f(G)$ is $\mathcal{T}_{c}(\mathcal{P}_{cc})$-dense in $cl_{\mathcal{T}_{c}(\mathcal{P}_{cc})}(f(G))$, the family $\{(x+V)\cap G| x\in f(G)\}$ is a $\sigma$-bistable $\mathcal{T}_{c}(\mathcal{P}_{cc})|_{G}$-open cover of $cl_{\mathcal{T}_{c}(\mathcal{P}_{cc})}(f(G))$. By the stable compactness of $cl_{\mathcal{T}_{c}(\mathcal{P}_{cc})}(f(G))$, there exist $\{a_{n},n\in \mathbb{N}\}\in p(1)$ and a sequence $\{G_{n},n\in \mathbb{N}\}$ of nonempty finite subsets of $f(G)$ such that  
	\begin{align*}
		cl_{\mathcal{T}_{c}(\mathcal{P}_{cc})}(f(G)) &\subseteq \bigcup \left\lbrace (x+V)\cap G\middle| \ x\in \sum_{n=1}^{\infty}I_{a_{n}}\sigma(G_{n})\right\rbrace \\
		&\subseteq \bigcup \left\lbrace x+V\middle| \ x\in \sum_{n=1}^{\infty}I_{a_{n}}\sigma(G_{n})\right\rbrace \\
		&=\sum_{n=1}^{\infty}I_{a_{n}}\sigma(G_{n})+V.
	\end{align*}
	Let $K=\sum_{n=1}^{\infty}I_{a_{n}}\operatorname{Conv}_{L^{0}}(G_{n})$. By Lemma \ref{lemm.compact convex}, $K$ is a nonempty stably compact and $L^{0}$-convex subset of $G$. Moreover,
	$$f(K)\subseteq f(G)\subseteq cl_{\mathcal{T}_{c}(\mathcal{P}_{cc})}(f(G))\subseteq \sum_{n=1}^{\infty}I_{a_{n}}\sigma(G_{n})+V\subseteq K+V.$$
	By Lemma \ref{lemm.T1}, there exists $x_{V}\in K\subseteq G$ such that $f(x_{V})\in x_{V}+V$.
	
	\par
	For each $V\in \mathcal{V}(\mathcal{P}_{cc})$, let 
	$$A_{V}:=\{x\in G\mid f(x)-x\in V\}$$ 
	and 
	$$F_{V}=cl_{\mathcal{T}_{c}(\mathcal{P}_{cc})}(f(A_{V})).$$
	By the first step, $A_{V}\neq\emptyset$. Furthermore, it is clear that  $A_{V}$ is $\sigma$-stable. By part (2) of \cite[Theorem 2.28]{GMPTY2026}, $F_{V}$ is a nonempty $\sigma$-stable closed subset of $cl_{\mathcal{T}_{c}(\mathcal{P}_{cc})}(f(G))$.
	
	\par 
	For any $V_{1},\ldots,V_{m}\in\mathcal{V}(\mathcal{P}_{cc})$, there exists $V_{0}\in\mathcal{V}(\mathcal{P}_{cc})$ such that $V_{0}\subseteq\bigcap_{i=1}^{m} V_{i}$. Then $\emptyset\neq F_{V_{0}}\subseteq\bigcap_{j=1}^{m}F_{V_{j}}$, which implies that $\{F_{V}:V\in \mathcal{V}(\mathcal{P}_{cc})\}$ has the finite intersection property. By part (3) of Proposition \ref{prop.stable compactness}, $\bigcap_{V\in\mathcal{V}(\mathcal{P}_{cc})}F_{V}\neq \emptyset$. 
	
	\par
	Choose $y_{0}\in\bigcap_{V\in\mathcal{V}(\mathcal{P}_{cc})}F_{V}$. We claim that $f(y_{0})=y_{0}$. If not, since $(E,\mathcal{T}_{c}(\mathcal{P}_{cc}))$ is Hausdorff, there exists $U\in \mathcal{V}(\mathcal{P}_{cc})$ such that
	$$(y_{0}+U)\cap(f(y_{0})+U)=\emptyset.$$
	By the continuity of $f$ and addition under $\mathcal{T}_{c}(\mathcal{P}_{cc})$, there exists $U_{1}\in\mathcal{V}(\mathcal{P}_{cc})$ such that $U_{1}+U_{1}\subseteq U$ and 
	\begin{equation}\label{eq.Schauder-Tychonoff1}
		f[(y_{0}+U_{1}+U_{1})\cap G]\subseteq f(y_{0})+U.	
	\end{equation}
	Since $y_{0}\in F_{U_{1}}$, there exists $x\in A_{U_{1}}$ such that $f(x)\in y_{0}+U_{1}$. Then
	$$x-y_{0}=(f(x)-y_{0})-(f(x)-x)\in U_{1}-U_{1}=U_{1}+U_{1}.$$
	It follows from \eqref{eq.Schauder-Tychonoff1} that $f(x)\in f(y_{0})+U$.
	
	\par 
	On the other hand, 
	$$f(x)\in y_{0}+U_{1}\subseteq y_{0}+U.$$ 
	Thus
	$$f(x)\in(y_{0}+U)\cap(f(y_{0})+U)=\emptyset,$$
	which is a contradiction. Therefore $f(y_{0})=y_{0}$, and the proof is complete.
\end{proof}

\begin{proof}[Proof of Theorem \ref{thm.Fan-Iokhvidov}]
	Fix $V\in \mathcal{V}(\mathcal{P}_{cc})$. Define a mapping $T_{V}:G\rightarrow 2^{G}\setminus \{\emptyset\}$ by
	$$T_{V}(x):=\{y\in G \mid f(x)-y\in K+V\},~\forall x\in G.$$
	Since $f(G)\subseteq G+K$ and both $G$ and $K+V$ are $L^{0}$-convex, $T_{V}(x)$ is a nonempty $L^{0}$-convex subset of $G$ for any $x\in G$. By the $\sigma$-stability of $G$, $K$, $V$ and $f$,  $T_{V}$ is $\sigma$-stable. Moreover, since $K+V$ is $\sigma$-stable and $\mathcal{T}_{c}(\mathcal{P}_{cc})$-open and $f$ is $\sigma$-stable and $\mathcal{T}_{c}(\mathcal{P}_{cc})$-continuous,  $T_{V}^{-1}(y)$ is $\sigma$-stable and $\mathcal{T}_{c}(\mathcal{P}_{cc})|_{G}$-open for each $y\in T_{V}(G)$. By Theorem~\ref{thm.Fan fixed point}, there exists $x\in G$ such that $x\in T_{V}(x)$.
	
	\par 
	Let $D=\{f(x)-x\mid x\in G\}$. Since $f-\mathrm{id}$ is $\sigma$-stable and $\mathcal{T}_{c}(\mathcal{P}_{cc})$-continuous, $D$ is stably compact by \cite[Corollary 2.35]{GMPTY2026}. By the first step,
	$$D\cap(K+V)\neq\emptyset,~\forall V\in\mathcal{V}(\mathcal{P}_{cc}).$$
	We claim that $D\cap K\neq\emptyset$. If not, by Lemma~\ref{lemm.hereditarily disjoint degree}, there exists $a\in B_{\mathcal{F}}$ with $a>0$ such that
	$$I_{c}D\cap I_{c}K=\emptyset,~\forall\,0<c\leq a.$$
	By Lemma~\ref{lemm.separation1}, there exists $V_{0}\in\mathcal{V}(\mathcal{P}_{cc})$ such that
	\begin{equation}\label{eq.Fan-Iokhvidov1}
		I_{a}(D+V_{0})\cap I_{a}(K+V_{0})=\emptyset.
	\end{equation}
	On the other hand, $D\cap(K+V_{0})\neq\emptyset$ and $\theta\in V_{0}$ imply that 
	$$(D+V_{0})\cap (K+V_{0})\neq \emptyset,$$
	which contradicts \eqref{eq.Fan-Iokhvidov1}. Therefore $D\cap K\neq\emptyset$, and hence there exists $x_{0}\in G$ such that $f(x_{0})-x_{0}\in K$.
\end{proof}

\subsection{Proofs of the results in Subsection \ref{sec.2.7}}\label{sec.3.6}

\begin{definition}[\cite{GZWW2021}]
	Let $(E,\mathcal{P})$ be an $RLC$-module over $\mathbb{K}$ with base $(\Omega,\mathcal{F},P)$ and $G$ a $\mathcal{T}_{\varepsilon,\lambda}(\mathcal{P})$-closed $L^{0}$-convex subset of $E$. $G$ is $L^{0}$-convexly compact if any family of $\mathcal{T}_{\varepsilon,\lambda}(\mathcal{P})$-closed $L^{0}$-convex subsets of $G$ has a nonempty intersection whenever this family has the finite intersection property.
\end{definition}

\begin{proof}[Proof of Theorem \ref{thm.Hartman-Stampacchia}]
	Since $E$ is $\mathcal{T}_{\varepsilon,\lambda}$-complete and $G$ is $\mathcal{T}_{\varepsilon,\lambda}$-closed and $L^{0}$-convex, $G$ is $\sigma$-stable. Define a mapping $T:G\rightarrow 2^{G}\setminus \{\emptyset\}$ by 
	$$T(x):=\{y\in G:\langle f(y),y-x\rangle \leq 0\},~\forall x\in G.$$ 
	It suffices to show that $\bigcap_{x\in G}T(x)\neq \emptyset $.
	
	\par 
	Let $H:G\times G\rightarrow L^{0}(\mathcal{F})$ be defined by 
	$$H(y,x)=\langle f(y),x\rangle, \forall (y,x)\in G\times G.$$
	Then, $H$ is $\sigma$-stable by the $\sigma$-stability of $f$. For any fixed $y\in G$, since $G$ is $L^{0}$-convex, $H(y,\cdot)$ is $L^{0}$-convex. Furthermore, it is clear that
	$$H(y,y)\leq H(y,x) \Leftrightarrow  \langle f(y),y\rangle\leq \langle f(y),x\rangle \Leftrightarrow \langle f(y),y-x\rangle\leq 0.$$
	Thus, by Lemma \ref{lemm.example KKM}, $T$ is a random KKM mapping. Define $\Gamma:G\rightarrow 2^{G}\setminus \{\emptyset\}$ by $$\Gamma(x)=\{y\in G:\langle f(x),y-x\rangle \leq 0\},~\forall x\in G.$$ 
	Clearly, $\Gamma$ is $\sigma$-stable.  
	
	\par 
	Fix $x\in G$ and let $y\in T(x)$. Then $$\langle f(y),y-x\rangle \leq 0.$$ Since $f$ is monotone, we have 
	$$\langle f(y),y-x\rangle-\langle f(x),y-x\rangle=\langle f(y)-f(x),y-x\rangle\geq 0,$$ 
	namely, 
	$$\langle f(x),y-x\rangle\leq \langle f(y),y-x\rangle\leq 0,$$ 
	implying $y\in \Gamma(x)$. Hence, 
	$$T(x)\subseteq \Gamma(x),~\forall x\in G.$$
	
	\par 
	For any $y_{0}\in \bigcap_{x\in G}\Gamma(x)$, fix $x\in G$. If $x=y_{0}$, the desired inequality is trivial.
	Suppose that $x\neq y_{0}$ and let 
	$$ z_{\lambda}:=\lambda x+(1-\lambda)y_{0},~\forall \lambda\in L^{0}(\mathcal{F},[0,1]).$$ 
	Since $G$ is $L^{0}$-convex, $z_{\lambda}\in G$ and $y_{0}\in \Gamma(z_{\lambda})$ for any $\lambda\in L^{0}(\mathcal{F},[0,1])$, which implies that 
	$$\lambda \langle f(z_{\lambda}),y_{0}-x\rangle\leq 0,~\forall \lambda\in L^{0}(\mathcal{F},[0,1]).$$ 
	Let $\{\lambda_{\alpha},\alpha\in \Lambda\}$ be a net in $L^{0}(\mathcal{F},(0,1])$ that converges in $\mathcal{T}_{c}$ to $0$. Since $z_{\lambda_{\alpha}}\in l(x,y_{0})\cap G$ and 
	$$\|z_{\lambda_{\alpha}}-y_{0}\|=\lambda_{\alpha}\|x-y_{0}\|\xrightarrow{\mathcal{T}_{c}}0,$$
	the $\mathcal{T}_{c}$-hemicontinuity of $f$ yields
	$$f(z_{\lambda_{\alpha}})\xrightarrow{\mathcal{T}_{c}}f(y_{0}).$$
	By the random Cauchy--Schwarz inequality \cite[Proposition 4.1]{Guo1999}, 
	$$\langle f(z_{\lambda_{\alpha}}),y_{0}-x\rangle\xrightarrow{\mathcal{T}_{c}}\langle f(y_{0}),y_{0}-x\rangle.$$	
	Since $\{\xi\in L^{0}(\mathcal{F})\mid \xi\leq 0\}$ is $\mathcal{T}_{c}(|\cdot|)$-closed and $\langle f(z_{\lambda_{\alpha}}),y_{0}-x\rangle\leq 0$ for each $\alpha$, we obtain
	$$\langle f(y_{0}),y_{0}-x\rangle\leq 0.$$
	Thus, $y_{0}\in T(x)$ for any $x\in G$, and $\bigcap_{x\in G}\Gamma(x)\subseteq \bigcap_{x\in G}T(x) $. On the other hand, it is clear that $\bigcap_{x\in G}T(x) \subseteq\bigcap_{x\in G}\Gamma(x)$. Consequently, we have 
	$$\bigcap_{x\in G}T(x) =\bigcap_{x\in G}\Gamma(x).$$
	
	\par 
	For each $x\in G$, $\Gamma(x)$ is nonempty, $\sigma$-stable, closed and $L^{0}$-convex, since $x\in\Gamma(x)$ and the mapping $y\mapsto\langle f(x),y-x\rangle$ is continuous and preserves $L^0$-convex combinations. Moreover, $\Gamma$ is a random KKM mapping because it is $\sigma$-stable and $T(x)\subseteq\Gamma(x)$ for every $x\in G$. By Theorem \ref{thm.KKM2}, the family $\{\Gamma(x):x\in G\}$ has the finite intersection property.
	
	\par 	
	The $\mathcal{T}_{\varepsilon,\lambda}$-complete $RIP$ module $E$ is random reflexive (for the definition of random reflexivity, see \cite{GL2005})  by the random Riesz representation theorem \cite[Lemma 4.1]{Guo1999}. Since $G$ is closed, a.s. bounded and $L^{0}$-convex, \cite[Corollary 2.23]{GZWW2021} implies that $G$ is $L^{0}$-convexly compact. Consequently,
	$$\bigcap_{x\in G}T(x)=\bigcap_{x\in G}\Gamma(x)\neq\emptyset.$$
\end{proof}

\par 
In 2020, Guo et al. established the random Kirk fixed point theorem for nonexpansive self-mappings of $L^0$-convexly compact sets with random normal structure in $\mathcal{T}_{\varepsilon,\lambda}$-complete $RN$ modules \cite[Theorem~3.8]{GZWG2020}. Corollary~\ref{coro.Browder-Goedhe-Kirk} below is a special case of this theorem. Indeed, a $\mathcal{T}_{\varepsilon,\lambda}$-complete $RIP$ module $E$ is random reflexive, so every nonempty closed a.s. bounded $L^0$-convex subset $G$ of $E$ is $L^0$-convexly compact by \cite[Corollary~2.23]{GZWW2021}. Moreover, the random uniform convexity of $E$ ensures that $G$ has random normal structure by \cite[Theorem~3.6]{GZWG2020}. Here we give an alternative proof of this corollary by applying Theorem~\ref{thm.Hartman-Stampacchia}.

\begin{corollary}\label{coro.Browder-Goedhe-Kirk}
	Let $(E,\langle\cdot,\cdot\rangle)$ be a $\mathcal{T}_{\varepsilon,\lambda}$-complete $RIP$ module over $\mathbb{R}$ with base $(\Omega,\mathcal{F},P)$, $G$ a nonempty closed a.s. bounded $L^{0}$-convex subset of $E$ and $f:G\rightarrow G$ a nonexpansive mapping, i.e., $$\|f(x)-f(y)\|\leq \|x-y\|,~\forall x,y\in G.$$ Then $f$ has a fixed point.
\end{corollary}
\begin{proof}
	Let $g(x)=x-f(x)$ for any $x\in G$. It is easy to check that $g$ is $\sigma$-stable (by \cite[Lemma 2.11]{GZWG2020}), $\mathcal{T}_{c}(\|\cdot\|)$-continuous and monotone. By Theorem \ref{thm.Hartman-Stampacchia}, there exists $y_{0}\in G$ such that 
	$$\langle g(y_{0}),y_{0}-x\rangle\leq 0,~\forall x\in G.$$ 
	Taking $x=f(y_{0})$ yields $\|y_{0}-f(y_{0})\|^{2}\leq0$, and hence $y_{0}=f(y_{0})$.
\end{proof}

\begin{corollary}\label{coro.Nikodym}
	Let $(E,\langle\cdot,\cdot\rangle)$ be a $\mathcal{T}_{\varepsilon,\lambda}$-complete $RIP$ module over $\mathbb{R}$ with base $(\Omega,\mathcal{F},P)$ and $G$ be a nonempty closed a.s. bounded $L^{0}$-convex subset of $E$. Then for each $x_{0}\in E$ there exists a unique $y_{0}\in G$ such that $$\|x_{0}-y_{0}\|=\bigwedge_{x\in G} \|x_{0}-x\|.$$
\end{corollary}
\begin{proof}
	Since every $RIP$ module is random strictly convex (for the definition of random strict convexity, see \cite[Definition 3.2]{GZ2010}), uniqueness follows from \cite[Proposition 3.2]{GZ2010}. Define a mapping $f:G\rightarrow E$ by $$f(x)=x-x_{0},~\forall x\in G.$$ It is clear that  $f$ is $\sigma$-stable, continuous and monotone. By Theorem \ref{thm.Hartman-Stampacchia}, there exists $y_{0}\in G$ such that 
	$$\langle y_{0}-x_{0},y_{0}-x\rangle\leq 0,~\forall x\in G,$$ 
	which, combined with the parallelogram law \cite[Proposition~4.1]{Guo1999}, implies that
	$$\|x_{0}-y_{0}\|=\bigwedge_{x\in G} \|x_{0}-x\|.$$
\end{proof}

\section{Concluding remark}\label{sec.4}

The proofs in the preceding sections start from the random Brouwer fixed point theorem \cite[Theorem~2.1]{TMGC2026}. The solid arrows in the diagram below summarize the implications established in those proofs. The dashed arrows indicate the implications proved in this section. The random Tychonoff fixed point theorem represents the five fixed point theorems in Subsection~\ref{sec.2.6}, whose individual relations are given in (i).

\begin{center}
	\resizebox{\linewidth}{!}{%
		\begin{tikzpicture}[
			y=7.5mm,
			every node/.style={font=\normalsize,inner sep=5pt},
			implication/.style={
				-{Stealth[length=2.2mm]},
				line width=.65pt
			},
			converse/.style={implication,dashed},
			arrowlabel/.style={font=\small,fill=white,inner sep=2pt}
			]
			
			\node (T) at (0,5.4) {Random Tychonoff fixed point theorem};
			\node (B) at (0,3.4) {Random Brouwer fixed point theorem};
			\node (K0) at (3.5,1.7) {Random KKM lemma};
			\node (K) at (3.5,-1.7) {Random KKM mapping principle};
			\node (F) at (0,-3.4) {Random Fan--Browder fixed point theorem};
			\node (I) at (-3.5,-1.7) {Random intersection theorem};
			\node (N) at (-3.5,1.7) {Conditional Nash equilibrium theorem};
			\node (C) at (0,0) {Random coincidence theorem};
			\node (MP) at (0,-4.9) {Random minimax principle};
			\node (MI) at (0,-6.5) {Random minimax inequality};
			\node (H) at (0,-8.1) {Random Hartman--Stampacchia theorem};
			
			\draw[implication]
			(B.east) to[out=0,in=90] (K0.north);
			\draw[implication]
			(K0.south) -- (K.north);
			\draw[implication]
			(K.south) -- (K.south |- F.north);
			\draw[implication]
			(I.south |- F.north) -- (I.south);
			\draw[implication]
			(I.north) -- (N.south);
			\draw[implication]
			(F.north) -- (C.south);
			\draw[implication]
			(F.south) -- (MP.north);
			
			\coordinate (R) at ([xshift=10mm]K.east);
			\draw[line width=.65pt]
			(K.east) -- (R);
			\draw[implication,rounded corners=6pt]
			(R) -- (R |- T.east) -- (T.east);
			\draw[implication]
			(R |- MI.east) -- (MI.east);
			\draw[implication,rounded corners=6pt]
			(R) -- (R |- H.east) -- (H.east);
			\fill (R) circle[radius=.8pt];
			\fill (R |- MI.east) circle[radius=.8pt];
			
			\draw[converse]
			(T.south) -- node[arrowlabel,right] {(i)} (B.north);
			
			\coordinate (L) at ([xshift=-10mm]N.west);
			\draw[converse,rounded corners=6pt]
			(MI.west)
			-- (L |- MI.west)
			-- node[arrowlabel,left,pos=.55] {(ii)} (L |- B.west)
			-- (B.west);
			
			\draw[converse]
			(C.north) -- node[arrowlabel,right] {(iii)} (B.south);
			\draw[converse]
			(N.north) to[out=90,in=180]
			node[arrowlabel,left,pos=.5] {(iv)} (B.west);
			
		\end{tikzpicture}%
	}
\end{center}

\par 
Let $G$ be a nonempty $\sigma$-stable a.s. bounded $\mathcal{T}_{c}(\|\cdot\|)$-closed $L^0$-convex subset of $L^0(\mathcal{F},\mathbb{R}^n)$ and $T:G\to G$ be a
$\sigma$-stable $\mathcal{T}_{c}(\|\cdot\|)$-continuous mapping. By \cite[Corollary~2.51]{GMPTY2026}, $G$ is stably compact. We next show that each of the results corresponding to the dashed arrows yields a fixed point of $T$.

\par 
\noindent\textit{\textbf{Implication (i).}}  The proofs in Subsection~\ref{sec.3.5} give the following implications:
$$
\begin{aligned}
	\text{Random KKM mapping principle}
	&\implies \text{Theorem~\ref{thm.inward1}}
	\implies
	\begin{cases}
		\text{Theorem~\ref{thm.tychonoff1}},\\
		\text{Theorem~\ref{thm.inward2}},
	\end{cases}
	\\[1ex]
	\text{Random Fan--Browder fixed point theorem}
	&\implies
	\begin{cases}
		\text{Lemma~\ref{lemm.T1}}
		\implies \text{Theorem~\ref{thm.Schauder-Tychonoff}},\\[1ex]
		\text{Theorem~\ref{thm.Fan-Iokhvidov}}.
	\end{cases}
\end{aligned}
$$
Each of these five theorems gives a fixed point of $T$. Indeed, Theorem~\ref{thm.tychonoff1} applies directly in the random Euclidean space. For Theorem~\ref{thm.inward2}, the inward condition holds since $T(G)\subseteq G$. Further,  $cl_{\mathcal{T}_{c}(\|\cdot\|)}(T(G))$ is a $\sigma$-stable closed subset of $G$ and hence stably compact, so Theorem~\ref{thm.Schauder-Tychonoff} applies.
Finally, taking $K=\{\theta\}$ in Theorem~\ref{thm.Fan-Iokhvidov} yields $x_0\in G$ with $T(x_0)\in x_0+\{\theta\}$, that is, $T(x_0)=x_0$.
 
\par
\noindent\textit{\textbf{Implication (ii).}}
Define $f:G\times G\rightarrow L^{0}(\mathcal{F})$ by 
$$f(x,y)=\|y-T(y)\|-\|x-T(y)\|,~\forall x,y\in G.$$
Then $f$ is $\sigma$-stable and $\mathcal{T}_{c}(\|\cdot\|)$-continuous, $f(\cdot,y)$ is $L^0$-concave for each $y\in G$, and $f(x,x)=0$. By Theorem~\ref{thm.Minimax inequality}, there exists
$y_0\in G$ such that 
$$\bigvee_{x\in G} f(x,y_0)\leq 0.$$
Taking $x=T(y_0)$ gives $\|y_0-T(y_0)\|\leq 0$, and hence $T(y_0)=y_0$.

\par
\noindent\textit{\textbf{Implication (iii).}}
For each $\varepsilon\in L^0_{++}(\mathcal F)$, define $A_\varepsilon,B_\varepsilon:G\to 2^G\setminus\{\emptyset\}$ by
$$\begin{aligned}
	A_\varepsilon(x)
	&=\{y\in G\mid \|y-x\|<\varepsilon~\text{on}~1\},\\
	B_\varepsilon(x)
	&=\{y\in G\mid \|y-T(x)\|<\varepsilon~\text{on}~1\}.
\end{aligned}$$
Then $A_\varepsilon$ and $B_\varepsilon$ satisfy the hypotheses of Theorem~\ref{thm.coincidence}. Hence, there exists $x_\varepsilon\in G$ such that
$$\|x_\varepsilon-T(x_\varepsilon)\|<2\varepsilon~\text{on}~1.$$
By part (3) of Lemma \ref{lemm.tools}, the $\mathcal{T}_{c}(\|\cdot\|)$-continuous $\sigma$-stable function $x\mapsto\|x-T(x)\|$ attains its minimum on $G$, namely, there exists $x_{0}\in G$ such that
$$\|x_0-T(x_0)\|<2\varepsilon~\text{on}~1$$   
for any $\varepsilon\in L^0_{++}(\mathcal F)$. Thus, $T(x_{0})=x_{0}$. 

\par\medskip
\noindent\textit{\textbf{Implication (iv).}} Let $X_1=X_2=G$ and $X=X_1\times X_2$. Define $f_1,f_2:X\to L^0(\mathcal{F})$ by
$$f_1(x,y)=-\|x-y\|,~\forall (x,y)\in X,$$
$$f_2(x,y)=-\|y-T(x)\|,~\forall (x,y)\in X.$$
Then $f_1$ and $f_2$ are $\sigma$-stable and continuous. For each $(x,y)\in X$, $f_1(\cdot,y)$ and $f_2(x,\cdot)$ are $L^0$-concave and hence $L^0$-quasiconcave.
By Theorem~\ref{thm.Nash equilibrium}, there exists $(x_0,y_0)\in X$ such that
$$
\begin{aligned}
	0\geq f_1(x_0,y_0)
	&=\bigvee\{f_1(x,y_0)\mid x\in G\}
	\geq f_1(y_0,y_0)=0,\\
	0\geq f_2(x_0,y_0)
	&=\bigvee\{f_2(x_0,y)\mid y\in G\}
	\geq f_2(x_0,T(x_0))=0.
\end{aligned}
$$
Thus $x_0=y_0=T(x_0)$.

\par\medskip
Together with the forward implications, these return arguments establish the equivalence of the principles lying on the cycles in the diagram. The random minimax principle(Theorem~\ref{thm.Minimax principle}) and the random Hartman--Stampacchia theorem (Theorem~\ref{thm.Hartman-Stampacchia}) are included as further consequences of the same framework. Their reverse implications are not established by the preceding return arguments.


\end{document}